\documentclass{article}    

\usepackage{latexsym}
\usepackage{a4wide}
\usepackage{amscd}
\usepackage{graphics}
\usepackage{amsmath}
\usepackage{amssymb}
\input xy
\xyoption{all}

\newcommand{\N}{\mathbb{N}}                     % the natural numbers
\newcommand{\Z}{\mathbb{Z}}                     % the integer numbers
\newcommand{\R}{\mathbb{R}}                     % the real line
\newcommand{\C}{\mathbb{C}}                     % the complex plane
\newcommand{\T}{\mathbb{T}}                     % the torus
\newcommand{\set}[2]{\left\{{#1}\mid{#2}\right\}}       % the set
\newcommand{\dist}{\mathrm{dist\,}}             % distance
\newcommand{\coker}{\mathrm{coker\,}}           % Cokernel
\newcommand{\Span}{\mathrm{span\,}}             % span
\newcommand{\ind}{\mathrm{ind\,}}               % Fredholm index
\newcommand{\codim}{\mathrm{codim}}           % codimension
\newcommand{\supp}{\mathrm{supp\,}}             % support
\newcommand{\ran}{\mathrm{ran\,}}   % range
\newcommand{\crit}{\mathrm{crit}}

\newcommand{\Det}{\mathrm{Det}}

\newtheorem{thm}{\sc Theorem}[section]      % numbered within each section

\newtheorem{lem}[thm]{\sc Lemma}            % numbered along with Theorem
\newtheorem{prop}[thm]{\sc Proposition}     % numbered along with Theorem
\newtheorem{defn}{\sc Definition}[section]  % numbered within each section

\newtheorem{rem}{\sc Remark}[section]

\title{On the Floer homology of cotangent bundles}

\author{Alberto Abbondandolo\thanks{Scuola Normale Superiore, Piazza 
dei Cavalieri 7,
56125 Pisa, Italy.}\hspace{32pt} Matthias
Schwarz\thanks{Mathematisches Institut, Universit\"at Leipzig,
04109 Leipzig, Germany. The second author was partially supported
by DFG grant no. SCHW892/2-1.}}

\date{February 24, 2004}

\begin{document}

\maketitle

\begin{abstract}
This paper concerns Floer homology for periodic orbits and for a
Lagrangian intersection problem on the cotangent bundle $T^*M$ of
a compact orientable manifold $M$. The first result is a new
$L^{\infty}$ estimate for the solutions of the Floer equation,
which allows to deal with a larger - and more natural - class of
Hamiltonians. The second and main result is a new construction of
the isomorphism between the Floer homology and the singular
homology of the free loop space of $M$, in the periodic case, or
of the based loop space of $M$, in the Lagrangian intersection
problem. The idea for the construction of such an isomorphism is
to consider a Hamiltonian which is the Legendre transform of a
Lagrangian on $TM$, and to construct an isomorphism between the
Floer complex and the Morse complex of the classical Lagrangian
action functional on the space of $W^{1,2}$ free or based loops on
$M$.
\end{abstract}

\section*{Introduction}

Let $M$ be a compact orientable smooth manifold without boundary.
Points in its cotangent bundle $T^*M$ will be denoted by $(q,p)$,
with $q\in M$ and $p\in T_q^*M$. The manifold $T^* M$ has a
canonical 1-form, the Liouville form $\theta=pdq$, whose
differential is the standard symplectic form $\omega= dp \wedge
dq$. A time-dependent and 1-periodic Hamiltonian function $H:\T
\times T^* M \rightarrow \R$, $\T=\R/\Z$, produces the
time-dependent Hamiltonian vector field $X_H$ on $T^* M$, defined
by $\omega(X_H,\cdot)=-dH$, whose 1-periodic orbits are the
critical points of the Hamiltonian action functional
\[
\mathcal{A}(x) = \int_{\T} x^* (\theta) - \int_0^1 H(t,x(t))\, dt
\]
on the space of smooth loops on $T^* M$. We shall denote the set
of all 1-periodic orbits of $X_H$ by $\mathcal{P}(H)$, and we shall
assume that every 1-periodic orbit is non-degenerate (a
condition which holds for a generic choice of $H$, in several
reasonable senses).

Let us assume that $H$ behaves as a positive quadratic form in the $p$
variables for $|p|$ large. More precisely, let us assume the following 
conditions:
\begin{description}
\item[(H1)] $dH(t,q,p)\left[p \frac{\partial}{\partial p}\right] -
H(t,q,p)\geq h_0 |p|^2- h_1$, for some constants $h_0>0$ and $h_1\geq 0$;
\item[(H2)] $|\nabla_q H(t,q,p)| \leq h_2 (1+|p|^2)$, $|\nabla_p
  H(t,q,p) |\leq h_2 (1+|p|)$, for some constant $h_2\geq 0$.
\end{description}
These assumptions are stated in terms of some metric on $M$, but
they are actually independent of the choice of the metric, in the
sense that if the metric on $M$ is changed then (H1) and (H2) will
still hold, with different constants $h_0,h_1,h_2$. 
Notice
also that these conditions involve only the behavior of $H$ for
$|p|$ large, and they impose no restriction on the behavior of $H$
on compact subsets of $T^*M$. Physical Hamiltonians of the form
\[
H(t,q,p) = \frac{1}{2} |T(t,q)p-A(t,q)|^2 + V(t,q),
\]
where the symmetric tensor $T^*T$ is everywhere positive, satisfy
(H1) and (H2). Let us endow $T^*M$ with a time-dependent
1-periodic $\omega$-compatible almost complex structure $J$,
assumed to be close enough to the almost complex structure induced
by some metric on $M$.

Under these assumptions, the elements of $\mathcal{P}(H)$ generate
a chain complex of Abelian groups, the Floer complex
$\{CF_*(H),\partial_*(H,J)\}$. Indeed, (H1) and (H2) easily imply
that the set of $x\in \mathcal{P}(H)$ with $\mathcal{A}(x)\leq a$
is finite for every real $a$, and we will prove that the set of
solutions $u: \R \times \T \rightarrow T^*M$ of the Floer equation
\begin{equation}
\label{floer}
\partial_s u -  J(t,u) (\partial_t u - X_H(t,u)) = 0 \quad
\forall (s,t)\in \R \times \T,
\end{equation}
with action bound $|\mathcal{A}(u(s,\cdot))|\leq a$ for every $s\in \R$,
is bounded in $L^{\infty}$. The latter statement improves the
$L^{\infty}$ estimates proved by Cieliebak \cite{cie94} under more
restrictive - and not metric independent - assumptions on the
Hamiltonian. Essentially, the reason why we get better estimates is
that we deal directly with the Cauchy-Riemann operator, instead than
differentiating and using the maximum principle for the Laplace
operator.

Then the methods of standard Floer theory for
compact symplectic manifolds can be applied. Actually, the
situation is somehow simpler than the general compact case,
because the presence of the Lagrangian foliation given by the fibers of $T^*M$
implies that the Conley-Zehnder index $\mu_{CZ}(x)$ of every $x\in
\mathcal{P}(H)$ is a well-defined integer, and because the
exactness of $\omega$ excludes the possible lack of compactness of
solutions of (\ref{floer}) coming from the phenomenon of bubbling
of $J$-holomorphic spheres. Therefore, the equation (\ref{floer})
can be successfully seen as the negative gradient equation for the
action functional $\mathcal{A}$, and the analogue of Morse theory
for such a functional can be developed: if $CF_k(H)$ denotes the
free Abelian group generated by the elements of $\mathcal{P}(H)$
of Conley-Zehnder index $k$, one constructs a boundary operator
\[
\partial_k(H,J) : CF_k(H) \rightarrow CF_{k-1}(H)
\]
by counting the number of solutions of (\ref{floer}) which connect
elements of $\mathcal{P}(H)$ with Conley-Zehnder index $k$ to
those with index $k-1$. The isomorphism class of the chain complex
$\{CF_*(H),\partial_*(H,J)\}$ - called the Floer complex of
$(T^*M,H,J)$ - does not depend on the choice of the almost complex
structure $J$. The homology of the Floer complex
does not depend on the Hamiltonian $H$, as long as $H$ satisfies (H1)
and (H2) (see \cite{oh97} for other possible choices, such as the
case of compactly supported Hamiltonians). Thus it makes sense to
talk about the Floer homology of $T^* M$, which we denote by
$HF_*(T^*M)$.

Unlike the compact case where the homology of the Floer complex is
just the singular homology of the underlying symplectic manifold,
the homology of the Floer homology of $T^*M$ can be fairly
complicated. Indeed, Viterbo \cite{vit95} has shown that it is
isomorphic to the singular homology of $\Lambda(M)$, the free loop
space of $M$. His proof makes use of generating functions and of
Traynor's homology for generating functions \cite{tra94}. A
complete proof is contained in \cite{vit96}. See also \cite{mo97},
and \cite{mo98}. Actually, in Viterbo's work cohomology is
considered, and the use of $L^{\infty}$ estimates
for wide classes of Hamiltonians is avoided, by
considering the Floer cohomology of the compact symplectic
manifold with convex contact-type boundary $\set{(q,p)\in T^*M}{|p|\leq
r}$, and then taking a limit of the Floer groups for $r\rightarrow
+\infty$.

Another beautiful approach to the question of the isomorphism
between the Floer homology of $T^*M$ and the singular homology of
$\Lambda(M)$, has been recently proposed by Salamon and Weber
\cite{sw03}. Their idea consists in considering a Hamiltonian of
the form
\[
H(t,q,p) = \frac{1}{2} |p|^2 + V(t,q),
\]
and almost complex structures which in the horizontal-vertical
splitting of $TT^*M$ given by some metric have the form
\[
J_{\epsilon} = \left( \begin{array}{cc} 0 & \epsilon \\
-1/\epsilon & 0 \end{array} \right).
\]
Writing $u=(q,p)$, and making the change of variable $s\mapsto
s/\epsilon$, the Floer equation (\ref{floer}) becomes
\begin{equation}
\label{floer2}
\begin{cases} \partial_s q - \nabla_t p - \nabla V(t,q) = 0, & \\
\epsilon^2 \nabla_s p + \partial_t q - p = 0. & \end{cases}
\end{equation}
As $\epsilon$ tends to 0 one obtains, at least formally, that
$p=\partial_t q$ and $q$ solves the heat equation
\begin{equation}
\label{heat}
\partial_s q - \nabla_t \partial_t q - \nabla V(t,q) = 0,
\end{equation}
where $s$ plays the role of time, and $t$ of space. This is a
better equation than the Floer equation, because the corresponding
Cauchy problem is well-posed in suitable spaces, and one can
derive the equivalent of Morse theory for its flow, obtaining a
chain complex whose homology can be proved to be isomorphic to the
singular homology of $\Lambda(M)$. A rigorous formulation of the
asymptotics for $\epsilon$ going to 0 then shows that for any
$a\in \R$ there is some small $\epsilon$ for which the Floer
subcomplex given by solutions of action not exceeding $a$ is
isomorphic to the corresponding subcomplex of the heat flow
equation, just because there is a one-to-one correspondence
between the relevant solutions of the corresponding PDE's. A limit
for $a\rightarrow +\infty$ then allows to conclude.

\medskip

The aim of this paper is to present a third construction of the
isomorphism between the Floer homology of $T^* M$ and the singular
homology of $\Lambda(M)$. Taking advantage of the freedom in the
choice of the Hamiltonian, provided that it satisfies (H1) and
(H2), we choose $H$ to be the Legendre transform of a Lagrangian
function $L:\T \times TM \rightarrow \R$, which is assumed to be
strongly convex on the fibers $T_q M$:
\begin{description}
\item[(L1)] $\nabla_{vv} L(t,q,v) \geq \ell_0 I$, for some constant $\ell_0>0$,
\end{description}
and to have bounds on the second derivatives, analogous to (H2):
\begin{description}
\item[(L2)] $|\nabla_{qq} L(t,q,v)| \leq \ell_1 (1+|v|^2)$, $|\nabla_{qv}
  L(t,q,v)| \leq \ell_1 (1+|v|)$, and 
$|\nabla_{vv} L(t,q,v)| \leq
  \ell_1$, for some constant $\ell_1\geq 0$.
\end{description}
The Hamiltonian function takes the form
\[
H(t,q,p) = \max_{v\in T_q M} ( p[v] - L(t,q,v)),
\]
and for each $(t,q,p)$ the above maximum is achieved at a unique
point $v(t,q,p)$. The Legendre transform $(t,q,p)\mapsto
(t,q,v(t,q,p))$ establishes a one-to-one correspondence
$(q,p)\mapsto (q,\dot{q})$ between the solutions of the first
order Hamiltonian system on $T^* M$ given by $H$, and the
solutions of the second order Lagrangian system on $M$ given by
$L$, which can be written in local coordinates as
\[
\frac{d}{dt} \partial_v L(t,q(t),\dot{q}(t)) = \partial_q
L(t,q(t),\dot{q}(t)).
\]
In the latter formulation,
the set of 1-periodic solutions - denoted by $\mathcal{P}(L)$ - is
the set of critical points of the Lagrangian action functional
\[
\mathcal{E}(q) = \int_0^1 L(t,q(t),\dot{q}(t))\, dt
\]
on the space of smooth loops on $M$. Developing Morse theory for
$\mathcal{E}$ is considerably simpler than developing it for
$\mathcal{A}$. Indeed, $\mathcal{E}$ is smooth on the Hilbert
manifold $\Lambda^1(M)$, the space of loops on $M$ of Sobolev
class $W^{1,2}$, satisfies the Palais-Smale condition (as proved
by Benci \cite{ben86}), and is bounded below. Moreover, classical
results by Duistermaat \cite{dui76} show that the Morse index
$m(q)$ of every $q\in \mathcal{P}(L)$ coincides with the
Conley-Zehnder index $\mu_{CZ}(x)$ of the corresponding solution
$x\in \mathcal{P}(H)$ (which is therefore always non-negative in
the case of Hamiltonians which are strictly convex in the $p$
variables). The finiteness of the indices, the Palais-Smale
condition, and the lower bound on $\mathcal{E}$ make it possible
to apply infinite dimensional Morse theory as developed by Palais
\cite{pal63}. Actually, it is convenient to use also here the
Morse complex approach: if we denote by $CM_k(\mathcal{E})$ the
free Abelian group generated by the critical points of
$\mathcal{E}$ of Morse index $k$, we obtain a boundary
homomorphism
\[
\partial_k (\mathcal{E},g) : CM_k(\mathcal{E}) \rightarrow
CM_{k-1} (\mathcal{E}),
\]
by introducing a Morse-Smale Riemannian metric $g$ on
$\Lambda^1(M)$ and by counting the solutions of the corresponding
negative gradient equation
\begin{equation}
\label{morse}
\gamma^{\prime} = - \nabla \mathcal{E} (\gamma),
\end{equation}
connecting the critical points of Morse index $k$ to those of
index $k-1$. The fact that this homomorphism is a boundary
operator, and the fact that the homology of the chain complex
$\{CM_*(\mathcal{E}) ,\partial_*(\mathcal{E},g)\}$ - called the
Morse complex of $( \mathcal{E},g)$ - is isomorphic to the
singular homology of $\Lambda^1(M)$, just comes from the fact that
$\{CM_*(\mathcal{E}) ,\partial_*(\mathcal{E},g)\}$ turns out to be
the chain complex associated to a suitable cellular filtration of
$\Lambda^1(M)$ (see \cite{ama04m} for a complete exposition).
Notice also that the inclusion $\Lambda^1(M)\hookrightarrow
\Lambda(M)$ is a homotopy equivalence, hence the two spaces have
isomorphic singular homology.

The isomorphism between the Floer homology of $T^*M$ and the
singular homology of $\Lambda(M)$ is then an immediate consequence
of the following stronger result, which is the main theorem of
this paper (see Theorem \ref{prec} for a more complete and precise
statement):

\medskip 

\noindent {\bf Theorem.} {\em
There exists a chain-complex isomorphism}
\[
\Theta  : \{CM_*(\mathcal{E}),\partial_*(\mathcal{E},g)\}
\longrightarrow \{CF_*(H), \partial_*(H,J) \}.
\]

\medskip

Let us sketch the idea of the proof. We already know that there is
a one-to-one correspondence between the generators of
$CM_k(\mathcal{E})$ and those of $CF_k(H)$. However this
correspondence need not produce a chain homomorphism. In order to
produce such a chain homomorphism, we shall consider the set of
solutions $\gamma(s)$, $s\in ]-\infty,0]$, which flow off some
$q\in \mathcal{P}(L)$ by the negative gradient flow equation
(\ref{morse}), and for $s=0$ can be lifted to a loop in $T^* M$
which is the trace at $s=0$ of a solution $u(s,t)$, $(s,t)\in
[0,+\infty[ \times \T$ of the Floer equation (\ref{floer}),
converging to some $x\in \mathcal{P}(H)$ for $s\rightarrow
+\infty$. In other words, for $q\in \mathcal{P}(L)$ and $x\in
\mathcal{P}(H)$, we shall consider the set
\[
\mathcal{M}^+(q,x) =  \{u:[0,+\infty]\times \T \rightarrow T^*M\, | \,
 u \mbox{ solves (\ref{floer}), } u(+\infty,\cdot) = x, 
\mbox{ and } \tau^* u(0,\cdot) \in W^u(q)\}.  
\]
Here $\tau^*:T^*M \rightarrow M$ denotes the standard projection, and
$W^u(q)\subset \Lambda^1(M)$ is the unstable manifold of $q$, that is
the set of $w\in \Lambda^1(M)$ such that the solution $\gamma(s)$ of
(\ref{morse}) with $\gamma(0)=w$ converges to $q$ for $s\rightarrow
-\infty$. 
The above problem is a Fredholm one, because the unstable manifold of
$q$ is finite dimensional, and because we are imposing Lagrangian
boundary conditions and non-degenerate asymptotic conditions on
the Cauchy-Riemann type equation (\ref{floer}) on $[0,+\infty[ \times
\T$. Indeed, we shall
prove that for a generic choice of $J$,
$\mathcal{M}^+(q,x)$ is a smooth manifold of dimension
$m(q)-\mu_{CZ}(x)$. In particular when $m(q)=\mu_{CZ}(x)=k$,
$\mathcal{M}^+(q,x)$ is a discrete set, and we shall define the
homomorphism
\[
\Theta_k : CM_k(\mathcal{E}) \rightarrow CF_k(H),
\]
by counting the elements of such sets. A crucial ingredient is of
course the compactness of $\mathcal{M}^+(q,x)$, which will be a
consequence of the following simple but important action estimate: if
$x(t)=(q(t),p(t))$ is a loop on $T^* M$, then
\begin{equation}
\label{idis} \mathcal{A}(x) \leq \mathcal{E}(q),
\end{equation}
the equality holding if and only if $p$ is related to $\dot{q}$ by
the Legendre transform. In particular,
$\mathcal{A}(x)=\mathcal{E}(q)$ holds for $x$ and $q$
corresponding to the same solution. 
The above estimate guarantees
that if $u\in \mathcal{M}^+(q,x)$ then $\mathcal{A}(u(s,\cdot)) \leq
\mathcal{E}(q)$ for every $s\geq 0$, the starting point to get compactness.

A standard gluing argument shows that $\Theta$ is a chain
homomorphism. The fact that $\Theta$ is an isomorphism is again a
consequence of (\ref{idis}), together with its differential version
\begin{equation}
\label{dis2} d^2 \mathcal{A}(x) [\zeta,\zeta] \leq d^2 \mathcal{E}(q)
[D\tau^*(x)\zeta,D\tau^*(x) \zeta]
\quad \forall \zeta\in C^{\infty}(x^*(TT^*M)),
\end{equation}
for every $x\in \mathcal{P}(H)$ and $q\in\mathcal{P}(L)$
corresponding to the same periodic solution. Indeed, (\ref{idis})
implies that $\mathcal{M}^+(q,x)=\emptyset$ if $\mathcal{E}(q)\leq
\mathcal{A}(x)$, unless $q$ and $x$ correspond to the same
periodic solution, in which case $\mathcal{M}^+(q,x)$ consists
just of the stationary solution $u(s,t)=x(t)$. Such a
solution is a regular one (in the sense of transversality theory)
because of (\ref{dis2}). We conclude that, if we order the
generators of $CM_k(\mathcal{E})$ and $CF_k(H)$ by increasing
action, the homomorphism $\Theta_k$ is given by a square matrix
which is lower triangular and has $\pm 1$ on each diagonal entry.
Hence $\Theta$ is an isomorphism.

\medskip

Notice that we cannot construct an isomorphism $CF_*(H)
\rightarrow CM_*(\mathcal{E})$ using moduli spaces of
solutions analogue to $\mathcal{M}^+(q,x)$. Indeed, we would not
obtain a Fredholm problem, and we would not dispose of action
estimates guaranteeing the compactness property.

\medskip

In the proof sketched above, we are coupling Floer theory for
$\mathcal{A}$ with Morse theory for $\mathcal{E}$ on
$\Lambda^1(M)$. We should however stress the fact that it
is not essential that on the Morse side we have a $W^{1,2}$
negative gradient flow: another Morse-Smale flow, having
$\mathcal{E}$ as a Lyapunov function would suffice. The $W^{1,2}$
option is the obvious one, and it makes particularly easy the step
from the Morse complex to the singular homology of the loop space,
but one could try other possibilities. For instance, one may
couple Floer theory with the heat flow equation (\ref{heat}),
avoiding all the analysis on the behavior of the solutions of the
Floer equation (\ref{floer2}) as $\epsilon$ tends to 0.

\medskip

A similar construction works for the fixed ends case: given
$q_0,q_1\in M$, we look at solutions $q$ of the Lagrangian
system such that $q(0)=q_0$, $q(1)=q_1$. On the Hamiltonian side,
this means that we are looking at solutions $x$ such that
$x(0)\in T_{q_0}^*M$ and $x(1) \in T_{q_1}^*M$. The Floer equation
is the same as (\ref{floer}), but this time $u$ is defined on the
strip $\R \times [0,1]$, with boundary conditions $u(s,0)\in
T_{q_0}^* M$, $u(s,1)\in T_{q_1}^*M$, for every $s\in \R$. Again,
one finds a Floer complex and a Morse complex, the latter one
being associated to the Lagrangian action functional $\mathcal{E}$
on $\Omega^1(M,q_0,q_1)$, the Hilbert manifold of $W^{1,2}$ curves
$[0,1]\rightarrow M$ connecting $q_0$ to $q_1$. Such a manifold
is homotopically equivalent to $\Omega(M)$, the based loop space
of $M$, hence we obtain that the Floer homology of the fixed ends
problem is isomorphic to the singular homology of $\Omega(M)$.

\medskip

We wish to emphasize the fact that the assumption that the
manifold $M$ should be orientable is made only to have some
technical simplifications - mainly in the choice of suitable
preferred symplectic trivializations of the tangent bundle of $T^*
M$ along the solutions of the Hamiltonian system - but it could be
easily dropped.

\medskip

In a forthcoming paper we will prove that the isomorphism defined
here from the Floer homology of $T^*M$ - in  the periodic cases -
to the singular homology of the free loop space of $M$ is actually
a ring isomorphism: it relates the pair-of-pants product in Floer
homology to the Chas-Sullivan loop product on the singular
homology of $\Lambda(M)$. In the fixed ends problem, this
isomorphism relates the $Y$ product (the analogue of the
pair-of-pants product with strips) in Floer homology to the
classical Pontrjagin product on the singular homology of
$\Omega(M)$. Related results are proven by Ralph Cohen \cite{coh04}
and Antonio Ramirez \cite{ram04}. 

\paragraph{Acknowledgements.} We wish to thank the
Max Planck Institute for Mathematics in the Sciences of Leipzig
and the Courant Institute of New York, and in particular Tobias
Colding and Helmut Hofer, for their kind hospitality.  We are also
indebted to Ralph Cohen, Yong-Geun Oh, Dietmar Salamon, 
Claude Viterbo, and Joa Weber for many fruitful discussions.

\tableofcontents

\section{The Floer complex of the Hamiltonian action functional}

\subsection{Hamiltonian dynamical systems on cotangent bundles}

Let $M$ be a connected compact orientable smooth manifold of dimension $n$.
Points in the cotangent bundle $T^*M$ will be denoted
by $(q,p)$, with $q\in M$, $p\in T_q^* M$, and
$\tau^*:T^*M \rightarrow M$ will denote the standard projection.
The cotangent bundle $T^*M$ carries the following canonical
structures: the Liouville 1-form $\theta$ and the Liouville vector
field $\eta$, which are defined by
\[
\theta(\zeta) = p ( D\tau^*(x) [\zeta]) = d\theta (\eta,\zeta) \quad \forall
\zeta\in T_x T^*M, \; x=(q,p)\in T^*M,
\]
and the symplectic structure $\omega=d\theta$. In local coordinates
$(q,p)$ of $T^* M$ we have
\[
\theta_{\mathrm{local}} = p\, dq, \quad \eta_{\mathrm{local}} = p
\frac{\partial}{\partial p}, \quad \omega_{\mathrm{local}} = dp \wedge
dq.
\]
The vertical space
\[
T_x^v T^*M = \ker D \tau^*(x) \cong T_q^* M, \quad x=(q,p)\in T^* M,
\]
is a Lagrangian subspace of $(T_x T^*M,\omega_x)$.

A 1-periodic Hamiltonian $H$, i.e.\ a smooth function $H : \T \times
T^*M \rightarrow \R$, $\T = \R/\Z$, determines a 1-periodic vector
field, the Hamiltonian vector field $X_H$ defined by
\[
\omega( X_H(t,x),\zeta) = -dH(t,x) [\zeta], \quad \forall \zeta\in T_x T^*M.
\]
In local coordinates, the Hamiltonian equation
\begin{equation}
\label{ham}
\dot{x}(t) = X_H(t,x(t)),
\end{equation}
takes the classical physical form
\[
\left\{ \begin{array}{l} \dot{q} = \partial_p H(t,q,p), \\ \dot{p} = -
  \partial_q H(t,q,p). \end{array} \right.
\]
The integral flow of the vector field $X_H$ will be denoted by
$\phi_H^t$. We will be interested in the set
$\mathcal{P}_{\Lambda}(H)$ of 1-periodic solutions of (\ref{ham}), and
  in the set $\mathcal{P}_{\Omega}(H)$ of solutions
  $x:[0,1]\rightarrow T^* M$ of (\ref{ham}) such that $x(0)\in
  T_{q_0}^* M$ and $x(1)\in T_{q_1}^*M$, for two fixed
  points\footnote{In this latter case, it is not necessary to assume
    $H$ to be 1-periodic in time and $M$ to be compact.
    However we shall keep these
    assumptions in order to have a uniform presentation.}
$q_0,q_1\in M$. In each of these cases\footnote{The symbols
    $\Lambda$ and $\Omega$ will appear as subscripts of many objects
    we are going to introduce. We will omit such a subscript whenever
    we wish to consider both situations at the same time.} we shall
  make one of the following non-degeneracy assumptions:
\begin{description}

\item[(H0)$_{\mathbf{\Lambda}}$] every solution $x\in
  \mathcal{P}_{\Lambda} (H)$ is non-degenerate, meaning that 1 is not
  an eigenvalue of $D\phi_H^1 (x(0))\in GL(T_{x(0)} T^*M)$;

 \item[(H0)$_{\mathbf{\Omega}}$] every solution $x\in
  \mathcal{P}_{\Omega} (H)$ is non-degenerate, meaning that the image
  of $T_{x(0)}^v T^* M$ by $D \phi_H^1(x(0))$ has intersection $(0)$
  with $T_{x(1)}^v T^*M$.

\end{description}

\noindent The above conditions imply that the set
$\set{x(0)}{x\in \mathcal{P}(H)}$ is discrete in $T^*M$.

\begin{rem}
Assumption (H0) holds for a generic choice of $H$, in several
reasonable senses (see for instance \cite{sz92,web02}). It is
worth remarking that if $H$ satisfies the additional condition
\[
\partial_{pp} H(x)\in GL(T_x^v T^*M) \quad \forall x\in T^*M,
\]
then $H+V$ satisfies (H0) for a residual set of potentials $V\in
C^{\infty}( \T \times M, \R)$. If moreover $\partial_{pp} H>0$, then
$H$ satisfies (H0)$_{\Omega}$ for a set of $(q_0,q_1)\in M \times M$
having full measure.
\end{rem}

\subsection{The Maslov index}

Let $\R^{2n} = \R^n \times \R^n$ be endowed with its
standard Euclidean product, with its standard symplectic structure
\[
\omega_0 = dp \wedge dq, \quad (q,p) \in \R^n \times \R^n,
\]
and with its standard complex structure
\[
J_0 = \left( \begin{array}{cc} 0 & I \\ -I & 0 \end{array} \right),
\]
so that $\omega_0(\zeta_1,\zeta_2) = J_0 \zeta_1 \cdot \zeta_2$. We denote by
$Sp(2n)$ the group of symplectic automorphisms of
$(\R^{2n},\omega_0)$, by $\mathcal{L}(n)$ the space of Lagrangian
subspaces of $(\R^{2n},\omega_0)$, and by $\lambda_0$ the vertical
Lagrangian subspace $\lambda_0= (0) \times \R^n$.

We recall that the {\em Conley-Zehnder index} assigns an integer
$\mu_{CZ}(\gamma)$ to every path of symplectic automorphisms $\gamma$
belonging to the space
\begin{equation}
\label{nd}
\set{\gamma \in C^0([0,1],Sp(2n)) }{\gamma(0)=I \mbox{ and 1 is not an
    eigenvalue of } \gamma(1)}.
\end{equation}
See \cite{sz92}, section 3. For future reference, we recall that the
Conley-Zehnder indices of the paths
$\gamma_1,\gamma_2:[0,1]\rightarrow Sp(2)$, $\gamma_1(t) =
\left( \begin{array}{cc} e^{\alpha t} & 0 \\ 0 & e^{-\alpha t}
\end{array} \right)$, with $\alpha\in \R \setminus \{0\}$,
$\gamma_2(t) =
e^{t\theta J_0}$, with $\theta\in \R\setminus 2\pi\Z$, are the integers
\begin{equation}
\label{czi}
\mu_{CZ}(\gamma_1) = 0, \quad
\mu_{CZ}(\gamma_2) = 2 \left\lfloor \frac{\theta}{2\pi} \right\rfloor +
1.
\end{equation}
A related notion is the {\em relative
  Maslov index} of a pair of Lagrangian paths, which assigns a
  half-integer $\mu(\lambda_1,\lambda_2)$ to every pair of continuous
  paths $\lambda_1,\lambda_2 : [0,1]\rightarrow \mathcal{L}(n)$. See
  \cite{rs93}. For future reference, we recall that if $n=1$ and $\gamma(t) =
  e^{t\theta J_0}$, with $\theta\in \R\setminus \pi \Z$, there holds
\begin{equation}
\label{rmi}
\mu(\gamma \lambda_0,\lambda_0) = \frac{1}{2} + \left\lfloor
\frac{\theta}{\pi} \right\rfloor.
\end{equation}

\begin{lem} \label{triv0}
(i) Assume that $M$ is orientable, and let $x\in
\mathcal{P}_{\Lambda}(H)$. Then the symplectic vector bundle
$x^*(TT^*M)$ admits a symplectic trivialization
\[
\Phi: \T \times \R^{2n} \rightarrow x^*(TT^*M)
\]
such that
\begin{equation}
\label{triv} \Phi(t)\lambda_ 0 = T_{x(t)}^v T^*M \quad \forall
t\in\T.
\end{equation}
(ii) Let $x\in \mathcal{P}_{\Omega}(H)$. Then the symplectic
vector bundle $x^*(TT^*M)$ admits a symplectic trivialization
\[
\Phi: [0,1] \times \R^{2n} \rightarrow x^*(TT^*M)
\]
such that
\begin{equation}
\label{triv2} \Phi(t)\lambda_ 0 = T_{x(t)}^v T^*M \quad \forall
t\in [0,1].
\end{equation}
\end{lem}

\begin{proof} (i) Since $M$ is orientable, the vector bundle $x^*(T^v
T^*M) \cong (\tau^* \circ x)^* (T^* M)$ is orientable, hence
trivial. Let
\[
\Psi:\T \times \R^n \rightarrow x^*(T^v T^*M)
\]
be a trivialization, and let $J$ be a $\omega$-compatible complex
structure on $x^*(TT^*M)$ (meaning that $\omega(\cdot,J\cdot)$ is
an inner product on $x^*(TT^*M)$). Then
\[
T_{x(t)} T^* M = J(t) T_{x(t)}^v T^* M \oplus T_{x(t)}^v T^*M,
\]
and the trivialization
\[
\Phi:\T \times \R^n \oplus \R^n \rightarrow x^*(TT^*M), \quad
\Phi(t) = (-J(t)\Psi(t)J_0) \oplus \Psi(t),
\]
is symplectic (actually unitary) and maps $\lambda_0$ into the
vertical subbundle.

(ii) Starting from the fact that the vector bundle $x^*(T^v T^*M)$
is trivial because $[0,1]$ is contractible, the construction is
identical to the one shown in (i).
\end{proof}

Let $x\in \mathcal{P}_{\Lambda}(H)$. We can use the symplectic
trivialization $\Phi$ provided by the above lemma to transform the
differential of the Hamiltonian flow along $x$ into a path in
$Sp(2n)$,
\[
\gamma_{\Phi} (t) = \Phi(t)^{-1} \, D\phi_H^t(x(0)) \, \Phi(0),
\]
which belongs to the space (\ref{nd}), thanks to (H0)$_{\Lambda}$.

Similarly if $x\in \mathcal{P}_{\Omega}(H)$, the symplectic
trivialization $\Phi$ provided by the above lemma produces the
path in $\mathcal{L}(n)$,
\[
\lambda_{\Phi} (t) = \Phi(t)^{-1} \, D \phi_H^t(x(0))  [T_{x(0)}^v
T^*M]
\]
such that $\lambda_{\Phi}(0)=\lambda_0$ and $\lambda_{\Phi}(1)
\cap \lambda_0 = (0)$, thanks to (H0)$_{\Omega}$.

Denote by $Sp(2n,\lambda_0)$ the subgroup of the symplectic group
consisting of those automorphisms which preserve the vertical
Lagrangian subspace $\lambda_0$:
\[
Sp(2n,\lambda_0) := \set{A\in Sp(2n)}{A \lambda_0 = \lambda_0} 
= \set{ \left( \begin{array}{cc} A_1 & 0 \\ A_2 & A_3 \end{array}
\right) }{A_1^* A_3 = I, \; A_1^* A_2 = A_2^* A_1}.
\]
It is easily seen that $Sp(2n,\lambda_0)$ is continuously
retractable onto its closed subgroup
\[
Sp(2n,\lambda_0) \cap U(n) = \set{ \left( \begin{array}{cc} R & 0
\\ 0 & R^* \end{array} \right) }{R\in O(n)},
\]
on which the determinant map $\det : U(n) \rightarrow S^1$ takes
the values $\pm 1$. It follows that $Sp(2n,\lambda_0)$ and
$Sp(2n,\lambda_0)\cap U(n)$ have two connected components, and
that the inclusions
\[
Sp(2n,\lambda_0) \hookrightarrow Sp(2n), \quad Sp(2n,\lambda_0)
\cap U(n) \hookrightarrow U(n),
\]
induce the zero homomorphism between fundamental groups.

\begin{lem}
(i) If $x\in \mathcal{P}_{\Lambda}(H)$, the Conley-Zehnder
$\mu_{CZ}(\gamma_{\Phi})$ does not depend on the symplectic
trivialization $\Phi$ satisfying (\ref{triv}).

(ii) If $x\in \mathcal{P}_{\Omega}(H)$, the relative Maslov index
$\mu(\lambda_{\Phi},\lambda_0)$ does not depend on the symplectic
trivialization $\Phi$ satisfying (\ref{triv2}).
\end{lem}

\begin{proof}
Let $x\in \mathcal{P}_{\Lambda}(H)$ and let $\Phi,\Psi$ be two
symplectic trivializations satisfying (\ref{triv}). Then
\[
\gamma_{\Psi}(t) = \alpha(t) \gamma_{\Phi}(t) \alpha(0)^{-1},
\]
for some
\[
\alpha: \T \rightarrow Sp(2n,\lambda_0).
\]
Since the inclusion $Sp(2n,\lambda_0)\hookrightarrow Sp(2n)$
induces the zero homomorphism between fundamental groups,
$\alpha(0) \gamma_{\Phi} \alpha(0)^{-1}$ and $\gamma_{\Psi}$ are
homotopic by a homotopy which fixes the end-points. The homotopy
and the naturality property of the Conley-Zehnder index imply that
$\mu_{CZ}(\gamma_{\Phi}) = \mu_{CZ} (\gamma_{\Psi})$.

Now let $x\in \mathcal{P}_{\Omega}(H)$ and let $\Phi,\Psi$ be two
symplectic trivializations satisfying (\ref{triv2}). Then
\[
\lambda_{\Psi} (t) = \alpha(t) \lambda_{\Phi}(t),
\]
for some $\alpha:[0,1]\rightarrow Sp(2n,\lambda_0)$, and by the
naturality of the relative Maslov index,
\[
\mu(\lambda_{\Phi}, \lambda_0) = \mu (\alpha \lambda_{\Phi}, \alpha
  \lambda_0) = \mu(\lambda_{\Psi}, \lambda_0).
\]
\end{proof}

The above lemma allows us to give the following:

\begin{defn} The Maslov index of a periodic solution $x\in
  \mathcal{P}_{\Lambda}(H)$ is the integer $\mu_{\Lambda}(x) :=
  \mu_{CZ}(\gamma_{\Phi})$, where $\Phi$ is a symplectic trivialization
  of $x^*(TT^*M)$ satisfying (\ref{triv}).

\noindent The Maslov index of a solution $x\in \mathcal{P}_{\Omega}(H)$ is the
integer $\mu_{\Omega}(x) := \mu(\lambda_{\Phi},\lambda_0) - n/2$,
where $\Phi$ is a symplectic trivialization of $x^*(TT^*M)$ 
satisfying (\ref{triv2}).
\end{defn}

Indeed, since $\lambda_{\Phi}(0)\lambda_0=\lambda_0$ and
$\lambda_{\Phi}(1)\cap \lambda_0 = (0)$, the number $\mu
(\lambda_{\Phi},\lambda_0) - n/2$ is an integer (see \cite{rs93},
Corollary 4.12).

\begin{rem}
The sign of the symplectic form and the $n/2$-shift have been chosen
in order for the above Maslov indices to coincide with the Morse
indices of corresponding critical points of the Lagrangian action
functional (see Theorem \ref{compind}). Notice that in the case of a
contractible periodic orbit $x\in \mathcal{P}_{\Lambda}(H)$, if
$\overline{x} : D \rightarrow T^*M$ is an extension of $x$ to the disc $D$,
there are symplectic trivializations of $\overline{x}^*(TT^*M)$
satisfying (\ref{triv}), so $\mu_{\Lambda}(x)$ coincides with the
Conley-Zehnder index in standard Floer theory for contractible
periodic orbits.
\end{rem}

\begin{rem}
It is actually possible to drop the assumption on the
orientability of $M$ also in the case of periodic solutions. In
the non-orientable case indeed, one can still single out a special
class of symplectic trivializations for which the Conley-Zehnder
index coincides with the Morse index of the Lagrangian action
functional (see \cite{web02}). For sake of simplicity, in this
paper we deal only with the orientable case.
\end{rem}

\subsection{The $\mathbf{L^2}$-gradient of the Hamiltonian action functional}

Denote by $\Lambda^1(T^*M)$ the space of all loops $x:\T \rightarrow
T^*M$ of Sobolev class $W^{1,2}$, and denote by $\Omega^1(T^*M, q_0,q_1)$ the
space of all paths $x:[0,1]\rightarrow T^*M$ of Sobolev class
$W^{1,2}$ such that $x(0)\in T_{q_0}^*M$ and $x(1)\in
T_{q_1}^*M$. These spaces have canonical Hilbert manifold
structures. The action functional
\[
\mathcal{A} (x) = \mathcal{A}_H (x) := 
\int x^* (\theta - H dt) = \int_0^1 \left(
\theta(\dot{x}) - H(t,x) \right) \, dt
\]
is smooth on $\Lambda^1(T^*M)$ and on $\Omega^1(T^*M, q_0,q_1)$. The
differential of $\mathcal{A}$ on both manifolds takes the form
\begin{equation}
\label{diff}
d\mathcal{A}(x) [\zeta] = \int_0^1 \left( \omega(\zeta,\dot{x}) -
dH(t,x)[\zeta] \right)\, dt = \int_0^1 \omega( \zeta, \dot{x} - X_H(t,x))
\, dt,
\end{equation}
so the critical points of $\mathcal{A}|_{\Lambda^1}$ are the elements of
$\mathcal{P}_{\Lambda}(H)$, while the critical points of
$\mathcal{A}|_{\Omega^1}$ are the elements of
$\mathcal{P}_{\Omega}(H)$. However, variational methods using the
gradient flow of $\mathcal{A}$ with respect to some metric compatible
with the $W^{1,2}$ topology are not suitable for classifying the
critical points of $\mathcal{A}$, because the Morse index of every
critical point is infinite, and because of lack of compactness
(the Palais-Smale condition would not hold, but see \cite{hof85b} for
a possible approach in this direction). It
was Floer's idea to overcome these difficulties by studying the
$L^2$-gradient equation for $\mathcal{A}$. More precisely, let $J$ be
a smooth almost complex structure on $T^*M$, 1-periodic
in the time variable $t$, and compatible with $\omega$, meaning that
\[
\langle \zeta_1, \zeta_2 \rangle_{J_t} := \omega (\zeta_1, J(t,x) \zeta_2), \quad
\zeta_1,\zeta_2 \in T_x T^*M, \; x\in T^*M,
\]
is a loop of Riemannian metrics on $T^*M$. We can rewrite (\ref{diff})
as
\[
d\mathcal{A}(x) [\zeta] = \int_0^1 \langle \zeta, - J(t,x)(\dot{x} -
X_H(t,x) \rangle_{J_t} \, dt,
\]
and we will denote by $\nabla_J \mathcal{A}$ the gradient of
$\mathcal{A}$ with respect to the $L^2$ inner product given by the
periodic metric $\langle \cdot, \cdot \rangle_{J_t}$, namely
\[
\nabla_J \mathcal{A} (x) = - J(t,x) ( \dot{x} - X_H(t,x)).
\]
We will be interested in the {\em negative gradient equation}, that is
in the Cauchy-Riemann type PDE
\begin{equation}
\label{cr}
\partial_s u - J(t,u) (\partial_t u - X_H(t,u)) = 0,
\end{equation}
where
\begin{equation}
\label{bdry}
\begin{split}
u\in C^{\infty} (\R \times \T, T^*M) \mbox{ in the $\Lambda$ case}, \\
u\in C^{\infty} (\R \times [0,1],T^* M), \; u(\R \times \{0\}) \subset
T_{q_0}^* M, \;  u(\R \times \{1\}) \subset T_{q_1}^* M, \mbox{ in the
  $\Omega$ case}.
\end{split}
\end{equation}
Given $x^-,x^+\in \mathcal{P}(H)$, we will denote by
$\mathcal{M}(x^-,x^+) = \mathcal{M}(x^-,x^+;H,J)$ the set of all
solutions of (\ref{cr},\ref{bdry}) such that
\[
\lim_{s\rightarrow \pm \infty} u(s,t) = x^{\pm}(t) \mbox{ uniformly in
  $t$.}
\]
As usual, we shall add the subscript $\Lambda$ or $\Omega$ when we
wish to distinguish between the periodic and the fixed-ends problem.
The elements of $\mathcal{P}(H)$ are the stationary solutions of
(\ref{cr},\ref{bdry}), and $\mathcal{A}$ is strictly decreasing on all the
other solutions. So $\mathcal{M}(x,x)$ contains only the element $x$,
and
\begin{equation}
\label{decr}
\mathcal{M}(x^-,x^+) = \emptyset \quad \mbox{if } \mathcal{A}(x^-)\leq
\mathcal{A}(x^+) \mbox{ and } x^-\neq x^+.
\end{equation}
Clearly, $\mathcal{M}_{\Lambda}(x^-,x^+)\neq \emptyset$ implies that
the loops $x^-$ and $x^+$ are homotopic, which is equivalent to saying
that their projections onto $M$, $\tau^* \circ x^-$ and $\tau^* \circ
x^+$, are homotopic. Similarly, $\mathcal{M}_{\Omega}(x^-,x^+)\neq
\emptyset$  implies that
the paths $x^-$ and $x^+$ are homotopic within the space of paths
having end-points on $T_{q_0}^* M$ and on $T_{q_1}^* M$,
which is equivalent to saying
that their projections onto $M$ are homotopic with fixed end-points.

We conclude this section by describing the standard functional setting which
allows to see $\mathcal{M}(x^-,x^+)$ as the set of zeros of a smooth
section of a Banach bundle. Let us fix a number $r>2$, and recall that
the maps on two-dimensional domains of Sobolev class $W^{1,r}$ are
H\"older continuous.
Fix two solutions $x^-,x^+\in \mathcal{P}(H)$, which are homotopic in the sense
explained above.

In the $\Lambda$ case, we define $\mathcal{B}_{\Lambda}=
\mathcal{B}_{\Lambda}(x^-,x^+)$ as the set of all maps $u: \R \times
\T \rightarrow T^* M$ of Sobolev class $W^{1,r}_{\mathrm{loc}}$ such
that there is $s_0\geq 0$ for which
\begin{equation}
\label{asymp}
u(s,t) = \begin{cases} \exp_{x^-(t)} (\zeta^-(s,t)) & \forall s\leq
  -s_0, \\ \exp_{x^+(t)} (\zeta^+(s,t)) & \forall s\geq s_0, \end{cases}
\end{equation}
where $\zeta^-$ and $\zeta^+$ are $W^{1,r}$ sections of 
the bundles ${x^-}^*(TT^*M)
\rightarrow ]-\infty,-s_0[ \times \T$ and ${x^+}^*(TT^*M)
\rightarrow ]s_0,+\infty[ \times \T$, respectively. Here $\exp$ denotes the
  exponential map with respect to some metric on $T^*M$, and the space
  of $W^{1,r}$ sections is also defined in terms of this metric. Then
  $\mathcal{B}_{\Lambda}$ can be given the structure of a smooth Banach
  manifold, and the tangent space at $u\in \mathcal{B}_{\Lambda}$ is
  identified with the space of $W^{1,r}$ sections of $u^*(TT^*M)$.

Similarly, $\mathcal{B}_{\Omega}= \mathcal{B}_{\Omega} (x^-,x^+)$
will be the Banach manifold of all
maps $u:\R \times [0,1]\rightarrow T^* M$ which are of Sobolev class
$W^{1,r}$ on every compact subset of $\R \times [0,1]$, such that
$u(\R \times \{0\}) \subset T_{q_0}^* M$, $u(\R \times \{1\}) \subset
T_{q_1}^* M$, and such that there is $s_0\geq 0$ for which
(\ref{asymp}) holds, $\zeta^-$ and $\zeta^+$ being $W^{1,r}$ sections of
${x^-}^*(TT^*M) \rightarrow ]-\infty,-s_0[ \times [0,1]$ and ${x^+}^*(TT^*M)
\rightarrow ]s_0,+\infty[ \times [0,1]$, respectively.

Denote by $\mathcal{W}_{\Lambda}= \mathcal{W}_{\Lambda} (x^-,x^+)$
(resp.\ $\mathcal{W}_{\Omega}= \mathcal{W}_{\Omega} (x^-,x^+)$)
the Banach bundle over $\mathcal{B}_{\Lambda}$ (resp.\
$\mathcal{B}_{\Omega}$) whose fiber $\mathcal{W}_u$ at $u$ is the
space of $L^r$ sections of $u^*(TT^*M)$. Then
$\mathcal{M}(x^-,x^+)$ is the set of zeros of the smooth section
\[
\partial_{J,H} : \mathcal{B} \rightarrow \mathcal{W}, \quad u \mapsto
\partial_s u + \nabla_J \mathcal{A}(u) = \partial_s u - J(t,u)
(\partial_t u - X_H(t,u)).
\]
Indeed, standard elliptic regularity results imply that the zeros
of the above section are smooth maps. Moreover, the non-degeneracy
assumption (H0) guarantees that $u(s,t)\rightarrow x^{\pm}(t)$ and
$\partial_s u(s,t)\rightarrow 0$ for $s\rightarrow \pm \infty$
exponentially fast, uniformly in $t$ (see e.g.
\cite{sch95,sal99}).

Denote by $\overline{\R}$ the extended real line
$\R\cup\{-\infty,+\infty\}$, with the differentiable structure induced
by the bijection $[-\pi/2,\pi/2]\rightarrow \overline{\R}$, $s
\mapsto \tan s$ for $s\in ]-\pi/2,\pi/2[$, $\pm\pi/2
    \mapsto\pm \infty$. A map $u\in \mathcal{M}_{\Lambda}(x^-,x^+)$
    (resp.\ $u\in \mathcal{M}_{\Omega}(x^-,x^+)$)
extends to a
    smooth map on $\overline{\R} \times \T$ (resp.\ $\overline{\R}
    \times [0,1]$), which we shall also denote by $u$.

\begin{lem}
\label{ltriv} (i) Let $u\in \mathcal{M}_{\Lambda}(x^-,x^+)$, and
let
\[
\Phi^{\pm} : \T \times \R^{2n} \rightarrow {x^{\pm}}^*(TT^*M)
\]
be two unitary trivializations such that $\Phi^{\pm}(t)\lambda_0 =
T_{x^{\pm}(t)}^v T^*M$ for every $t\in \T$. Then there exists a
smooth unitary trivialization
\[
\Phi: \overline{\R} \times \T \times \R^{2n} \rightarrow
u^*(TT^*M)
\]
such that $\Phi(\pm \infty,t) = \Phi^{\pm}(t)$ for every $t\in
\T$.

(ii) Let $u\in \mathcal{M}_{\Omega}(x^-,x^+)$, and let
\[
\Phi^{\pm} : [0,1] \times \R^{2n} \rightarrow {x^{\pm}}^*(TT^*M)
\]
be two unitary trivializations such that $\Phi^{\pm}(t)\lambda_0 =
T_{x^{\pm}(t)}^v T^*M$ for every $t\in [0,1]$, and such that the
isomorphisms
\[
\Phi^{\pm}(t)|_{\lambda_0} : \lambda_0 \rightarrow
T_{x^{\pm}(t)}^v T^*M \cong T_{\tau^* \circ x^{\pm}(t)}^* M
\]
are orientation preserving. Then there exists a smooth unitary
trivialization
\[
\Phi: \overline{\R} \times [0,1] \times \R^{2n} \rightarrow
u^*(TT^*M)
\]
such that $\Phi(s,t) \lambda_0 = T_{u(s,t)}^v T^*M$ for every
$(s,t)\in \overline{\R} \times [0,1]$, and $\Phi(\pm \infty,t) =
\Phi^{\pm}(t)$ for every $t\in [0,1]$.
\end{lem}

\begin{proof} (i) By the same construction used in the proof of Lemma
\ref{triv0}, we can find a smooth unitary trivialization
\[
\Psi: \overline{\R} \times \T \times \R^{2n} \rightarrow u^*(TT^*
M)
\]
such that $\Psi(s,t)\lambda_0 = T_{u(s,t)}^v T^*M$. Consider the
loops in $Sp(2n,\lambda_0)\cap U(n)$
\[
\alpha^{\pm}(t) := \Psi(\pm \infty,t)^{-1} \Phi^{\pm}(t).
\]
Since the inclusion $Sp(2n,\lambda_0)\cap U(n) \hookrightarrow
U(n)$ induces the zero homomorphism between fundamental groups,
and since $U(n)$ is connected, we can find a homotopy
\[
\alpha: \overline{\R} \times \T \rightarrow U(n)
\]
such that $\alpha(\pm \infty,t) = \alpha^{\pm}(t)$ for every $t\in
\T$. Then the unitary trivialization
\[
\Phi(s,t) := \Psi(s,t) \alpha(s,t), \quad (s,t) \in \overline{\R}
\times \T,
\]
has the required asymptotics.

(ii) By the same construction used in the proof of Lemma
\ref{triv0}, we can find a smooth unitary trivialization
\[
\Psi: \overline{\R} \times [0,1] \times \R^{2n} \rightarrow
u^*(TT^* M)
\]
such that $\Psi(s,t)\lambda_0 = T_{u(s,t)}^v T^*M$. Since the
isomorphisms $\Phi^{\pm}(t)|_{\lambda_0}$ are both orientation
preserving, the paths
\[
\alpha^{\pm}: [0,1] \rightarrow Sp(2n,\lambda_0) \cap U(n), \quad
\alpha^{\pm}(t) := \Psi(\pm \infty,t)^{-1} \Phi^{\pm}(t),
\]
take values into the same connected component of $Sp(2n,\lambda_0)
\cap U(n)$. Therefore we can find a homotopy
\[
\alpha: \overline{\R} \times [0,1] \rightarrow Sp(2n,\lambda_0) \cap
U(n)
\]
such that $\alpha(\pm \infty,t) = \alpha^{\pm}(t)$ for every $t\in
[0,1]$. Then the unitary trivialization
\[
\Phi(s,t) := \Psi(s,t) \alpha(s,t), \quad (s,t) \in \overline{\R}
\times [0,1],
\]
maps $\lambda_0$ into the vertical subbundle, and has the required
asymptotics.
\end{proof}

We denote by $W^{1,r}_{\lambda_0}$ the Sobolev space of
$\R^{2n}$-valued maps $v$ taking values in $\lambda_0$ on the boundary:
\[
W^{1,r}_{\lambda_0} (\R \times ]0,1[,\R^{2n}) := W^{1,r}_0(\R \times
]0,1[,\R^{n}) \times W^{1,r} (\R\times ]0,1[,\R^n).
\]
Let $u\in \mathcal{M}_{\Lambda}(x^-,x^+)$ (resp.\ $u\in
\mathcal{M}_{\Omega}(x^-,x^+)$),
and let $\Phi$ be a trivialization of the bundle
$u^*(TT^*M)$ as in the lemma above. Then $\Phi$ defines a conjugacy
between the fiberwise derivative of the section $\partial_{J,H}$ at $u$,
\[
D_f\partial_{J,H} (u) : T_u \mathcal{B} \rightarrow \mathcal{W}_u,
\]
and a bounded operator
\begin{eqnarray*}
D_{S,\Lambda} : W^{1,r} (\R \times \T, \R^{2n}) \rightarrow L^r (\R \times \T,
\R^{2n}) \\ \Bigl(\mbox{resp. } D_{S,\Omega} : W^{1,r}_{\lambda_0}
(\R \times ]0,1[, \R^{2n})
\rightarrow L^r (\R \times ]0,1[, \R^{2n}) \Bigr)
\end{eqnarray*}
of the form
\[
D_S v = \partial_s v - J_0 \partial_t v - S(s,t) v.
\]
Here $S$ is a smooth family of endomorphisms of $\R^{2n}$ - 1-periodic
in $t$ in the $\Lambda$ case - such that the limits
\[
S^{\pm}(t) = \lim_{s\rightarrow \pm \infty} S(s,t)
\]
are symmetric. Moreover, the solution of
\[
\frac{d}{dt} \gamma^{\pm}(t) = J_0 S^{\pm}(t) \gamma^{\pm}(t), \quad
\gamma^{\pm}(0)=I,
\]
is easily seen to be 
\[
\gamma^{\pm}(t) = \Phi(\pm \infty ,t)^{-1} \, D\phi_H^t(x^{\pm}(0)) \,
\Phi(\pm\infty,0).
\]
Finally, the requirements of Lemma \ref{ltriv} guarantee that
$\Phi(\pm \infty,\cdot)$ are symplectic trivializations of
${x^{\pm}}^*(TT^*M)$ satisfying (\ref{triv}) (resp.\
(\ref{triv2})). The following theorem is then an immediate
consequence of Theorem 7.42 and Remarks 7.44, 7.46 in \cite{rs95}
(together with the estimates of Lemma 2.4 in \cite{sal99} to deal
with the case $r>2$, see also \cite{sch95}):

\begin{thm}
\label{indice}
If $x^-,x^+\in \mathcal{P}_{\Lambda}(H)$ and
$u\in \mathcal{M}_{\Lambda}(x^-,x^+)$, the fiberwise derivative of
$\partial_{J,H}$ at $u$ is a Fredholm operator of index
\[
\ind D_f\partial_{J,H} (u) = \mu_{\Lambda} (x^-) - \mu_{\Lambda}(x^+).
\]
If $x^-,x^+\in \mathcal{P}_{\Omega}(H)$ and
$u\in \mathcal{M}_{\Omega}(x^-,x^+)$, the fiberwise derivative of
$\partial_{J,H}$ at $u$ is a Fredholm operator of index
\[
\ind D_f\partial_{J,H} (u) = \mu_{\Omega} (x^-) - \mu_{\Omega}(x^+).
\]
\end{thm}

\subsection{Coherent orientations}
\label{scoh}

The aim of this section is to show how the manifolds
$\mathcal{M}(x^-,x^+)$ can be oriented in a way which is coherent
with gluing. The construction is a particular case of the
procedure described in \cite{fh93} for an arbitrary symplectic
manifold. However, since the fact that we are dealing with the
cotangent bundle of an oriented manifold allows some slight
simplifications, we carry out the construction explicitly.

If $E$ and $F$ are real Banach spaces, $\mathrm{Fred}(E,F)$ will
denote the space of Fredholm linear operators from $E$ to $F$,
endowed with the operator norm topology. It is a Banach manifold,
being an open subset of the Banach space of all linear continuous
operators from $E$ to $F$. The manifold $\mathrm{Fred}(E,F)$ is
the base space of a smooth real line bundle, the determinant
bundle $\Det(\mathrm{Fred}(E,F))$, with fibers
\[
\Det(A) := \Lambda^{\max}(\ker A) \otimes (\Lambda^{\max}(\coker
A))^* , \quad \forall A\in \mathrm{Fred}(E,F),
\]
where $\Lambda^{\max}(V)$ denotes the exterior algebra of top
degree of the real finite dimensional vector space $V$. See
\cite{qui85} for the construction of the smooth bundle
structure\footnote{Actually, in \cite{qui85} the dual of this
object is considered. The convention we use here is more common in
symplectic geometry.} of $\Det(\mathrm{Fred}(E,F))$. Two
isomorphisms $\Phi:E \cong E^{\prime}$ and $\Psi:F \cong
F^{\prime}$ induce a canonical smooth line bundle isomorphism
$\Det(E,F) \cong \Det(E^{\prime},F^{\prime})$ lifting the
diffeomorphism $\mathrm{Fred}(E,F) \rightarrow
\mathrm{Fred}(E^{\prime},F^{\prime})$, $A \mapsto \Psi A
\Phi^{-1}$. If $E_0\subset E$ is a closed finite codimensional
linear subspace and $A\in \mathrm{Fred}(E,F)$, the restriction
$A|_{E_0}$ belongs to $\mathrm{Fred}(E_0,F)$, and the exact
sequence
\[
0 \rightarrow \ker A|_{E_0} \rightarrow \ker A \rightarrow E/E_0
\stackrel{A}{\rightarrow} F/A(E_0) = \coker A|_{E_0} \rightarrow
F/A(E) = \coker A \rightarrow 0
\]
determines a canonical isomorphism (see Lemma 18 in \cite{fh93})
\begin{equation}
\label{co0} \Det (A) \cong \Det(A|_{E_0}) \otimes \Lambda^{\max}
(E/E_0).
\end{equation}

Fix some $r>1$. Let $\Sigma_{\Lambda}$ be the set of operators
\[
D_{S,\Lambda} : W^{1,r}(\R \times \T,\R^{2n}) \rightarrow L^r(\R
\times \T,\R^{2n})
\]
of the form
\begin{equation}
\label{co1}
v \mapsto \partial_s v - J_0 \partial_t v - S(s,t) v,
\end{equation}
with $S\in C^0(\overline{\R} \times \T, \mathfrak{gl}(2n))$ such
that $S(\pm \infty,t)$ are symmetric for every $t\in \T$, and the
paths of symplectic matrices $\gamma_S^-$ and $\gamma_S^+$ solving
\[
\frac{d}{dt} \gamma_S^{\pm} (t) = J_0 S(\pm \infty,t)
\gamma_S^{\pm} (t), \quad \gamma_S^{\pm}(0) = I,
\]
satisfy
\begin{equation}
\label{co2} 1 \notin \sigma (\gamma_S^{\pm} (1)).
\end{equation}
Then $D_{S,\Lambda}$ is a Fredholm operator of index
\[
\ind D_{S,\Lambda} = \mu_{CZ} (\gamma_S^-) - \mu_{CZ}
(\gamma_S^+).
\]
Similarly, $\Sigma_{\Omega}$ will denote the set of operators
\[
D_{S,\Omega} : W^{1,r}_{\lambda_0} ( \R \times ]0,1[,\R^{2n})
\rightarrow L^r  (\R \times ]0,1[,\R^{2n}),
\]
of the form (\ref{co1}), where $S\in C^0(\overline{\R}\times
[0,1], \mathfrak{gl}(2n))$, $S(\pm \infty,t)$ symmetric, is such
that the paths of symplectic matrices $\gamma_S^-$ and
$\gamma_S^+$ satisfy
\begin{equation}
\label{co3} \gamma_S^{\pm}(1) \lambda_0 \cap \lambda_0 = (0).
\end{equation}
They are Fredholm operators of index
\[
\ind D_{S,\Omega} = \mu(\gamma_S^-(\cdot) \lambda_0,\lambda_0) -
\mu(\gamma_S^+(\cdot) \lambda_0,\lambda_0).
\]
The paths of symmetric matrices $S(\pm \infty,\cdot)$ satisfying
(\ref{co2}) (resp.\ (\ref{co3})) will be called simply {\em
non-degenerate paths}. As usual, we shall omit the subscript
$\Lambda$ or $\Omega$ when we wish to consider both situations at
the same time.

The determinant bundle over $\mathrm{Fred}(W^{1,r},L^r)$ restricts
to a line bundle over $\Sigma$, which is non-trivial on some
connected components of $\Sigma$. Actually, the determinant bundle
may be non-trivial on such simple sets as
\[
\set{A^{-1} D_S A}{A\in Sp(2n,\lambda_0) \cap U(n)},
\]
for a fixed $D_S\in \Sigma$. Indeed, if $n$ is even we can find a
path $\alpha:[0,1]\rightarrow Sp(2n,\lambda_0) \cap U(n)$ such
that $\alpha(0)=I$ and $\alpha(1)=-I$, and it is easy to show (see
Theorem 2 of \cite{fh93}) that the restriction of the determinant
bundle to the loop
\[
\set{\alpha(\lambda)^{-1} D_{S,\Lambda} \alpha(\lambda)}{\lambda
\in [0,1]}
\]
is trivial if and only if $\ind D_{S,\Lambda}$ is even.

However, the determinant bundle becomes trivial when we fix the
asymptotics: if $S^+$ and $S^-\in C^0(\T,Sym(2n))$ (resp.\
$C^0([0,1],Sym(2n))$) are non-degenerate paths, we can consider
the subset of $\Sigma$,
\[
\Sigma(S^-,S^+) := \set{D_S \in \Sigma}{S(\pm \infty,t) =
S^{\pm}(t)},
\]
consisting of those operators having fixed asymptotics. The set
$\Sigma(S^-,S^+)$ is contractible (\cite{fh93}, Proposition 7), so
the determinant bundle restricts to a trivial bundle on it, which
we denote by $\Det(\Sigma(S^-,S^+))$.

Two orientations $o(S_1,S_2)$ and $o(S_2,S_3)$ of
$\Det(\Sigma(S_1,S_2))$ and $\Det(\Sigma(S_2,S_3))$, respectively,
induce in a canonical way an orientation
\[
o(S_1,S_2) \#\,\, o(S_2,S_3)
\]
of $\Det(\Sigma(S_1,S_3))$ (see \cite{fh93}, section 3). Such an
orientation is associative, meaning that
\[
(o(S_1,S_2) \#\,\, o(S_2,S_3)) \#\, o(S_3,S_4) = o(S_1,S_2) \#\,
(o(S_2,S_3) \#\, o(S_3,S_4)).
\]
A {\em coherent orientation for $\Sigma$} is a set of orientations
$o(S^-,S^+)$ of $\Det(\Sigma(S^-,S^+))$ for each pair $(S^-,S^+)$
of non-degenerate paths, such that
\begin{equation}
\label{co4} o(S_1,S_3) = o(S_1,S_2) \#\, o(S_2,S_3),
\end{equation}
for each triplet $(S_1,S_2,S_3)$ of non-degenerate paths. The
existence of coherent orientations for $\Sigma_{\Lambda}$ is
established in \cite{fh93}, Theorem 12. The construction for
$\Sigma_{\Omega}$ is identical.

Let us fix unitary trivializations $\Phi_x$ of $x^*(TT^*M)$
satisfying $\Phi_x(t)\lambda_0 = T_{x(t)}^v T^*M$, for each $x\in
\mathcal{P}(H)$. In the $\Omega$ case we also require the
isomorphisms $\Phi_x(t)|_{\lambda_0}$ to be orientation
preserving. Correspondingly, we obtain the non-degenerate path
$S_x \in C^0(\T,Sym(2n))$ (resp.\ $C^0([0,1],Sym(2n))$) such that
\[
\gamma_{S_x}(t) = \Phi_x(t)^{-1} D \phi_H^t (x(0)) \Phi_x(0).
\]
Let us fix also a coherent orientation for $\Sigma$. We shall see
that these data determine an orientation of
\[
\Det(D_f \partial_{J,H}(u)),
\]
the determinant of the fiberwise derivative of the section
\begin{equation}
\label{sez}
\partial_{J,H}: \mathcal{B}(x^-,x^+) \rightarrow
\mathcal{W}(x^-,x^+)
\end{equation}
at every $u\in \mathcal{M}(x^-,x^+)$, for
every pair $x^-,x^+\in \mathcal{P}(H)$.

Let $x^-,x^+\in \mathcal{P}(H)$ and let $u\in
\mathcal{M}(x^-,x^+)$. By Lemma \ref{ltriv}, we can find a smooth
unitary trivialization $\Phi_u$ of $u^*(TT^*M)$ agreeing with
$\Phi_{x^-}$ and $\Phi_{x^+}$ for $s=-\infty$ and $s=+\infty$,
respectively. In the $\Omega$ case we also require that
$\Phi_u(s,t) \lambda_0 = T_{u(s,t)}^v T^*M$ for every $(s,t)\in
\overline{\R}\times [0,1]$. Then $D_f \partial_{J,H}(u)$ is
conjugated by $\Phi_u$ to an operator $D_S$ belonging to
$\Sigma(S_{x^-},S_{x^+})$. So we have a canonical isomorphism
\[
\Det(D_f \partial_{J,H}(u)) \cong \Det (D_S),
\] 
and $\Det(D_f \partial_{J,H}(u))$ inherits an orientation from
$o(S_{x^-},S_{x^+})$.

Changing the trivialization $\Phi_u$ by another one with the same
properties changes $D_S$ by a unitary conjugacy $\Psi D_S
\Psi^{-1}$, where $\Psi\in C^0(\overline{\R} \times \T,U(n))$
(resp.\ $C^0(\overline{\R} \times [0,1],Sp(2n,\lambda_0) \cap
U(n))$) is the identity for $s=\pm \infty$, so Lemma 13 in
\cite{fh93} implies that the orientation of $\Det(D_f
\partial_{J,H}(u))$ does not depend on the choice of $\Phi_u$.

Such an orientation varies continuously with $u$, so if the
section (\ref{sez})
is transverse to the zero section,
$\mathcal{M}(x^-x^+)$ is a finite dimensional manifold and
\[
\Lambda^{\max} (T_u \mathcal{M}(x^-,x^+)) = \Lambda^{\max} (\ker
D_f \partial_{J,H}(u)) = \Det(D_f \partial_{J,H}(u)) \otimes \R
\]
is oriented, meaning that $\mathcal{M}(x^-,x^+)$ is oriented.

In particular when $\mu(x^-)-\mu(x^+)=1$, $\mathcal{M}(x^-,x^+)$
is an oriented one-dimensional manifold. Since translation of the $s$
variable defines a free $\R$ action on it,
$\mathcal{M}(x^-,x^+)$ consists of lines. Denoting by $[u]$
the equivalence class of $u$ in the zero-dimensional manifold
$\mathcal{M}(x^-,x^+)/\R$, we define
\[
\epsilon([u]) \in \{-1,+1\}
\]
to be $+1$ if the $\R$-action is orientation preserving on the
connected component of $\mathcal{M}(x^-,x^+)$ containing $u$, $-1$
in the opposite case.

\begin{rem}
In our construction, the orientations of moduli spaces of
solutions of the Floer equation depend on the choice of suitable
trivializations of $x^*(TT^*M)$, for every $x\in \mathcal{P}(H)$,
and of a coherent orientation for $\Sigma$. This approach is
possible because here we can find trivializations of $u^*(TT^*M)$
with prescribed asymptotics, something which is not possible for
an arbitrary symplectic manifold. So here we do not need to
introduce the notion of a coherent orientation for the symplectic
vector bundle $TT^*M \rightarrow T^*M$, as in \cite{fh93}. The use
of a coherent orientation for $TT^*M \rightarrow T^*M$ would allow
to drop the orientability assumption on $M$.
\end{rem}

\subsection{$\mathbf{L^{\infty}}$ estimates}
\label{cs}

In order to have $L^{\infty}$ bounds on the set of solutions of
(\ref{cr},\ref{bdry}) with bounded action,
further assumptions on the Hamiltonian $H$ and
on the almost complex structure $J$ are
needed. Let us fix a metric $\langle \cdot ,\cdot \rangle$ on $M$. We
shall denote by the same symbol the induced metric on $TM$ and on
$T^* M$, and by
$\nabla$ the corresponding Levi-Civita covariant derivation.
All the $L^r$ and Sobolev norms we will use refer to this metric.
This metric determines an isometry $TM \rightarrow
T^*M$, and a direct summand of the vertical bundle
$T^v T^* M$, the horizontal bundle $T^h T^* M$, together with
isomorphisms
\[
T_x T^*M = T^h_x T^*M \oplus T^v_x T^*M \cong T_q M \oplus T_q^* M
\cong T_q M \oplus T_q M,
\quad x=(q,p)\in T^* M.
\]
There is a preferred $\omega$-compatible almost complex structure
$\widehat{J}$ on $T^* M$, which in the above splitting has the form
\[
\widehat{J} =
\left( \begin{array}{cc} 0 & I \\ -I & 0 \end{array}
\right).
\]
The Liouville and the Hamiltonian vector field can be written as
\[
\eta (q,p) = (0,p), \quad
X_H(t,q,p) = \widehat{J} \nabla H(t,q,p) = (\nabla_p H(t,q,p),-
\nabla_q H(t,q,p)),
\]
where $\nabla_q$ and $\nabla_p$ denote the horizontal and the vertical
components of the gradient.
We shall make the following assumptions (recall that $\eta$ denotes
the Liouville vector field on $T^*M$):

\begin{description}

\item[(H1)] there exist $h_0>0$ and $h_1\geq 0$ such that
\[
dH(t,q,p) [\eta] - H(t,q,p) \geq h_0 |p|^2 - h_1,
\]
for every $(t,q,p)\in \T \times T^* M$;

\item[(H2)] there exists $h_2\geq 0$ such that
\[
|\nabla_q H(t,q,p)| \leq h_2 (1+|p|^2), \quad |\nabla_p H(t,q,p)|
 \leq h_2 (1+|p|),
\]
for every $(t,q,p)\in \T \times T^* M$.

\end{description}

Condition (H1) is assumed also in \cite{cie94}, and it is 
a condition of quadratic growth at infinity: thanks
to the compactness of $M$, it easily implies the estimate
\[
H(t,q,p) \geq \frac{1}{2} h_0 |p|^2 - h_3,
\]
for a suitable constant $h_3$. Condition (H1)
does not depend on the choice of the metric on $M$: if $\langle \cdot,
\cdot \rangle_*$ is another metric, by the compactness of $M$ $|\xi|
\leq c |\xi|_*$, so if $H$ satisfies (H1) with respect to $|\cdot|$
with constants $h_0,h_1$, it also satisfies (H1) with respect to
$|\cdot|_*$ with constants $h_0/c^2$ and $h_1$.

We will show that also (H2) does not depend on the metric, by checking
that $H$ satisfies (H2) if and
only if for any coordinate system $(q^1,\dots,q^n)\in \R^n$ on
$U\subset M$ - inducing the coordinate system
$(q^1,\dots,q^n,p_1,\dots,$ $p_n) \in \R^n \times {\R^n}^*$ on $T^*U
\subset T^*M$ - there is $a\geq 0$ such that 
\begin{equation}
\label{h2loc}
|\partial_{q^i} H(t,q,p)| \leq a (1+|p|^2), \quad
|\partial_{p_i} H(t,q,p)| \leq a (1+|p|), \quad \forall i = 1.\dots,n.
\end{equation}
Here $|\cdot|$ denotes any norm, for instance the Euclidean one, on
${\R^n}^*$. It is readily seen that if (\ref{h2loc}) holds for $H$ and
$\psi$ is a change of coordinates on $\R^n$, then (\ref{h2loc}) holds
for $H(t,\psi(q),p \circ D\psi(q)^{-1})$ (with a different constant $a$),
hence this local condition is independent of the choice of the
coordinate system. 

Let $K:TT^*M \rightarrow T^*M$ be the connection associated to the
metric $\langle \cdot , \cdot \rangle$. Then the horizontal and
vertical components of the gradient of $H$ are
\begin{equation}
\label{iom1}
\nabla_q H = D\tau^* \nabla H, \quad \nabla_p H = K \nabla H.
\end{equation}
In the coordinate system $(q^1,\dots,q^n,p_1,\dots p_n)$ the
connection $K$ has the form 
\begin{equation}
\label{iom2}
K (\xi,\zeta) = \zeta - B \xi, \quad (\xi,\zeta)\in \R^n \times {\R^n}^*,
\end{equation}
where
$B=B(q,p)\in \mathrm{Hom}(\R^n,{\R^n}^*)$ is symmetric and
depends linearly on $p$. 
If the symmetric operator $G(q)\in \mathrm{Hom}(\R^n,{\R^n}^*)$
represents the metric on $M$ in the local coordinates
$(q^1,\dots,q^n)$, the induced metric on $T^*M$ has the local
expression
\[
\langle (\xi_1,\zeta_1), (\xi_2,\zeta_2) \rangle = (G \xi_1) \xi_2 +
(\zeta_1 - B \xi_1) (G^{-1} (\zeta_2 - B\xi_2)),
\]
which can be rewritten in matrix form as
\begin{equation}
\label{iom3}
\tilde{G} = \left( \begin{array}{cc} G + B G^{-1} B & - B G^{-1} \\ - G^{-1} B
    & G^{-1} \end{array} \right).
\end{equation}
The inverse of this matrix is
\begin{equation}
\label{iom4}
\tilde{G}^{-1} =
\left( \begin{array}{cc} G^{-1} & G^{-1}B \\ B G^{-1} 
    & G + BG^{-1} B \end{array} \right).
\end{equation}
Since the local expression of the gradient of $H$ is 
$\nabla H = \tilde{G}^{-1} \partial H$, by (\ref{iom1}), (\ref{iom2}),
(\ref{iom3}), (\ref{iom4}),
\[
\left\{ \begin{array}{l} \nabla_q H = G^{-1} \partial_q H + G^{-1}
  B \partial_p H , \\ \nabla_p H = G \partial_p H,
\end{array} \right.
\quad \left\{ \begin{array}{l} \partial_q H = G \nabla_q H - B
  G^{-1} \nabla_p H, \\ \partial_p H = G^{-1} \nabla_p H.
\end{array} \right.
\]
Since $M$ is compact and $B$ depends linearly on $p$, the above
formulas show the equivalence of (H2) and (\ref{h2loc}).   

Condition (H2) is weaker than the corresponding growth condition appearing
in \cite{cie94}. Physical Hamiltonians of the form
\[
H(t,q,p) = \frac{1}{2} |T(t,q)p-A(t,q)|^2 + V(t,q)
\]
satisfy (H1) and (H2), provided that the symmetric tensor $T^*T$
is everywhere positive. Condition (H2) implies the estimate
\begin{equation}
\label{gs}
|X_H(t,q,p)| = |\nabla H(t,q,p)| \leq h_4 (1+|p|^2),
\end{equation}
for a suitable constant $h_4$.
Here is a first important consequence of assumptions (H1) and (H2):

\begin{lem}
\label{lem0}
Assume that $H$ satisfies (H0), (H1), and (H2). Then
for every $a\in \R$, the set of solutions $x\in \mathcal{P}(H)$ such
that $\mathcal{A}(x)\leq a$ is finite.
\end{lem}

\begin{proof}
Let $x=(q,p)\in \mathcal{P}(H)$ be such that $\mathcal{A}(x)\leq
a$. Then by (H1),
\begin{eqnarray*}
a \geq \mathcal{A}(x) = \int_0^1 \left( \theta(\dot{x}) - H(t,x)
\right)\, dt = \int_0^1 \left( d\theta (\eta,X_H(t,x)) -
H(t,x)\right)\, dt \\
= \int_0^1 \left( dH(t,q,p) [\eta] - H(t,q,p) \right)\, dt \geq
h_0 \|p\|_{L^2}^2 - h_1,
\end{eqnarray*}
so $\mathcal{P}(H)\cap \{\mathcal{A}\leq a\}$ is bounded in $L^2$.
By (\ref{gs}) we also have
\[
|\dot{x}| = |X_H(t,x)| \leq h_4 (1+|p|^2),
\]
from which we conclude that $\mathcal{P}(H)\cap \{\mathcal{A}\leq a\}$
is bounded in $W^{1,1}$, hence in $L^{\infty}$. In particular, the set
  $\set{x(0)}{x\in \mathcal{P}(H), \; \mathcal{A}(x)\leq a}$ is
  pre-compact in $T^*M$, and being discrete by (H0), it must be finite.
\end{proof}

We shall prove that if $H$
satisfies (H1), (H2), and $J$ is close to $\widehat{J}$, then the solutions of
(\ref{cr},\ref{bdry}) with bounded action are uniformly bounded in
$L^{\infty}$. We shall need the following interpolation inequality:

\begin{lem}
\label{interp}
There exists $C>0$ such that
\[
\|\varphi\|_{L^4(\R \times ]0,1[)}^4 \leq C
\|\varphi\|_{L^2(\R \times ]0,1[)}^2
\|\varphi\|_{W^{1,2}(\R \times ]0,1[)}^2
\]
for every $\varphi\in W^{1,2}(\R\times ]0,1[)$.
\end{lem}

\begin{proof}
This is an easy consequence of the interpolation inequality
\begin{equation}
\label{maz}
\|\psi\|_{L^4(\R^2)}^4 \leq C_0 \|\psi\|_{L^2(\R^2)}^2 \|\nabla
\psi\|_{L^2(\R^2)}^2 \quad\forall \psi\in W^{1,2}(\R^2),
\end{equation}
proved in \cite{maz85}, section 1.4.7.
Indeed, if $\varphi\in W^{1,2}(\R \times ]0,1[)$, by reflection
along the lines $\R
\times \{0\}$, and $\R \times \{1\}$ we obtain a function
$\tilde{\varphi}\in W^{1,2}(\R\times ]-1,2[)$ such that
\[
\| \tilde{\varphi}\|_{L^2(\R\times ]-1,2[)}^2 = 3 \|
\varphi\|_{L^2(\R\times ]0,1[)}^2, \quad \| \nabla
\tilde{\varphi}\|_{L^2(\R\times ]-1,2[)}^2 = 3 \| \nabla
\varphi\|_{L^2(\R\times ]0,1[)}^2.
\]
Let $\chi\in C^{\infty}(\R)$ be a function such that $\chi=1$ on
$[0,1]$, $\supp \chi\subset ]-1,2[$, $0\leq \chi\leq 1$, and
$|\chi^{\prime}| \leq 2$. Then the function $\psi(s,t) = \chi(t)
\tilde{\varphi}(s,t)$ belongs to $W^{1,2}(\R^2)$, so (\ref{maz})
leads to
\[
\|\varphi\|_{L^4(\R\times ]0,1[)}^4 \leq \|\psi\|_{L^4(\R^2)}^4
\leq C_0 \|\psi\|_{L^2(\R^2)}^2 \|\nabla \psi\|_{L^2(\R^2)}^2,
\]
and the conclusion follows from the inequalities
\[
\|\psi\|_{L^2(\R^2)} \leq \|\tilde{\varphi}\|_{L^2(\R\times
]-1,2[)} = \sqrt{3} \|\varphi\|_{L^2(\R\times ]0,1[)},
\]
and
\begin{eqnarray*}
\|\nabla \psi\|_{L^2(\R^2)} = \|\chi \nabla \tilde{\varphi} +
\chi^{\prime} \tilde{\varphi} \|_{L^2(\R\times ]-1,2[)} \leq
\|\nabla \tilde{\varphi} \|_{L^2(\R \times ]-1,2[)} + 2
\|\tilde{\varphi}\|_{L^2(\R \times ]-1,2[)} \\
= \sqrt{3} \|\nabla \varphi\|_{L^2(\R\times ]0,1[)} + 2
\sqrt{3} \|\varphi\|_{L^2(\R\times ]0,1[)}.
\end{eqnarray*}
\end{proof}

The main step to prove the $L^{\infty}$ estimates is provided by the
following:

\begin{lem}
\label{gradest}
Assume that $H$ satisfies (H1), (H2), and that the $t$-dependent
1-periodic almost complex structure $J$ on $T^*M$ satisfies
$\|J\|_{\infty} < +\infty$. For every
pair of real numbers $a_1,a_2$ there exists a number $c$ such that
for every $u=(q,p)\in C^{\infty} (\R \times [0,1],T^*M)$ (resp.\
$u=(q,p) \in C^{\infty} (]0,+\infty[ \times [0,1],T^*M)\cap
 W^{1,r}(]0,1[\times ]0,1[,T^*M)$, with $r>2$) solving
\begin{equation}
\label{nscr}
\partial_s u - J(t,u) (\partial_t u - X_H(t,u)) =0,
\end{equation}
and such that $a_1\leq \mathcal{A}(u(s,\cdot)) \leq a_2$ for every $s\in
\R$ (resp.\ $s\in [0,+\infty[$), there holds
\[
\|p\|_{L^2(I \times ]0,1[)} \leq c |I|^{\frac{1}{2}}, \quad
\|\nabla p\|_{L^2(I \times ]0,1[)} \leq c (|I|^{\frac{1}{2}} +1),
\]
for every interval $I\subset \R$ (resp.\ $I\subset [0,+\infty[$).
\end{lem}

\begin{proof} We shall denote by $\mathcal{U}$ the set of solutions $u\in
C^{\infty} ( \R \times [0,1],T^*M)$ (respectively $u\in C^{\infty}
( ]0,+\infty[ \times [0,1] ,T^*M)\cap W^{1,r}(]0,1[\times
]0,1[,T^*M)$) of (\ref{nscr}) such that $a_1 \leq \mathcal{A}
(u(s,\cdot)) \leq a_2$ for every $s\in \R$ (resp.\ $s\in
[0,+\infty[$). Notice that in the case of $u$ defined on the
half-strip, our assumptions imply that $u(0,\cdot)\in
W^{1-1/r,r}(]0,1[,T^*M)$ (see \cite{ada75} section 7.56), from
which it easily follows that the function $s\mapsto
\mathcal{A}(u(s,\cdot))$ is continuous on $[0,+\infty[$.

\medskip

{\sc Claim 1.} There exists $c_1$ such that for every $u\in
\mathcal{U}$, 
\[
\|\partial_s u\|_{L^2(\R
  \times ]0,1[)} \leq c_1 \quad \mbox{(resp. }   
\|\partial_s u\|_{L^2(]0,+\infty[ \times ]0,1[)} \leq c_1\mbox{)}.
\]

\medskip

Indeed, for $s_0<s_1$,
\begin{eqnarray*}
\|\partial_s u\|_{L^2 (]s_0,s_1[ \times ]0,1[)}^2 = \int_{s_0}^{s_1}
  \int_0^1 |\partial_s u|^2 \, dt\,ds \leq \|J^{-1}\|_{\infty}^2
  \int_{s_0}^{s_1} \int_0^1 |\partial_s u|^2_{J_{t}} \, dt\, ds \\
  = \|-J \|_{\infty}^2 \int_{s_0}^{s_1} \int_0^1 \langle -
  \nabla_{J} \mathcal{A}(u(s,\cdot))(t), \partial_s u(s,t)
  \rangle_{J_{t}} \, dt\, ds = - \|J\|^2_{\infty}
  \int_{s_0}^{s_1} d\mathcal{A} (u(s,\cdot)) [\partial_s
  u(s,\cdot)] \, ds \\ = \|J\|^2_{\infty} (
  \mathcal{A}(u(s_0,\cdot)) - \mathcal{A}(u(s_1,\cdot))) \leq
  \|J\|_{\infty}^2 (a_2-a_1). 
\end{eqnarray*}

\medskip

{\sc Claim 2.} There exists $c_2$ such that $\|p(s,\cdot)\|_{L^2 (]0,1[)}
\leq c_2 (1+ \|\partial_s u (s,\cdot)\|_{L^2(]0,1[)} )$ for every
$u=(q,p) \in \mathcal{U}$ and every $s\in \R$ (resp.\ $s\in [0,+\infty[$).

\medskip

Indeed, since $u$ solves (\ref{nscr}),
\begin{eqnarray*}
\theta(\partial_t u) = \theta(X_H(t,u) - J(t,u) \partial_s
u) = \omega(\eta(u), X_H(t,u) - J(t,u) \partial_s u) \\
= dH(t,u) [\eta(u)] - \langle \widehat{J}(u) \eta(u), J(t,u)
\partial_s u \rangle.
\end{eqnarray*}
Then by (H1) and by the fact that $|\eta(q,p)|=|p|$,
\[
\theta(\partial_t u) - H(t,u) \geq h_0 |p|^2 - h_1 - \|J\|_{\infty}
|p| \, |\partial_s u|,
\]
and integrating over $[0,1]$ we find
\begin{eqnarray*}
a_2 \geq \mathcal{A}(u(s,\cdot)) = \int_0^1 ( \theta (\partial_t u
(s,t)) - H(t,u(s,t)) )\, dt \\ \geq h_0 \|p(s,\cdot)\|_{L^2(]0,1[)}^2 -
h_1 - \|J\|_{\infty} \|p(s,\cdot)\|_{L^2(]0,1[)} \|\partial_s u(s,\cdot)
\|_{L^2(]0,1[)},
\end{eqnarray*}
which implies Claim 2.

\medskip

{\sc Claim 3.} There exists $c_3$ such that
$\|p(s,\cdot)\|_{L^{\infty}(]0,1[)} \leq c_3 ( 1 + \|\partial_s u
(s,\cdot)\|^2_{L^2(]0,1[)} )$ for every $u=(q,p)\in \mathcal{U}$ and
every $s\in \R$ (resp.\ $s\in [0,+\infty[$).

\medskip

Indeed, by (\ref{gs}),
\begin{eqnarray*}
\| \partial_t p (s,\cdot) \|_{L^1(]0,1[)} \leq \|\partial_t u(s,\cdot)
\|_{L^1(]0,1[)} \leq \|X_H(\cdot,u(s,\cdot))\|_{L^1(]0,1[)} +
\|J(\cdot,u) \partial_s u(s,\cdot) \|_{L^1(]0,1[)} \\
\leq h_4(1+ \|p(s,\cdot)\|^2_{L^2(]0,1[)} ) + \|J\|_{\infty} \|\partial_s
u(s,\cdot) \|_{L^1(]0,1[)},
\end{eqnarray*}
which can be estimated by Claim 2 by
\[
\leq h_4 (1+c^2_2 (1+\|\partial_s u(s,\cdot)\|_{L^2(]0,1[)} )^2 ) +
\|J\|_{\infty} \|\partial_s u(s,\cdot)\|_{L^2(]0,1[)}.
\]
Therefore, the $W^{1,1}$ norm of $p(s,\cdot)$ on $]0,1[$ is bounded by
a quadratic function of $\|\partial_s u(s,\cdot)\|_{L^2(]0,1[)}$. Then the
same is true a fortiori for the $L^{\infty}$ norm of
$p(s,\cdot)$.

\medskip

{\sc Claim 4.} For every $\delta>0$ there is a number $m(\delta)$ with
the following property: for every $u=(q,p)\in \mathcal{U}$ the closed 
subset of $\R$ (resp.\ of $[0,+\infty[$),
\[
S_{\delta}(u) := \set{s}{\|p(s,\cdot)\|_{L^{\infty}(]0,1[)} \leq
  m(\delta)}
\]
has non-empty intersection with any closed interval $I\subset
\R$ (resp.\ $I\subset [0,+\infty[$) of length
$\delta$. 

\medskip

Indeed, for every $s_0\in \R$ (resp.\ $s_0\in [0,+\infty[$),
\begin{eqnarray*}
\min_{s\in [s_0,s_0+\delta]} \|\partial_s u(s,\cdot) \|_{L^2(]0,1[)}^2
\leq \frac{1}{\delta} \int_{s_0}^{s_0+\delta} \|\partial_s
u(s,\cdot)\|_{L^2(]0,1[)}^2 \, ds = \frac{1}{\delta} \|\partial_s
u\|_{L^2(]s_0,s_0+\delta[ \times ]0,1[)}^2 \leq \frac{1}{\delta} c_1^2
\end{eqnarray*}
because of Claim 1. Then Claim 3 implies Claim 4 with
\[
m(\delta) = c_3 \left( 1 + \frac{c_1^2}{\delta} \right).
\]

\medskip

{\sc Claim 5.} There exists $c_4$ such that $\|p(s,\cdot)\|_{L^2(]0,1[)}
\leq c_4$ for every $u=(q,p)\in \mathcal{U}$ and every $s\in \R$
(resp.\ $s\in [0,+\infty[$).

\medskip

Given $u=(q,p)\in \mathcal{U}$ and $s\in \R$ (resp.\ $s\in
[0,+\infty[$),
let $s_0$ be an element of $S_1(u)$ such that $|s-s_0|\leq 1$
(see Claim 4). Then  
\begin{eqnarray*}
\| p(s,\cdot) \|^2_{L^2(]0,1[)} = \|p(s_0,\cdot)\|_{L^2(]0,1[)}^2 +
\int_{s_0}^s \frac{d}{d\sigma} \|p(\sigma,\cdot)\|^2_{L^2(]0,1[)} \,
d\sigma \\ = \|p(s_0,\cdot)\|_{L^2(]0,1[)}^2 + 2 \int_{s_0}^s \int_0^1
\langle p(\sigma,t), \partial_s p(\sigma,t) \rangle\, dt\, d\sigma \\
\leq m(1)^2 + 2 \left| \int_{s_0}^s \|p(\sigma,\cdot)\|_{L^2(]0,1[)}^2 \,
d\sigma \right|^{\frac{1}{2}} \|\partial_s p\|_{L^2(]s_0,s[ \times ]0,1[)}.
\end{eqnarray*}
By Claim 1, $\|\partial_s p\|_{L^2(]s_0,s[ \times ]0,1[)} \leq \|\partial_s
u\|_{L^2(]s_0,s[ \times ]0,1[)} \leq c_1$, and using also Claim 2 we get
\begin{eqnarray*}
\|p(s,\cdot)\|_{L^2(]0,1[)}^2 \leq m(1)^2 + 2c_1 \left| \int_{s_0}^s
c_2^2(1+ \|\partial_s u(\sigma,\cdot)\|_{L^2(]0,1[)} )^2 \, d\sigma
\right|^{\frac{1}{2}} \\ \leq m(1)^2 + 2 c_1 c_2 \left( 2|s-s_0| + 2\|
\partial_s u \|_{L^2(]s_0,s[ \times ]0,1[)}^2 \right)^{\frac{1}{2}} \leq
m(1)^2 + 2c_1 c_2 (2+2c_1^2)^{\frac{1}{2}}.
\end{eqnarray*}

\medskip

{\sc Conclusion.} Let $u=(q,p)\in \mathcal{U}$ and let
$I\subset \R$ (resp.\ $I\subset
[0,+\infty[$) be an interval. By Claim 5,
\[
\|p\|_{L^2(I \times ]0,1[)}^2 = \int_I \|p(s,\cdot)\|_{L^2(]0,1[)}^2 \, ds
\leq c_4^2|I|,
\]
so it is enough to estimate the $L^2$ norm of $\nabla p$. By
(\ref{gs}),
\begin{eqnarray*}
|\nabla p|^2 \leq |\nabla u|^2 = |\partial_s u|^2 + | \partial_t u|^2
 = |\partial_s u|^2 + |X_H - J\partial_s u|^2  \leq
 (1+2\|J\|_{\infty}^2) |\partial_s u|^2 + 2h_4^2 (1+|p|^2)^2,
\end{eqnarray*}
which implies
\[
|\nabla p|^2 \leq b_1 (1 +|\partial_s u|^2 + |p|^4),
\]
for a suitable constant $b_1$. Integrating this inequality over
$]s_0,s_1[\times ]0,1[$ and using Claim 1 we get
\begin{equation}
\label{la22}
\|\nabla p\|_{L^2(]s_0,s_1[ \times ]0,1[)}^2 \leq b_1 (|s_1-s_0| + c_1^2)
+ b_1 \|p\|_{L^4(]s_0,s_1[ \times ]0,1[)}^4.
\end{equation}
Let $\delta$ be a positive number, to be fixed later, and let
$S_{\delta}(u)$ be the $\delta$-dense subset of $\R$ (resp.\ of
$[0,+\infty[$) provided by Claim
4.  Let $s_1\in S_{\delta}(u)$, and set $s_0=s_1-\delta$ (resp.\ $s_0=
\max\{s_1-\delta,0\}$). The real valued function $(|p(s,t)|-m(\delta))^+$
vanishes for $s=s_1$, so by reflection along the line $\{s_0\}\times
\R$ we obtain the function
\[
\varphi(s,t) = \begin{cases} (|p(s,t)|-m(\delta))^+ & \mbox{if } s\in
  [s_0,s_1], \\ (|p(2s_0 -s,t)-m(\delta))^+ & \mbox{if } s\in
  [2s_0-s_1,s_0] \\ 0 & \mbox{if } s\in \R \setminus [2s_0-s_1,s_1],
\end{cases}
\]
which belongs to $W^{1,2}(\R \times ]0,1[)$ and satisfies
\begin{eqnarray}
\label{eq0}
\|\varphi\|_{L^2(\R \times ]0,1[)}^2 = 2
\|(|p|-m(\delta))^+\|_{L^2(]s_0,s_1[\times ]0,1[)}^2, \\
\label{eq1}
\|\varphi\|_{L^4(\R \times ]0,1[)}^4 = 2
\|(|p|-m(\delta))^+\|_{L^4(]s_0,s_1[\times ]0,1[)}^4,\\
\label{eq2} \begin{split}
\|\nabla\varphi\|_{L^2(\R\times ]0,1[)}^2 = 2
\|\nabla (|p|-m(\delta))^+\|_{L^2(]s_0,s_1[\times ]0,1[)}^2 \\ \leq
2\|\nabla |p| \|_{L^2(]s_0,s_1[\times ]0,1[)}^2 \leq 2 \|\nabla
p\|_{L^2(]s_0,s_1[ \times ]0,1[)}^2. \end{split}
\end{eqnarray}
Since $(a+b)^4\leq 8(a^4+b^4)$, by (\ref{eq1}) we have
\begin{eqnarray*}
\|p\|^4_{L^4 (]s_0,s_1[ \times ]0,1[)} =
\int_{\substack{]s_0,s_1[\times ]0,1[ \\ |p|\geq m(\delta)}} |p|^4 \,
ds\,dt + \int_{\substack{]s_0,s_1[\times ]0,1[ \\ |p|< m(\delta)}} |p|^4 \,
ds\,dt \\ \leq 8 \int_{\substack{]s_0,s_1[\times ]0,1[ \\ |p|\geq
    m(\delta)}} \left( (|p|-m(\delta))^4 + m(\delta)^4 \right)\,
ds\,dt + m(\delta)^4 |s_1-s_0| \\
\leq 8 \|(|p|-m(\delta))^+\|^4_{L^4(]s_0,s_1[ \times ]0,1[)}
+ 9 m(\delta)^4 |s_1-s_0| \leq 4 \| \varphi\|^4_{L^4 (\R \times
  ]0,1[)} + 9 m(\delta)^4 \delta.
\end{eqnarray*}
The interpolation estimate of Lemma \ref{interp} then implies
\begin{equation}
\label{la44}
\begin{split}
\|p\|_{L^4(]s_0,s_1[\times ]0,1[)}^4 \leq 4C
\|\varphi\|_{L^2(\R \times ]0,1[)}^2
\|\varphi\|_{W^{1,2}(\R \times ]0,1[)}^2 + 9 m(\delta)^4 \delta \\ =
4C \|\varphi\|^4_{L^2(\R \times ]0,1[)} + 4C \|\varphi\|^2_{L^2(\R
  \times ]0,1[)} \|\nabla \varphi\|^2_{L^2(\R \times ]0,1[)} + 9
m(\delta)^4 \delta.
\end{split}
\end{equation}
By (\ref{eq0}) and Claim 5,
\begin{eqnarray*}
\|\varphi\|_{L^2(\R \times ]0,1[)}^2 \leq 2\|p\|_{L^2(]s_0,s_1[\times
  ]0,1[)}^2 = 2\int_{s_0}^{s_1} \|p(s,\cdot)\|_{L^2(]0,1[)}^2\, ds  
\leq 2c_4^2  |s_1-s_0| \leq 2c_4^2 \delta.
\end{eqnarray*}
So (\ref{eq2}) and (\ref{la44}) imply
\[
\|p\|_{L^4(]s_0,s_1[\times ]0,1[)}^4 \leq 16 C
c_4^4 \delta^2 +
16 C c_4^2 \delta \|\nabla p\|_{L^2(]s_0,s_1[\times ]0,1[)}^2 + 9
m(\delta)^4 \delta.
\]
Therefore, by (\ref{la22}),
\begin{eqnarray*}
\|\nabla p\|_{L^2(]s_0,s_1[\times ]0,1[)}^2 \leq b_1 (\delta + c_1^2
+ 16 C c_4^4 \delta^2 + 9 m(\delta)^4 \delta)  + 16 b_1 C c_4^2 \delta
\|\nabla p\|_{L^2(]s_0,s_1[\times ]0,1[)}^2.
\end{eqnarray*}
Hence, if we choose $\delta$ to be $1/(32 b_1C c_4^2)$, the above
inequality implies
\[
\|\nabla p\|_{L^2(]s_0,s_1[\times ]0,1[)}^2 \leq 2b_1  (\delta + c_1^2
+ 16 C c_4^4 \delta^2 + 9 m(\delta)^4 \delta) =: b_2.
\]

We have proved that the square of the
$L^2$ norm of $\nabla p$ on the set $I
\times ]0,1[$ is bounded by $b_2$, if $I\subset \R$ (resp.\ $I\subset
[0,+\infty[$) is an
interval with the right-hand point in $S_{\delta}(u)$ and length at most
$\delta$. By the properties of $S_{\delta}(u)$, any interval in $\R$
(resp.\ in $[0,+\infty[$) of length less than $\delta$ can be covered
by two intervals with the right-hand point in $S_{\delta}(u)$ and
length at most $\delta$. Any bounded interval $I\subset \R$ (resp.\ $I
\subset [0,+\infty[$) can be covered by $\lceil |I|/\delta \rceil + 1$
intervals of length less than $\delta$, hence by $2(\lceil |I|/\delta
\rceil + 1)$ intervals with the right-hand point in $S_{\delta}(u)$
and length at most $\delta$. Therefore,
\[
\| \nabla p\|_{L^2(I\times ]0,1[)}^2 \leq 2 b_2 \Bigl(\left\lceil
  \frac{|I|}{\delta} \right\rceil + 1\Bigr) \leq 2 b_2
\left(\frac{|I|}{\delta} + 2 \right),
\]
concluding the proof.  
\end{proof}

We recall that $\lambda_0$ denotes the Lagrangian subspace $(0) \times
\R^n$ in the symplectic vector space $\R^n \times \R^n$, and that
$W^{1,r}_{\lambda_0}$ denotes the Sobolev space of $\R^{2n}$-valued
maps taking values in $\lambda_0$ on the boundary:
\[
W^{1,r}_{\lambda_0} (\Omega,\R^{2n}) := W^{1,r}_0(\Omega,\R^{n})
\times W^{1,r} (\Omega,\R^n).
\]
We recall the following consequences of the Calderon-Zygmund
inequalities for the Cauchy-Riemann operator:

\begin{prop}
\label{caldzyg} Let $\Omega$ be one of the following domains: the
cylinder $\R \times \T$, the strip $\R \times ]0,1[$, the
half-cylinder $]0,+\infty[ \times \T$, the half-strip $]0,+\infty[
\times ]0,1[$. For every $r>1$ there exists a constant
$C_r(\Omega)\geq 1$ such that
\[
\|\nabla v\|_{L^r(\Omega)} \leq C_r(\Omega) \|(\partial_s - J_0
  \partial_t) v\|_{L^r(\Omega)}
\]
for every $v\in W^{1,r}_{\lambda_0}(\Omega,\R^{2n})$.
\end{prop}

Indeed, one can start by proving the estimate 
\[
\| \nabla \varphi\|_{L^r} \leq c(r) \|\overline{\partial} \varphi\|_{L^r},
\quad \forall \varphi\in C^{\infty}_c(\R \times \T,\C), \;\;\forall
r\in ]1,+\infty[, 
\]
by the usual argument involving the fundamental
solution of the Cauchy-Riemann operator
$\overline{\partial}=\partial_s + i \partial_t$ (see e.g.\ \cite{hz94}
or \cite{ms04}, Appendix B). By Schwarz reflection, we obtain an
analogous estimate for $\varphi\in C^{\infty}([0,+\infty[ \times \T,\C)$
with real boundary conditions, and for $\varphi\in C^{\infty}(\R
\times [0,1],\C)$ with real boundary conditions. A second reflection
yields to the analogous estimate for $\varphi\in C^{\infty}([0,+\infty[
\times [0,1],\C)$ with real boundary conditions. Proposition
\ref{caldzyg} follows, by identifying $\R^{2n}$ with $\C^n$ and $J_0$
with $-i$. 

It will be useful to view $M$ as a submanifold of $\R^N$, for some
large $N$, by means of an isometric embedding $M\hookrightarrow
\R^N$, as given by Nash's theorem. Such an embedding induces also
isometric embeddings of $TM$ and $T^* M$ into $\R^{2N}$, and it is
easy to see that $\widehat{J}$ is the restriction of $J_0$.

\begin{thm}
\label{c0est}
Assume that $H$ satisfies (H1), (H2). Then there exists a number
$j_0>0$ such that, if the $t$-dependent 1-periodic
almost complex structure $J$ on $T^*M$ satisfies
$\|J-\widehat{J}\|_{\infty} < j_0$, then for every $a_1,a_2\in \R$
there holds: 
\begin{enumerate}

\item the set of solutions $u=(q,p)\in C^{\infty}(\R \times \T, T^*M)$ of
\begin{equation}
\label{nncr}
\partial_s u - J(t,u) (\partial_t u - X_H(t,u)) =0,
\end{equation}
such that $a_1\leq \mathcal{A}(u(s,\cdot)) \leq a_2$ for any $s\in
\R$, is bounded in $L^{\infty}(\R \times \T,T^*M)$;

\item the set of solutions $u=(q,p)\in C^{\infty}(\R \times [0,1], T^*M)$ of
(\ref{nncr}) such that $q(s,0)=q_0$, $q(s,1)=q_1$, and $a_1\leq
  \mathcal{A}(u(s,\cdot)) \leq a_2$
for any $s\in \R$ is bounded in $L^{\infty}(\R \times [0,1],T^*M)$.

\end{enumerate}

\noindent Furthermore, if $r>2$ there exists $j_1=j_1(r)$ such that if
$\|J-\widehat{J}\|_{\infty} < j_1$ then for every
$a_1,a_2,a_3\in \R$ there holds:

\begin{enumerate}
\setcounter{enumi}{2}

\item the set of solutions 
\[
u=(q,p)\in C^{\infty}(]0,+\infty[ \times
\T, T^*M)\cap W^{1,r}(]0,1[\times \T,T^*M)\] 
of (\ref{nncr}) such that
$a_1\leq \mathcal{A}(u(s,\cdot)) \leq a_2$ for any $s\in [0,+\infty[$,
and 
\[
\|q(0,\cdot)\|_{W^{1-1/r,r}(\T,\R^N)}\leq a_3
\] 
is bounded in $L^{\infty}([0,+\infty[ \times \T,T^*M)$;

\item the set of solutions 
\[
u=(q,p)\in C^{\infty}(]0,+\infty[ \times
[0,1], T^*M)\cap W^{1,r}(]0,1[\times ]0,1[,T^*M)
\] 
of (\ref{nncr}) such that
$q(s,0)=q_0$, $q(s,1)=q_1$,
$a_1\leq \mathcal{A}(u(s,\cdot)) \leq a_2$ for any $s\in [0,+\infty[$,
and 
\[
\|q(0,\cdot)\|_{W^{1-1/r,r}(]0,1[,\R^N)}\leq a_3
\] 
is bounded in $L^{\infty}([0,+\infty[ \times [0,1],T^*M)$.

\end{enumerate}
\end{thm}

\begin{proof}
Using the above mentioned embedding, equation (\ref{nncr}) can be written as
\begin{equation}
\label{nnncr}
(\partial_s - J_0 \partial_t)u = (J-J_0)\partial_t u - JX_H(t,u).
\end{equation}
Let $\chi\in C^{\infty}(\R)$ be a function such
that $\chi=1$ on $[0,1]$, $\supp \chi\subset ]-1,2[$, $0\leq
\chi\leq 1$, and $|\chi^{\prime}|\leq 2$.

\medskip

{\sc Proof of} (i) {\sc and} (ii). Let $u=(q,p)$ be a solution meeting
the requirements of (i) or (ii).
Let $r>2$, $k\in \Z$ and let
$\bar{q}(t)=tq_1 + (1-t)q_0$ in case (ii), $\bar{q}(t)\equiv 0$ in
case (i). Set $v(s,t) = \chi(s-k) (q(s,t)-\bar{q}(t),p(s,t))$,
$(s,t)\in \R \times [0,1]$. By (\ref{nnncr}), $v$ satisfies
\[
(\partial_s - J_0 \partial_t) v = (J(t,u)-J_0) \partial_t v + \chi^{\prime}
(q-\bar{q},p) - \chi (0,\bar{q}^{\prime}) - \chi J(t,u) X_H(t,u),
\]
on either $\R \times \T$ - case (i) - or $\R \times [0,1]$ - case
(ii) - in which case we also have $v(s,0)=v(s,1)=0$ for every $s\in \R$.
Moreover $v$ has compact support, so by Proposition \ref{caldzyg}
\begin{equation}
\label{last}
\begin{split}
\| \nabla v \|_{L^r(\R \times ]0,1[)} \leq C_r \|(\partial_s - J_0
\partial_t) v\|_{L^r(\R \times ]0,1[)}  \leq C_r \|J-J_0\|_{\infty} \|
\partial_t v\|_{L^r(\R \times ]0,1[)} \\
+ 4\cdot 3^{\frac{1}{r}} C_r d + 2
\|p\|_{L^r(]k-1,k+2[ \times ]0,1[)} + C_r 3^{\frac{1}{r}}
\|\overline{q}^{\prime} \|_{\infty} + C_r
\|J\|_{\infty} \|X_H(\cdot,u)\|_{L^r(]k-1,k+2[ \times ]0,1[)},
\end{split}
\end{equation}
where $d=\max\set{|z|}{z\in M}$ and $C_r = C_r(\R \times \T)$ (resp.\
$C_r=C_r(\R \times ]0,1[)$).
If $\|J-\widehat{J}\|_{\infty}<+\infty$, $\|J\|_{\infty}$ is bounded, so
Lemma \ref{gradest} implies
\begin{equation}
\label{lr}
\|p\|_{L^r(]k-1,k+2[\times ]0,1[)} \leq
S_r \|p\|_{W^{1,2}(]k-1,k+2[ \times ]0,1[)} \leq S_r c
(3+(1+\sqrt{3})^2)^{\frac{1}{2}} \leq 4 S_r c,
\end{equation} 
where $S_r$ is the norm of the
continuous embedding 
\[
W^{1,2}(]0,3[ \times ]0,1[) \hookrightarrow L^r(]0,3[\times ]0,1[).
\] 
By (\ref{gs}) we also
have 
\begin{equation}
\label{nhg}
\begin{split}
\|X_H(\cdot,u)\|_{L^r(]k-1,k+2[ \times \T)} \leq h_4(3^{\frac{1}{r}} +
\|p\|^2_{L^{2r}(]k-1,k+2[\times ]0,1[)} ) \\ \leq  h_4(3^{\frac{1}{r}} +
S_{2r}^2 \|p\|^2_{W^{1,2} (]k-1,k+2[\times ]0,1[)} ) \leq
h_4(3^{\frac{1}{r}} + 16 S_{2r}^2 c^2). 
\end{split}
\end{equation}
If $\|J-J_0\|_{\infty} = \|J-\widehat{J}\|_{\infty} < 1/C_r$,
estimate (\ref{last}) together with (\ref{lr}) and (\ref{nhg}),
implies that $\nabla v$ is uniformly bounded in
$L^r(\R \times ]0,1[)$. Therefore,
$u$ is uniformly bounded in $W^{1,r}(]k,k+1[\times ]0,1[)$, and since
$r>2$, also in $L^{\infty} ([k,k+1]\times [0,1])$. Since $k\in \Z$ was
arbitrary we have a uniform bound for $u$ in $L^{\infty}(\R \times
[0,1])$. Therefore, statement (i) (resp.\ (ii)) holds with
\[
j_0 = \sup_{r\in ]2,+\infty[} 1/C_r(\R \times \T) \quad \mbox{(resp.\ }
j_0 = \sup_{r\in ]2,+\infty[} 1/C_r(\R \times ]0,1[) \mbox{ )}.
\]

\medskip

{\sc Proof of } (iii) {\sc and } (iv). The above argument for
$k\geq 1$ yields to a uniform bound for $u$ in
$L^{\infty}([1,+\infty[\times [0,1])$. There remains to find a uniform
bound for $u$ in $L^{\infty}([0,1]\times [0,1])$.

By the theory of Sobolev traces (see \cite{ada75} section 7.56), there
exists a number $b_r$ such that
every $f\in W^{1-1/r,r}(\{0\} \times ]0,1[)$ has an extension $\tilde{f}\in
W^{1,r}(]0,+\infty[\times ]0,1[)$ such that
\[
\| \tilde{f} \|_{W^{1,r}(]0,+\infty[\times ]0,1[)} \leq b_r
  \|f\|_{W^{1-1/r,r}(\{0\} \times ]0,1[)}.
\]
Therefore, there exists a map $\tilde{q} \in
W^{1,r}(]0,+\infty[\times ]0,1[,\R^N)$ such that $\tilde{q}(0,t)=q(0,t)$ and
\[
\|\tilde{q}\|_{W^{1,r}(]0,+\infty[ \times ]0,1[)} \leq b_r
  \|q(0,\cdot)\|_{W^{1-1/r,r}(]0,1[)} \leq b_r a_3.
\]
In the case (iii) we can assume $\tilde{q}$ to be 1-periodic in $t$,
while in the case (iv) we can assume that $\tilde{q}(s,0)=q_0$,
$\tilde{q}(s,1)=q_1$, for every $s\in [0,1]$. The map $w(s,t) =
\chi(s) (q(s,t) - \tilde{q}(s,t),p(s,t))$ satisfies
\[
(\partial_s - J_0 \partial_t) w = (J(t,u)-J_0) \partial_t w + \chi^{\prime}
(q-\tilde{q},p) - \chi (\partial_s \tilde{q}, \partial_t \tilde{q}) -
\chi J(t,u) X_H(t,u),
\]
Since
\[
\|\chi (\partial_s \tilde{q}, \partial_t \tilde{q})\|_{L^r(]0,+\infty[
  \times ]0,1[)}
  \leq  \|\tilde{q}\|_{W^{1,r}(]0,+\infty[ \times ]0,1[)} \leq 2b_r a_3,
\]
the same argument used above - involving Proposition \ref{caldzyg}
and the $W^{1,2}$ estimate of $p$ on $]0,2[\times
]0,1[$ provided by Lemma \ref{gradest} - allows to conclude.
\end{proof}

Once $L^{\infty}$ estimates are established, compactness
in $C^{\infty}_{\mathrm{loc}}$ follows by standard arguments. Here we
are interested in the following statement:

\begin{thm}
\label{comp}
Assume that $H$ satisfies (H1), (H2), and that $J$ satisfies
$\|J-\widehat{J}\|_{\infty} < j_0$. Then
for every $x^-,x^+\in \mathcal{P}_{\Lambda}(H)$ (resp.\ $x^-,x^+\in
\mathcal{P}_{\Omega} (H)$), the space
$\mathcal{M}_{\Lambda}(x^-,x^+)$ (respectively
$\mathcal{M}_{\Omega}(x^-,x^+)$) is pre-compact in
$C^{\infty}_{\mathrm{loc}} (\R\times \T,T^*M)$ (respectively in
$C^{\infty}_{\mathrm{loc}} (\R \times [0,1],T^*M)$).
\end{thm}

Indeed, solutions of (\ref{cr}) on $\R\times \T$ with bounded
action have uniform gradient bounds, because otherwise a
bubbling-off argument would produce a $J$-holomorphic sphere in
$T^* M$, which cannot exist because $\omega$ is exact. Solutions
of (\ref{cr}) on $\R \times [0,1]$ taking values in $T_{q_0}^* M$
for $t=0$ and in $T_{q_1}^* M$ for $t=1$,
and having bounded
action also have uniform gradient bounds: in this case the
bubbling-off argument could also produce a $J$-holomorphic disc
with boundary in either $T_{q_0}^* M$ or $T_{q_1}^* M$, which
cannot exist because the Liouville form $\theta$ vanishes on the
vertical subspaces.

Then elliptic bootstrap produces bounds on the derivatives
of every order. See for instance \cite{flo88d} or \cite{sal90} for
more details.

\begin{rem}
Notice that all the results of this section hold also by replacing (H2)
with the weaker condition (\ref{gs}), which could therefore replace
(H2) in the whole paper. However, condition (\ref{gs}) is somehow
unsatisfactory because it depends on the choice of the metric on $M$.
\end{rem}

\begin{rem}
Besides the conditions on $H$, the $L^{\infty}$ estimates for the
Floer equation require that $J$ belongs
to a suitable neighborhood of the set of almost complex structures on
$T^*M$ produced by metrics on $M$. It would be interesting to have a
better description of a set of $\omega$-compatible almost complex
structures for which the $L^{\infty}$ estimates hold.
\end{rem}

\subsection{Transversality}

Let $\mathcal{J}=\mathcal{J}(\langle \cdot,\cdot\rangle)$ be the set
of all $t$-dependent 1-periodic smooth $\omega$-compatible
almost complex structures on $T^*M$ such that
$\|J-\widehat{J}\|_{\infty} < +\infty$. The distance
\begin{eqnarray*}
\dist (J_1,J_2) = \|J_1-J_2\|_{\infty} +
\mathrm{dist}_{C^{\infty}_{\mathrm{loc}}} (J_1,J_2)
\end{eqnarray*}
makes $\mathcal{J}$ a complete metric space. Here
$\mathrm{dist}_{C^{\infty}_{\mathrm{loc}}}$ is the usual distance
\begin{eqnarray*}
\mathrm{dist}_{C^{\infty}_{\mathrm{loc}}} (J_1,J_2) = \sum_{r=1}^{\infty}
\sum_{\ell =0}^{\infty} 2^{-(r+\ell)} \frac{\|J_1 -
  J_2\|_{C^{\ell}(K_r)}}{1 + \|J_1 -J_2\|_{C^{\ell}(K_r)}}, \;
  \mbox{where } K_r = \set{(t,q,p)\in \T\times \T^*M}{|p|\leq r},
\end{eqnarray*}
inducing the $C^{\infty}_{\mathrm{loc}}$ topology. We denote by
$\mathcal{J}_{\mathrm{reg}}= \mathcal{J}_{\mathrm{reg}}(H)$
the subset of $\mathcal{J}$ consisting of those almost complex
structures $J$ for which the section
\[
\partial_{J,H} : \mathcal{B}(x^-,x^+) \rightarrow \mathcal{W}
(x^-,x^+)
\]
is transverse to the zero-section, for every $x^-,x^+\in
\mathcal{P}(H)$.  We recall that a subset of a topological space
is called residual if it contains a countable intersection of open
and dense sets. Baire theorem states that a residual subset of a
complete metric space is dense. The proof of the following result
is absolutely standard (see \cite{fhs96}):
\begin{thm}
\label{trans}
The set $\mathcal{J}_{\mathrm{reg}}(H)$ is residual in $\mathcal{J}$.
\end{thm}

\subsection{The Floer complex}

Let $H$ be a Hamiltonian on $\T \times T^*M$ satisfying (H0), (H1),
and (H2). Denote by $CF_{\Lambda,k}(H)$ (resp.\ $CF_{\Omega,k}(H)$)
the free Abelian group generated by the elements $x$ of
$\mathcal{P}_{\Lambda}(H)$ (resp.\ $\mathcal{P}_{\Omega}(H)$) with
Maslov index $\mu_{\Lambda}(x) = k$ (resp.\ $\mu_{\Omega}(x) =
k$). Notice that these groups need not be finitely generated. Since
the discussion will present no differences in the $\Lambda$ and in the
$\Omega$ case, we will omit the subscripts and deal with both
situations at the same time.

Let $j_0$ be the positive number given by Theorem \ref{c0est}.
By Theorem \ref{trans}, the set
\[
\mathcal{J}_{j_0,\mathrm{reg}}(H):=
\set{J\in \mathcal{J}_{\mathrm{reg}} (H)}{\|J - \widehat{J}\|_{\infty}
  < j_0}
\]
is non-empty. Let us fix some $J\in
\mathcal{J}_{j_0,\mathrm{reg}}(H)$. If $x,y\in \mathcal{P}(H)$
have index difference $\mu(x) - \mu(y)=1$, Theorem \ref{indice}
and transversality imply that $\mathcal{M}(x,y)$ is a
one-dimensional manifold. The compactness stated in Theorem
\ref{comp} and transversality imply that $\mathcal{M}(x,y)$
consists of finitely many lines. Then we can define the integer
$n(x,y)$ to be
\[
n(x,y) := \sum_{[u]\in \mathcal{M}(x,y)/\R} \epsilon([u]),
\]
the numbers $\epsilon([u])\in \{-1,+1\}$ having been defined in
section \ref{scoh}. The homomorphism
\[
\partial_k = \partial_k(H,J) : CF_k(H) \rightarrow CF_{k-1} (H)
\]
is defined in terms of the generators by
\[
\partial_k x = \sum_{\substack{y\in \mathcal{P}(H) \\ \mu(y) = k-1}}
n(x,y)  y, \quad \forall x\in \mathcal{P}(H), \; \mu(x) = k.
\]
Indeed, the above sum contains finitely many terms thanks to
(\ref{decr}) and Lemma \ref{lem0}.
A standard gluing argument shows that $\partial_{k-1}
\circ \partial_k =0$, so $\{CF_*(H),\partial_*(H,J)\}$ is a complex of
free Abelian groups, called the {\em Floer complex of $(H,J)$}.
The homology of such a complex is called the {\em Floer homology of
  $(H,J)$}:
\[
HF_k(H,J) = \frac{\ker(\partial_k:CF_k(H,J) \rightarrow
  CF_{k-1}(H,J))}{ \mathrm{ran} (\partial_{k+1}:CF_{k+1}(H,J) \rightarrow
  CF_k(H,J))}.
\]
The Floer complex splits into subcomplexes, one for each conjugacy
class of $\pi_1(M)$ in the $\Lambda$ case, one
for each element of $\pi_1(M)$ in the $\Omega$ case.
Moreover, the Floer complex has an $\R$-filtration
defined by the action functional: if $CF^a_k(H)$ denotes the subgroup
of $CF_k(H)$ generated by the $x\in \mathcal{P}(H)$ such that
$\mathcal{A}(x) <a$, the boundary operator $\partial_k$ maps $CF_k^a(H)$
into $CF_{k-1}^a(H)$, so $\{CF^a_*(H),\partial_*(H,J)\}$ is a
subcomplex. By Lemma \ref{lem0}, such a subcomplex is finitely generated.

Changing the orientation data (namely, the preferred unitary
trivializations of $x^*(TT^*M)$, for $x\in \mathcal{P}(H)$, and
the coherent orientation for $\Sigma$), we obtain an isomorphic
chain complex, the isomorphism being of the special form
\[
x\mapsto \sigma(x) x, \quad \forall x\in \mathcal{P}(H),
\]
where $\sigma(x)\in \{-1,+1\}$. What is less trivial is that a
different choice of the almost complex structure $J$ - an
operation which changes the Floer equation, and thus its solution
spaces - produces isomorphic Floer complexes, as the next result
shows:

\begin{thm} 
\label{chj}
If $J^0,J^1\in \mathcal{J}_{j_0,\mathrm{reg}}(H)$, there is an isomorphism of
complexes
\[
\phi_{01} : \{CF_*(H),\partial_*(H,J_0)\} \rightarrow
\{CF_*(H),\partial_*(H,J_1)\}, \quad x \mapsto \sum_{\substack{y\in
    \mathcal{P}(H)\\ \mu(y) = \mu(x)}} n_{01}(x,y) y,
\]
such that $n_{01}(x,x)=1$ and $n_{01}(x,y)=0$ if $\mathcal{A}(x) \leq
\mathcal{A}(y)$ and $x\neq y$, or if $x$ and $y$ are not homotopic.
Such an isomorphism is uniquely
defined up to chain homotopy. If $J^2$ is a third element of
$\mathcal{J}_{j_0,\mathrm{reg}}(H)$, the isomorphisms
$\phi_{12} \circ \phi_{01}$ and
$\phi_{02}$ are chain homotopic.
\end{thm}

In particular, the isomorphism $\phi_{01}$ preserves the
$\R$-filtration, and it is compatible with the splitting of the
Floer complex determined by the structure of $\pi_1(M)$. The above
result is due to Cornea and Ranicki, \cite{cr03} (in the case of
Floer homology for a class of compact symplectic manifolds). 
See \cite{oh97}, Lemma 6.3, for an earlier application of the same
argument.
Here we just sketch the proof.

\begin{proof}
Using the fact that the space $\mathcal{J}_{j_0}$, the
$L^{\infty}$-ball of $\mathcal{J}$
centered in $\widehat{J}$ of radius $j_0$, is contractible,
one can find a homotopy $(J_s)_{s\in \R}$ in $\mathcal{J}_{j_0}$ such that
$J_s=J^0$ for $s\leq 0$ and $J_s=J^1$ for $s\geq 1$, such that
counting solutions of
\begin{equation}
\label{hmtpy}
\partial_s u - J_s (t,u)(\partial_t u - X_H(t,u)) = 0
\end{equation}
between $x,y\in \mathcal{P}(H)$ of the same Maslov index, defines a
chain map
\[
\phi_{01}:\{CF_*(H),\partial_*(H,J_0)\} \rightarrow
\{CF_*(H),\partial_*(H,J_1)\}.
\]
If $u$ solves (\ref{hmtpy}),
\[
\frac{d}{ds} \mathcal{A}(u(s))=- \int_0^1 |\partial_s
u(s,t)|^2_{J_{s,t}}\, dt,
\]
so the only solutions of (\ref{hmtpy}) connecting a curve
$x$ with a curve $y$ with $\mathcal{A}(y)\geq \mathcal{A}(x)$ are
the stationary ones. Notice that transversality holds automatically at
stationary solutions $u(s,t)=x(t)$. 
Indeed, linearization along such a
solution yields to an operator of the form
\begin{equation}
\label{op}
\partial_s - (J_s \partial_t + S(s,t)) 
\end{equation}
where the self-adjoint operator $J_s \partial_t + S(s,\cdot)$ depends
on $s$, but represents always the same quadratic form - the second
differential of $\mathcal{A}$ at $x$ - with respect to
inner products varying with $s$. In this case, the operator (\ref{op})
is easily shown to be invertible.    

We conclude that the coefficients $n_{01}(x,y)$
satisfy the required assumptions. This means that, if we order the
elements of $\mathcal{P}(H)$ - the generators of $CF_*(H)$ - by
increasing action, $\phi_{01}$ is represented by a lower
triangular matrix with diagonal entries equal 1, so it is an
isomorphism. The other statements can be proved by introducing a
homotopy of homotopies. \end{proof}

Therefore, we can consider the Floer homology $HF_*(H) =
HF_*(H,J)$ as independent of $J$. A different choice of the
Hamiltonian, instead, produces chain homotopic complexes:

\begin{thm} 
\label{chh}
Let $H_0,H_1$ be Hamiltonians on $\T \times T^*M$ satisfying (H0),
(H1), and (H2), and let $J\in
\mathcal{J}_{j_0,\mathrm{reg}} (H_0) \cap \mathcal{J}_{j_0,\mathrm{reg}}(H_1)$.
Then there is a homotopy equivalence
\[
\psi_{01} : \{CF_*(H_0),\partial_*(H_0,J)\} \rightarrow
\{CF_*(H_1),\partial_*(H_1,J)\},
\]
uniquely determined up to chain homotopy.
If moreover $H_2$ is a third Hamiltonian satisfying the same
conditions, and such that $J \in \mathcal{J}_{j_0,\mathrm{reg}}(H_2)$, then
the chain maps $\psi_{12} \circ \psi_{01}$ and
$\psi_{02}$ are chain homotopic.
\end{thm}

In particular, $HF_*(H_0)\cong HF_*(H_1)$.
This result can be proved 
using the standard homotopy argument from Floer theory: one introduces
an $s$-dependent Hamiltonian $H:\R \times \T \times T^*M \rightarrow
\R$ such that $H(s,\cdot,\cdot)=H_0$ for $s\leq 0$ and
$H(s,\cdot,\cdot) = H_1$ for $s\geq 1$, and defines the chain map
$\psi_{01}$ by considering the solutions of the equation
\[
\partial_s u - J(t,u) (\partial_t u - X_{H} (s,t,u)) = 0,
\]
connecting elements of $\mathcal{P}(H_0)$ and $\mathcal{P}(H_1)$. The
only delicate point is the $L^{\infty}$ estimate for the solutions of
the above problem. This estimate can be achieved by adapting the arguments of
section \ref{cs}, provided that the Hamiltonians $H_0$ and $H_1$ are
close enough. Then the isomorphism between the Floer complexes of two
arbitrary Hamiltonians satisfying (H0), (H1) and (H2) can be constructed by
composing a finite number of isomorphisms. Details are contained in
the next section. 

\subsection{$\mathbf{L^{\infty}}$ estimates for homotopies}

Let $H_0$ and $H_1$ be Hamiltonians on $\T\times T^*M$ satisfying (H1)
and (H2). Up to choosing a smaller $h_0$ and larger $h_1$, $h_2$, we
may assume that $H_0$ and $H_1$ satisfy conditions (H1) and (H2) with
the same constants $h_0$, $h_1$, $h_2$. Let $\chi:\R\rightarrow [0,1]$
be a smooth function such that $\chi(s)=0$ for $s\leq 0$, $\chi(s)=1$
for $s\geq 1$, $0\leq \chi^{\prime} \leq 2$, and set
\begin{equation}
\label{sdh}
H:\R \times \T \times T^*M \rightarrow \R, \quad H(s,t,x) = \chi(s)
H_1(t,x) + (1-\chi(s)) H_0(t,x).
\end{equation}
Every Hamiltonian $H_s:= H(s,\cdot,\cdot)$ satisfies (H1) and (H2)
with constants $h_0$, $h_1$, $h_2$. 

We are going to show that if $H_0$ and $H_1$ are close enough, then
Lemma \ref{gradest} extends to the $s$-dependent Hamiltonian $H$. 

\begin{lem}
\label{ngradest}   
Assume that the $t$-dependent 1-periodic almost complex structure $J$
on $T^*M$ satisfies $\|J\|_{\infty}<+\infty$. 
There exists a positive number
$\epsilon=\epsilon(h_0, \|J\|_{\infty})$ such that
if $H_0$ and $H_1$ satisfy conditions (H1) and (H2) with constants
$h_0$, $h_1$, $h_2$, and
\begin{equation}
\label{oh1.5}
 |H_1(t,q,p)-H_0(t,q,p)| \leq h + \epsilon |p|^2 \quad \forall
 (t,q,p)\in \T \times T^*M,
\end{equation}
for some $h\geq 0$, then the following a priori estimate holds. For
every pair of real numbers $a_1,a_2$ there exists a number $c$ such
that for every $u=(q,p)\in C^{\infty}(\R \times [0,1],T^*M)$ solving
\begin{equation}
\label{oh2}
\partial_s u - J(t,u) (\partial_t u - X_H(s,t,u)) = 0,
\end{equation}
and such that
\begin{equation}
\label{oh3}
\mathcal{A}_{H_0}(u(s,\cdot)) \leq a_2 \;\; \forall s\leq
0, \quad \mathcal{A}_{H_1}(u(s,\cdot)) \geq a_1 \;\; \forall s\geq
1,
\end{equation}
there holds
\[
\|p\|_{L^2(I\times ]0,1[)} \leq c|I|^{\frac{1}{2}}, \quad
\|\nabla p\|_{L^2(I\times ]0,1[)} \leq c(1+|I|^{\frac{1}{2}}),
\]
for every interval $I\subset \R$. 
\end{lem}

\begin{proof}
Denote by $\mathcal{U}$ the set of solutions $u=(q,p)$ of (\ref{oh2})
satisfying the action estimate (\ref{oh3}). 

\medskip

{\sc Claim 0.} For every $u\in \mathcal{U}$ and every $s\in \R$
there holds
\begin{equation}
\label{oh4}
\mathcal{A}_{H_s}(u(s,\cdot)) \leq a_2 + 2h + 2\epsilon
\|p\|_{L^2(]0,1[\times ]0,1[)}^2.
\end{equation}

\medskip

The function $s\mapsto \mathcal{A}_{H_s}(u(s,\cdot))$ is decreasing on
$]-\infty,0]$ and on $[1,+\infty[$, so
\begin{eqnarray*}
\mathcal{A}_{H_s}(u(s,\cdot)) = \mathcal{A}_{H_0}(u(s,\cdot)) \leq a_2
\quad \forall s\leq 0, \\ \mathcal{A}_{H_s}(u(s,\cdot)) =
\mathcal{A}_{H_1}(u(s,\cdot)) \leq \mathcal{A}_{H_1}(u(1,\cdot)) 
\quad \forall s\geq 1,
\end{eqnarray*}
and it is enough to prove (\ref{oh4}) for $s\in [0,1]$. In this case,
by (\ref{oh1.5}), (\ref{oh2}), and (\ref{oh3}),
\begin{eqnarray*}
\mathcal{A}_{H_s}(u(s,\cdot)) = \mathcal{A}_{H_0}(u(0,\cdot)) +
\int_0^s \frac{d}{d\sigma} \bigl(
  \mathcal{A}_{H_{\sigma}}(u(\sigma,\cdot)) \bigr) \, d\sigma \\ 
= \mathcal{A}_{H_0}(u(0,\cdot)) + \int_0^s \left(
  d\mathcal{A}_{H_{\sigma}}(u(\sigma,\cdot)) [\partial_s
u(\sigma,\cdot)] - \int_0^1 \partial_s H(\sigma,t,u)\, dt \right) \,
d\sigma \\
=  \mathcal{A}_{H_0}(u(0,\cdot)) - \int_0^s \int_0^1 |\partial_s
u|_{J_t}^2 \, dt \, d\sigma  + \int_0^s \chi^{\prime}(\sigma)
\int_0^1 (H_0(t,u)-H_1(t,u))\, dt \, d\sigma \\
\leq a_2 + \|\chi^{\prime}\|_{\infty} \left( hs + \epsilon \int_0^s
  \|p(\sigma,\cdot) \|_{L^2(]0,1[)}^2 \, d\sigma \right) \leq a_2 + 2h
+ 2\epsilon \|p\|_{L^2(]0,1[\times ]0,1[)}^2 ,
\end{eqnarray*}
proving the claim.

\medskip

{\sc Claim 1.} For every $u\in \mathcal{U}$ there holds
\[
\| \partial_s u\|_{L^2(\R \times ]0,1[)}^2 \leq \|J\|_{\infty}^2
\left( a_2 - a_1 + 2h + 2\epsilon \|p\|^2_{L^2(]0,1[ \times ]0,1[)}
\right).
\]

\medskip

Indeed, by (\ref{oh1.5}), (\ref{oh2}), (\ref{oh3}),
\begin{eqnarray*}
\| \partial_s u\|_{L^2(\R \times ]0,1[)}^2 \leq \|J^{-1}\|_{\infty}^2
\int_{-\infty}^{+\infty} \int_0^1 |\partial_s u|_{J_t}^2 \, dt\, ds 
 = - \|-J\|_{\infty}^2 \int_{-\infty}^{+\infty} d\mathcal{A}_{H_s}
(u(s,\cdot)) [\partial_s u(s,\cdot)] \, ds \\
= - \|J\|_{\infty}^2 \int_{-\infty}^{+\infty} \left( \frac{d}{ds}
  \bigl( \mathcal{A}_{H_s}(u(s,\cdot)) \bigr) + \int_0^1 \partial_s
  H(s,t,u)\, dt \right) \, ds \\
= - \|J\|_{\infty}^2 \Bigl( \lim_{s\rightarrow +\infty}
  \mathcal{A}_{H_s}(u(s,\cdot)) - \lim_{s\rightarrow -\infty}
  \mathcal{A}_{H_s}(u(s,\cdot))  + \int_0^1 \chi^{\prime}(s)\int_0^1
  \Bigl( H_1(t,u) - H_0(t,u) \Bigr) \, dt \, ds \Bigr) \\
\leq \|J\|_{\infty}^2 \left( a_2 - a_1 + 2h + 2\epsilon
  \|p\|_{L^2(]0,1[\times ]0,1[)}^2 \right),
\end{eqnarray*}
as claimed.

\medskip

{\sc Claim 2.} For every $u\in \mathcal{U}$ and every $s\in \R$ 
there holds
\begin{equation}
\label{oh5}
\frac{h_0}{2} \|p(s,\cdot)\|_{L^2(]0,1[)}^2 - (h_1 + a_2 + 2h)  \leq
2\epsilon \|p\|_{L^2(]0,1[\times ]0,1[)}^2 + \frac{1}{2h_0}
\|J\|_{\infty}^2 \|\partial_s u(s,\cdot)\|_{L^2(]0,1[)}^2.
\end{equation}

\medskip

Since $H_s$ satisfies condition (H1), arguing as in the proof of Lemma
\ref{gradest}, Claim 2, we obtain
\[
\theta(\partial_t u) - H(s,t,u) \geq h_0 |p|^2 - h_1 - \|J\|_{\infty}
|p|\, |\partial_s u|.
\]
Then by Claim 0,
\begin{eqnarray*}
a_2 + 2h + 2\epsilon \|p\|_{L^2(]0,1[\times ]0,1[)}^2 \geq
\mathcal{A}_{H_s}(u(s,\cdot)) = \int_0^1 \left( \theta(\partial_t u) -
  H(s,t,u) \right) \, dt \\ 
\geq h_0 \|p(s,\cdot)\|_{L^2(]0,1[)}^2 - h_1 - \|J\|_{\infty}
\|p(s,\cdot)\|_{L^2(]0,1[)} \|\partial_s u(s,\cdot)\|_{L^2(]0,1[)} \\
\geq \frac{h_0}{2} \|p(s,\cdot)\|_{L^2(]0,1[)}^2 - h_1 -
\frac{1}{2h_0} \|J\|_{\infty}^2 \|\partial_s
u(s,\cdot)\|_{L^2(]0,1[)}^2,
\end{eqnarray*}
which is equivalent to (\ref{oh5}).

\medskip

Integrating the inequality (\ref{oh5}) over $]0,1[$ and using Claim 1,
we get 
\begin{eqnarray*}
\frac{h_0}{2} \|p\|_{L^2(]0,1[\times ]0,1[)}^2 - (h_1 + a_2 + 2h) 
\leq 2\epsilon  \|p\|_{L^2(]0,1[\times ]0,1[)}^2 + \frac{1}{2h_0}
\|J\|_{\infty}^2 \|\partial_s u\|_{L^2(]0,1[\times ]0,1[)}^2 \\
\leq 2\epsilon \|p\|_{L^2(]0,1[\times ]0,1[)}^2 +  \frac{1}{2h_0}
\|J\|_{\infty}^4 \left( a_2 - a_1 + 2h +
  2\epsilon  \|p\|_{L^2(]0,1[\times ]0,1[)}^2 \right),
\end{eqnarray*}
or equivalently,
\begin{eqnarray*}
\left( \frac{h_0}{2} - 2\epsilon \Bigl( 1 + \frac{1}{2h_0}
  \|J\|_{\infty}^4 \Bigr) \right)
\|p\|_{L^2(]0,1[\times ]0,1[)}^2 \leq h_1 + a_2 + 2h + \frac{1}{2h_0}
\|J\|_{\infty}^4  (a_2 - a_1 + 2h).
\end{eqnarray*}
Therefore, if $\epsilon=\epsilon(h_0, \|J\|_{\infty})$ satisfies 
\begin{equation}
\label{epsi}
\epsilon < \frac{h_0}{4} \left( 1 + \frac{1}{2h_0}\|J\|_{\infty}^4 
\right)^{-1} ,
\end{equation}
we deduce that $\|p\|_{L^2(]0,1[\times ]0,1[)}$ is uniformly bounded,
for $u=(q,p)\in \mathcal{U}$. Hence, Claim 1 and Claim 2 can be
improved, producing the following estimates.

\medskip

{\sc Claim 1$^{\prime}$.} There exists $c_1$ such that $\|\partial_s u\|_{L^2(\R
  \times ]0,1[)} \leq c_1$ for every $u\in \mathcal{U}$.

\medskip

{\sc Claim 2$^{\prime}$.} There exists $c_2$ such that
$\|p(s,\cdot)\|_{L^2(]0,1[)} \leq c_2 (1+ \|\partial_s
u(s,\cdot)\|_{L^2(]0,1[)} )$ for every $u\in \mathcal{U}$ and every
$s\in \R$.

\medskip

These are exactly the first two claims in the proof of Lemma
\ref{gradest}. The remaining part of the proof of that lemma extends
to the case of the $s$-dependent Hamiltonian without any change.
\end{proof}

The above lemma and the Calderon-Zygmund estimates of Proposition
\ref{caldzyg} imply the following $L^{\infty}$ estimates. The proof is
identical to the proof of Theorem \ref{c0est}.

\begin{thm}
\label{sc0est}
Assume that the Hamiltonians $H_0$ and $H_1$ satisfy (H1), (H2), and
(\ref{oh1.5}) with $\epsilon = \epsilon(h_0, \|J\|_{\infty})$ 
small enough as in (\ref{epsi}). 
Let $H$ be the $s$-dependent Hamiltonian defined in
(\ref{sdh}). Assume that the $t$-dependent 1-periodic almost complex
structure $J$ on $T^*M$ satisfies $\|J-\hat{J}\|_{\infty} < j_0$,
where $j_0$ is given by Theorem \ref{c0est}. Then for every
$a_1,a_2\in \R$ there holds:
\begin{enumerate}

\item the set of solutions $u=(q,p)\in C^{\infty}(\R \times \T,T^*M)$
  of
\begin{equation}
\label{aglo}
\partial_s u - J(t,u) (\partial_t - X_H(s,t,u)) = 0,
\end{equation}
such that $\mathcal{A}_{H_0}(u(s,\cdot))\leq a_2$ for every $s\leq 0$
and $\mathcal{A}_{H_1} (u(s,\cdot)) \geq a_1$ for every $s\geq 1$, is
bounded in $L^{\infty} (\R\times \T, T^*M)$;

\item the set of solutions $u=(q,p)\in C^{\infty}(\R \times [0,1],T^*M)$
  of (\ref{aglo}) such that $q(s,0)=q_0$, $q(s,1)=q_1$,
 $\mathcal{A}_{H_0}(u(s,\cdot))\leq a_2$ for every $s\leq 0$
and $\mathcal{A}_{H_1} (u(s,\cdot))$ $\geq a_1$ for every $s\geq 1$, is
bounded in $L^{\infty} (\R\times [0,1], T^*M)$.
\end{enumerate}
\end{thm}   

Let us conclude this section by sketching the proof of Theorem
\ref{chh}. Let $H_0$ and $H_1$ be Hamiltonians satisfying (H0), (H1),
and (H2). We may assume that (H1) and (H2) hold with the
same constants $h_0,h_1,h_2$. By the second condition of (H2) and the
compactness of $M$, there exists $h_3\geq 0$ such that
\begin{equation}
\label{lah3}
|H_0(t,q,p)|\leq h_3 ( 1+ |p|^2), \quad |H_1(t,q,p)|\leq h_3 ( 1+
|p|^2).
\end{equation}
Let $\lambda \in [0,1]$ and set $H_\lambda = \lambda H_1 + (1-\lambda)
H_0$. If $\lambda_0,\lambda_1\in [0,1]$, (\ref{lah3}) implies
\[
|H_{\lambda_1}(t,q,p) - H_{\lambda_0}(t,q,p)| = |(\lambda_1 -
\lambda_0) (H_1-H_0)(t,q,p)| \leq 2h_3 |\lambda_1 - \lambda_0|
(1+|p|^2).
\]
So, if $|\lambda_1 - \lambda_0|\leq  \epsilon/(2h_3)$ the Hamiltonians
$H_{\lambda_0}$ and $H_{\lambda_1}$ satisfy the assumptions of Theorem
\ref{sc0est}. In this case, the moduli spaces of solutions $u$ of equation
(\ref{aglo}) satisfying $u(-\infty,\cdot) \in
\mathcal{P}(H_{\lambda_0})$ and $u(+\infty,\cdot) \in
\mathcal{P}(H_{\lambda_1})$ can be used to define a chain map
\[
\psi_{\lambda_0 \lambda_1} : \{ CF_*(H_{\lambda_0},
\partial_*(H_{\lambda_0},J)\} \rightarrow 
\{ CF_*(H_{\lambda_1},\partial_*(H_{\lambda_1},J)\}.
\]
By the usual gluing argument, $\psi_{\lambda_1 \lambda_0}$ is a chain
homotopy inverse of $\psi_{\lambda_0 \lambda_1}$, which thus induces
an isomorphism at the homology level. The chain homotopy equivalence  
 \[
\psi_{01} : \{ CF_*(H_{0},
\partial_*(H_{0},J)\} \rightarrow 
\{ CF_*(H_{1},\partial_*(H_{1},J)\},
\] 
can then be defined as the composition
\begin{eqnarray*}
\psi_{01} = \psi_{\lambda_{k-1} \lambda_k}\circ  \dots \circ
\psi_{\lambda_0
  \lambda_1} \quad \mbox{where } 0=\lambda_0 < \lambda_1 < \dots <
\lambda_k =1, \\ |\lambda_{j} - \lambda_{j-1}|\leq
\frac{\epsilon}{2h_3} \;\; \forall j=1,\dots, k.
\end{eqnarray*}
Standard gluing arguments imply that the chain homotopy class of
$\psi_{01}$ does not depend on the choices we have made, and that
$\psi_{12} \circ \psi_{01}$ is chain homotopic to $\psi_{02}$,
concluding the proof of Theorem \ref{chh}. 

\section{The Morse complex of the Lagrangian action functional}

\subsection{Lagrangian dynamical systems}
\label{lds}

Let $M$ be a connected compact smooth manifold, the configuration space of a
Lagrangian dynamical system, assumed to be one-periodic in time.
Points in the tangent bundle $TM$ will be denoted
by $(q,v)$, with $q\in M$, $v\in T_q M$.
We will denote by $\tau:TM
\rightarrow M$ the standard projection, and by $T_{(q,v)}^v TM$ the vertical
subspace $\ker D\tau (q,v) \cong T_q M$ of $T_{(q,v)}TM$.
The Lagrangian $L: \T \times
TM \rightarrow \R$, $\T = \R/\Z$, will be a smooth function satisfying:

\begin{description}

\item[(L1)] there exists $\ell_0>0$ such that
\[
\nabla_{vv} L(t,q,v) \geq \ell_0 I
\]
for every $(t,q,v)\in \T \times TM$;

\item[(L2)] there exists $\ell_1\geq 0$ such that
\[
|\nabla_{vv} L(t,q,v)| \leq \ell_1, \quad |\nabla_{qv} L(t,q,v)| \leq \ell_1
 (1+|v|), \quad |\nabla_{qq} L(t,q,v) |\leq \ell_1(1+|v|^2)
\]
for every $(t,q,v)\in \T \times TM$.

\end{description}

Here we have fixed a Riemannian metric $\langle \cdot ,\cdot \rangle$
on $M$, with corresponding norm $|\cdot|$, and
$\nabla_{vv}$, $\nabla_{qv}$, $\nabla_{qq}$ denote the components of
the Hessian in the splitting of $TTM$ into the vertical and horizontal
part, given by the corresponding Levi-Civita connection. It
is easily seen that the above conditions do not depend on the choice
of the Riemannian metric. Physical Lagrangians of the form
\[
L(t,q,v) = \frac{1}{2} |T(t,q) v - A(t,q)|^2 - V(t,q)
\]
satisfy conditions (L1) and (L2), provided that the symmetric tensor
$T^*T$ is everywhere positive.

The strong convexity assumption\footnote{Here it would be enough
to
  assume that the map sending $v$ into the restriction of $dL(t,q,v)$
  to the vertical subspace is a
  diffeomorphism from $T_q M$ onto $T_q^*M$. The strict convexity in
  the $v$ variables will be important in order to guarantee the
  Palais-Smale condition (see Proposition \ref{palsma}).}
(L1) implies that $L$ defines a smooth
vector field $Y_L$ on $TM$. Indeed, we can define a $1$-periodic
Hamiltonian on $T^* M$ by means of the Legendre transform (see for
  instance \cite{man91}):
\begin{equation}
\label{leg}
H(t,q,p) = \max_{v\in T_q M} \left( p[v] - L(t,q,v) \right) =
p[v(t,q,p)] - L(t,q,v(t,q,p)),
\end{equation}
for every $(t,q,p) \in \T \times T^*M$, 
where the map $v$ is a component of the fiber-preserving
diffeomorphism
\[
\mathfrak{L}_L^{-1} : \T \times T^* M \rightarrow \T
\times TM, \quad (t,q,p) \mapsto (t,q,v(t,q,p)),
\]
the inverse of
\[
\mathfrak{L}_L : \T \times TM \rightarrow \T \times T^* M, \quad
(t,q,v) \mapsto \left(t,q,dL(t,q,v)|_{T_{(q,v)}^v TM}\right).
\]
Actually,
\begin{equation}
\label{equ}
H(t,q,p) = p[v] - L(t,q,v) \quad \mbox{if and only if} \quad p =
dL(t,q,v)|_{T^v_{(q,v)} TM}.
\end{equation}
The Hamiltonian $H$ and the canonical symplectic form $\omega$ on $T^*
M$ define the 1-periodic Hamiltonian vector field $X_H$ on $T^* M$,
and $Y_L$ is defined to be the pull-back of $X_H$ by the
diffeomorphism $\mathfrak{L}_{L}$. Equivalently, $Y_L$ is the
1-periodic Hamiltonian vector field on $TM$ determined by the
symplectic form $\mathfrak{L}_L^* \omega$ and by the Hamiltonian
$H\circ \mathfrak{L}_L$.

\begin{rem} Notice that if the Lagrangian has the form $L(t,q,v) = 1/2
  |v|^2 - V(t,q)$, then $\mathfrak{L}_L$ is the identity
  on $\T$ times the isomorphism $TM\rightarrow T^*M$ induced by the
  metric $\langle \cdot,\cdot \rangle$, and $H(t,q,p) = 1/2
  |p|^2 + V(t,q)$. In general, however, the restriction of
  $\mathfrak{L}_L$ to the fibers need not be linear.
\end{rem}

The integral curves $y:]a,b[ \rightarrow TM$ of the vector field
$Y_L$ are of the form $y(t)=(q(t),\dot{q}(t))$, where $q:]a,b[
\rightarrow M$ solves the second order ODE
\begin{equation}
\label{ls}
\nabla_t \left( \nabla_v L(t,q(t),\dot{q}(t)) \right) = \nabla_q
L(t,q(t),\dot{q}(t)).
\end{equation}
Here $\nabla_t$ denotes the covariant derivation along $q$, and
$\nabla_v$, $\nabla_q$ denote the vertical and the horizontal part of
the gradient of $L$.

We will be interested in the set $\mathcal{P}_{\Lambda}(L)$ of
1-periodic solutions of (\ref{ls}), and in the set
$\mathcal{P}_{\Omega} (L,q_0,q_1)$ of solutions $q:[0,1]\rightarrow M$
of (\ref{ls}) such that $q(0)=q_0$ and $q(1)=q_1$, for two fixed
points $q_0,q_1\in M$. In each of these cases we shall make one of the
following non-degeneracy assumptions:

\begin{description}

\item[(L0)$_{\mathbf{\Lambda}}$] every solution $q\in
\mathcal{P}_{\Lambda}(L)$ is
  non-degenerate, meaning that the differential of the time-one
  integral map of $Y_L$ at $q(0)$ does not have the eigenvalue 1;

\item[(L0)$_{\mathbf{\Omega}}$] every solution $q\in
  \mathcal{P}_{\Omega}(L,q_0,q_1)$ is non-degenerate, meaning that the
  differential of the time-one integral map of $Y_L$,
  \[
T_{(q(0),\dot{q}(0))} TM \rightarrow T_{(q(1),\dot{q}(1))} TM,
\]
 maps
  the vertical subspace at $(q(0),\dot{q}(0))$ into a subspace
  having intersection $(0)$ with the vertical subspace at
  $(q(1),\dot{q}(1))$.

\end{description}

These conditions can be stated in an equivalent way in terms of
the Jacobi vector fields along the solution $q$: (L0)$_{\Lambda}$
requires that there are no 1-periodic Jacobi vector fields, while
(L0)$_{\Omega}$ requires that there are no Jacobi vector fields
vanishing for $t=0$ and for $t=1$.

The Legendre transform $\mathfrak{L}_L$ provides us with a one-to-one
correspondence between the set of solutions
$\mathcal{P}_{\Lambda} (L)$ (resp.\ $\mathcal{P}_{\Omega} (L)$) of
the Lagrangian system and the set of solutions
$\mathcal{P}_{\Lambda}(H)$ (resp.\ $\mathcal{P}_{\Omega} (H)$) of the
Hamiltonian system. The non-degeneracy condition (L0) is equivalent to
its counterpart (H0).

\subsection{The variational setting}
Denote by $\Lambda^1(M)$ the
space of all loops $q:\T \rightarrow M$ of Sobolev class $W^{1,2}$, and by
$\Omega^1(M,q_0,q_1)$ the space of all paths $q:[0,1]\rightarrow M$ of
Sobolev class $W^{1,2}$ such that $q(0)=q_0$ and $q(1)=q_1$. These
spaces have Hilbert manifold structures (see \cite{kli82}
for this and for the other results cited in this section).
The tangent space of
$\Lambda^1(M)$ at $q$ is identified with the space of 1-periodic
$W^{1,2}$ tangent vector fields along $q$, while the tangent space of
$\Omega^1(M)$ at $q$ is identified with the space of
$W^{1,2}$ tangent vector fields along $q$ vanishing for $t=0$ and for
$t=1$.

The action functional
\[
\mathcal{E} (q) = \mathcal{E}_L (q) :=
\int_0^1 L(t,q(t),\dot{q}(t))\, dt
\]
is smooth on $\Lambda^1(M)$ and on $\Omega^1(M,q_0,q_1)$. Its
restrictions to these manifolds will be denoted by
$\mathcal{E}_{\Lambda}$ and $\mathcal{E}_{\Omega}$. The critical
points of $\mathcal{E}_{\Lambda}$ are the elements of
$\mathcal{P}_{\Lambda}(L)$, while the critical points of
$\mathcal{E}_{\Omega}$ are the elements of
$\mathcal{P}_{\Omega}(L,q_0,q_1)$. Condition (L0)$_{\Lambda}$
(resp.\ (L0)$_{\Omega}$) is equivalent to the fact that all the
critical points of $\mathcal{E}_{\Lambda}$ (resp.\ of
$\mathcal{E}_{\Omega}$) are non-degenerate. Moreover, condition
(L1) implies that all these critical points have finite Morse
indices, denoted by $m_{\Lambda}(q)$ and $m_{\Omega}(q)$. The
proof of the following result, essentially due to Duistermaat
\cite{dui76}, can be found in \cite{web02}, Theorem 1.2, for the
case of periodic orbits, and in \cite{rs95}, Proposition 6.38, for
the case of fixed end-points\footnote{In \cite{rs95} and in
  \cite{web02} there is a sign difference, due to the fact that in
  both papers the cotangent bundle is endowed with
  the symplectic form $-\omega$.}.

\begin{thm}
\label{compind}
Assume that $M$ is orientable, let $q\in \mathcal{P}_{\Lambda}(L)$,
and let $x\in \mathcal{P}_{\Lambda}(H)$, $(t,x(t)) =
\mathfrak{L}_L(t,q(t),\dot{q}(t))$, be the corresponding 1-periodic solution of
the Hamiltonian system on $T^* M$. Then
\[
m_{\Lambda}(q) = \mu_{\Lambda}(x).
\]
Let $q\in \mathcal{P}_{\Omega}(L,q_0,q_1)$ and let $x\in
\mathcal{P}_{\Omega}(H,q_0,q_1)$, $(t,x(t)) =
\mathfrak{L}_L(t,q(t),\dot{q}(t))$, be the corresponding solution of
the Hamiltonian system on $T^* M$. Then
\[
m_{\Omega}(q) = \mu_{\Omega}(x).
\]
\end{thm}

The following comparison between the Hamiltonian and the Lagrangian
action functionals follows immediately from the definition of the
Hamiltonian (\ref{leg}) and from (\ref{equ}):

\begin{lem}
\label{cmprsn}
If $x=(q,p):[0,1]\rightarrow T^* M$ is continuous, with $q$ of class
$W^{1,2}$, then
\[
\mathcal{A}(x) \leq \mathcal{E}(q),
\]
the equality holding if and only if $p(t)=dL(t,q(t),\dot{q}(t))|_{T_{
\dot{q}(t)}^v TM}$, that is if and only if $(t,q(t),p(t)) =
  \mathfrak{L}_L(t,q(t),\dot{q}(t))$ for every $t\in [0,1]$.
In particular, the Hamiltonian and the Lagrangian action functionals
coincide on the solutions of the two systems.
\end{lem}

Let $q$ be a solution in $\mathcal{P}_{\Lambda}(L)$ (resp.\
in $\mathcal{P}_{\Omega}(L)$), and
let $x=(q,p)$ be the corresponding
solution in $\mathcal{P}_{\Lambda}(H)$ (resp.\ in
$\mathcal{P}_{\Omega}(H)$). By the above lemma,
$\mathcal{A} \leq \mathcal{E} \circ \tau^*$ on
$\Lambda^1(T^* M)$ (resp.\ on $\Omega^1(T^* M)$), and
$\mathcal{A}(x)= \mathcal{E}(\tau^* \circ x)$.
So, taking also into account the fact that
$x$ is a critical point of $\mathcal{A}_{\Lambda}$ (resp.\
$\mathcal{A}_{\Omega}$), and $\tau^* \circ x =
q$ is a critical point of $\mathcal{E}_{\Lambda}$ (resp.\
$\mathcal{E}_{\Omega}$), we deduce the following:

\begin{lem}
\label{difdis}
Let $q$ be a solution in $\mathcal{P}_{\Lambda}(L)$ (resp.\
in $\mathcal{P}_{\Omega}(L)$), and
let $x=(q,p)$ be the corresponding
solution in $\mathcal{P}_{\Lambda}(H)$ (resp.\ in
$\mathcal{P}_{\Omega}(H)$). Then
\[
d^2 \mathcal{A} (x) [\zeta,\zeta] \leq d^2 \mathcal{E}
(q) [D \tau^*(x)[\zeta],D\tau^*(x)[\zeta]],
\]
for every $\zeta\in T_{x} \Lambda^1(T^* M)$ (resp.\ $\zeta\in
T_{x} \Omega^1(T^* M,q_0,q_1))$.
\end{lem}

Assumption (L1) implies that $L$ is bounded below, and so is
the action functional $\mathcal{E}$.

The metric $\langle \cdot, \cdot \rangle$ on $M$ induces a Riemannian
metric on the Hilbert manifolds $\Lambda^1(M)$ and
$\Omega^1(M,q_0,q_1)$, defined by
\[
\langle\langle \xi , \eta \rangle\rangle_1 = \int_0^1 \left( \langle
\nabla_t \xi,
  \nabla_t \eta \rangle + \langle \xi, \eta \rangle \right) \, dt,
\]
for $\xi,\eta$ elements of $T_q \Lambda^1(M)$ or of $T_q
\Omega^1(M,q_0,q_1)$. The corresponding distances on $\Lambda^1(M)$ and
on $\Omega^1(M,q_0,q_1)$ are compatible with the manifold topologies,
and they are complete. The following compactness result is proved in
\cite{ben86} (in the case of periodic orbits, the case of fixed
end-points is analogous):

\begin{prop}
\label{palsma}
The functional $\mathcal{E}$ satisfies the Palais-Smale condition on
the Riemannian manifold
$(\Lambda^1(M), \langle\langle \cdot , \cdot \rangle\rangle_1)$ and on
$(\Omega^1(M,q_0,q_1), \langle\langle \cdot , \cdot \rangle\rangle_1)$: every
sequence $(q_n)\subset \Lambda^1(M)$ (resp.\ $(q_n) \subset
\Omega^1(M,q_0,q_1)$) such that $\mathcal{E}(q_n)$ is bounded and
$\|\nabla \mathcal{E}_{\Lambda}(q_n)\|_1$ (resp.\ $\|\nabla
\mathcal{E}_{\Omega}(q_n)\|_1$) is infinitesimal, is compact.
\end{prop}

\subsection{The Morse complex}
\label{smc}

If a functional $f$ of class $C^h$, $h\in \{2,3,\dots,\infty\}$, on a
Hilbert Riemannian manifold $(\mathfrak{M}, \langle\langle \cdot, \cdot
\rangle\rangle)$ satisfies:
\begin{enumerate}
\item all the critical points $x$ of $f$ are non-degenerate and have
  finite Morse index $m(x)$;
\item $f$ is bounded below;
\item the Riemannian manifold
$(\mathfrak{M}, \langle\langle \cdot, \cdot \rangle\rangle)$ is complete;
\item $f$ satisfies the Palais-Smale condition on $(\mathfrak{M},
  \langle\langle \cdot, \cdot \rangle\rangle)$;
\end{enumerate}
we can associate to it a complex of Abelian groups, the {\em Morse
complex of $f$}, whose homology is isomorphic to the singular
homology of $\mathfrak{M}$. Its construction will be sketched in
this section. See \cite{ama04m} for full details.

Denote by $\crit(f)$ the set of critical points of $f$. Our assumptions
imply that the set
\begin{equation}
\label{fini}
\crit(f) \cap \{f\leq a\}
\end{equation}
is finite for every $a\in \R$. Denote by $\crit_k(f)$ the set of
critical points $x$ of $f$ of Morse index $m(x)=k$, and let $CM_k(f)$ be
the free Abelian group generated by the elements of
$\crit_k(f)$. Notice that $CM_k(f)$ may have infinite rank.

Denote by $\mathcal{G}$ the space of $C^{h-1}$ sections $G$ of the bundle
of endomorphisms of $T\mathfrak{M}$ such that
\[
G(p) \mbox{ is $\langle\langle \cdot , \cdot
\rangle\rangle$-symmetric for every $p\in \mathfrak{M}$, }
\|G\|_{C^{h-1}} < \infty \mbox{, and } \|G\|_{C^0} < 1,
\]
endowed with the (metrizable) topology of uniform convergence up to the
$(h-1)$-th derivative. If $G\in \mathcal{G}$,
\[
g(\xi,\eta) = g_G(\xi,\eta) := \langle\langle (I+G(p))
\xi, \eta \rangle\rangle ,
\quad \xi,\eta\in T_p \mathfrak{M},
\]
is a complete Riemannian metric on $\mathfrak{M}$, uniformly
equivalent to the original one, so $f$ satisfies the Palais-Smale
condition with respect to $g$. The gradient of $f$ with respect to
the metric $g$ will be denoted by $\nabla_g f$.

Let $\phi^t$ be the local integral flow of the vector field
$-\nabla_g f$. Its rest points are the critical points of $f$, and
$f$ is strictly decreasing on the non-constant orbits of $\phi^t$.
For $x\in \crit(f)$, let $T_x \mathfrak{M} = V^-(x) \oplus V^+(x)$
be the splitting of $T_x \mathfrak{M}$ corresponding to the
decomposition of the spectrum of the Hessian $\nabla^2_g f(x)$
into the negative and the positive part. By (i), $m(x) = \dim
V^-(x)$ is always finite. If $x\in \crit(f)$, the unstable and the
stable manifold of $x$,
\begin{eqnarray*}
W^u(x) = W^u(x;f,g) =
\set{p\in \mathfrak{M}}{\lim_{t\rightarrow -\infty} \phi^t(p)
  = x}, \\
W^s(x) = W^s(x;f,g) =
\set{p\in \mathfrak{M}}{\lim_{t\rightarrow
  +\infty} \phi^t(p)   = x},
\end{eqnarray*}
are images of $C^{h-1}$ embeddings $V^-(x)\hookrightarrow \mathfrak{M}$,
$V^+(x) \hookrightarrow \mathfrak{M}$, and $T_x W^u(x) = V^-(x)$, $T_x
W^s(x)= V^+(x)$. In particular,
\[
\dim W^u(x) = m(x), \quad \codim \, W^s(x) = m(x),
\]
$W^u(x)$ is orientable, and $W^s(x)$ is co-orientable, meaning that
its normal bundle is orientable. Assumptions (i-iv) on the functional
have the following consequences:

\begin{prop}
\label{mc}
Each unstable manifold $W^u(x)$ is pre-compact in $\mathfrak{M}$. The
closure of $W^u(x)$ is contained in
\[
W^u(x) \cup \bigcup_{\substack{y\in \crit(f)\\ f(y)<f(x)}} W^u(y).
\]
\end{prop}

\begin{prop}
\label{mtrans}
There is a residual set $\mathcal{G}_{\mathrm{reg}} \subset \mathcal{G}$ of
elements $G$ for which $W^u(x)$ and $W^s(y)$ intersect transversally,
whenever\footnote{The fact that we get transversal intersections only
  for index difference not exceeding $h-1$ is related to the fact that we
  are assuming $f$ to be only of class $C^h$. The possibility of
  keeping the regularity requirements low is important in nonlinear
  analysis, because
  functionals arising from smooth problems have often
  low regularity, and because in infinite dimensions $C^{h+1}$
  functionals are not dense in the space of $C^h$ ones, when $h\geq
  1$. Notice that
  $C^2$ regularity implies transversality
  up to index difference one, which is just what is needed for the
  construction of the Morse complex. In our case, the action
  functional $\mathcal{E}$ is smooth, so we get transversal
  intersections for arbitrary index difference, for a residual set of
  smooth metrics.}
 $m(x)-m(y)\leq h-1$. In particular, for $G\in
\mathcal{G}_{\mathrm{reg}}$:
\begin{enumerate}
\item if $m(x)\leq m(y)$ and $x\neq y$ then $W^u(x)\cap
  W^s(y)=\emptyset$;
\item if $m(x)-m(y)=1$ then $W^u(x)\cap W^s(y)$ is a one-dimensional
  manifold.
\end{enumerate}
\end{prop}

A metric $g=g_G$ coming from $G\in \mathcal{G}_{\mathrm{reg}}$ is
therefore a {\em Morse-Smale} metric for $f$. Let us fix such a metric $g$.
Consider the increasing sequence of open sets
\[
U_k = \set{\phi^t(p)}{t\geq 0, \; p\in \mathcal{U}(x), \; x\in
  \crit(f) \mbox{ with } m(x)\leq k}, \quad k\in \N,
\]
where $\mathcal{U}(x)$ is an open neighborhood of $x$. We also set
\[
U_{-1} = \emptyset, \quad \mbox{and} \quad U = \bigcup_{k\in \N} U_k.
\]
Using Proposition \ref{mtrans} (i) and the Palais-Smale condition,
it can be proved\footnote{Here we are assuming that for every
$k\in \N$ there are finitely many critical points of Morse index
$k$. In the case of infinitely many critical points with the same
index a stronger transversality assumption would be needed. In
that case, an easier way to construct the Morse complex is to deal
with each sublevel separately, and then take a direct limit. See
\cite{ama04m} for more details.} that if the neighborhoods
$\mathcal{U}(x)$ are suitably small, then the singular homology
groups of the pair $(U_k,U_{k-1})$ are
\begin{equation}
\label{cell}
H_j(U_k,U_{k-1}) =  \begin{cases} CM_k(f) & \mbox{if } j=k, \\ 0 &
  \mbox{if } j\neq k. \end{cases}
\end{equation}
Indeed, $H_k(U_k,U_{k-1})$ is the free Abelian group generated by the
relative homology classes of balls in $W^u(x)$, for $x\in \crit_k(f)$, chosen
to be so large that their boundary lies in $U_{k-1}$. Moreover, the
gradient flow of $f$ can be used to show that $U$ is a deformation
retract of $\mathfrak{M}$.
By (\ref{cell}), $\{U_k\}_{k\geq -1}$ is a cellular filtration of $U$
(see \cite{dol80}, section V.1),
and we define
\begin{eqnarray*}
\partial_k = \partial_k(f,g) : CM_k(f) \cong H_k(U_k,U_{k-1}) \rightarrow
H_{k-1}(U_{k-1}) \\ \rightarrow H_{k-1}(U_{k-1},U_{k-2}) \cong
CM_{k-1}(f)
\end{eqnarray*}
to be the associated cellular homomorphism.
By standard results about cellular filtrations, the above
homomorphisms are the data of a complex, the cellular complex
associated to the filtration $\{U_k\}_{k\geq -1}$, whose homology is isomorphic
to the singular homology of $U$:
\[
H_k (\{CM_*(f),\partial_*(f,g)\}) \cong H_k(U) \cong H_k(\mathfrak{M}).
\]
The complex $\{CM_*(f),\partial_*(f,g)\}$
is called the {\em Morse complex of $(f,g)$}.

Finally, let us describe what the boundary homomorphisms
$\partial_k(f,g)$ look like, in terms of the generators of
$CM_k(f)$. Let us fix an orientation of each unstable manifold
$W^u(x)$, in an arbitrary way. Consequently, we get a
co-orientation of each stable manifold $W^s(x)$. Since a
transversal intersection of an oriented submanifold and a
co-oriented submanifold has a canonical orientation, by
Proposition \ref{mtrans} (ii) we get an orientation of each
intersection $W^u(x)\cap W^s(y)$, in the case $m(x)-m(y)=1$. Let
$x\in \crit_k(f)$ and $y\in \crit_{k-1}(f)$. The compactness
expressed by Proposition \ref{mc} and the transversality expressed
by Proposition \ref{mtrans} imply that $W^u(x)\cap W^s(y)$
consists of finitely many flow lines. The flow line through $p$ -
denote it by $[p]$ - has the orientation defined above, and we
define $\epsilon([p])$ to be $+1$ if the tangent vector $-\nabla_g
f (p)$ is positively oriented, to be $-1$ in the opposite case.
Then the integer $n(x,y)$ is defined to be the sum of these
entries, $n(x,y) = \sum \epsilon([p])$, for $[p]$ varying among
all the flow lines connecting $x$ to $y$. It can be proved that
\[
\partial_k x = \sum_{y\in \crit_{k-1}(f)} n(x ,y)y.
\]
Notice that the above sum is finite because the set (\ref{fini}) is
finite. Furthermore, it can be proved that the isomorphism class of
the complex $\{CM_*(f),\partial_*(f,g_G)\}$ does not depend on $G\in
\mathcal{G}_{\mathrm{reg}}$:

\begin{prop}
\label{indep}
Assume that $g_0$ and $g_1$ are Morse-Smale metrics for $f$,
uniformly equivalent to $\langle\langle \cdot,\cdot
\rangle\rangle$. Then there exists a chain isomorphism
\[
\phi_{01} : \{CM_*(f),\partial_*(f,g_0)\}
\rightarrow \{CM_*(f),\partial_*(f,g_1)\}, \quad
x \mapsto \sum_{\substack{y\in \crit(f)\\ m(y)=m(x)}} n_{01} (x,y)y,
\]
such that $n_{01}(x,x)=1$ and $n_{01}(x,y)=0$ if $f(x)\leq f(y)$
and $x\neq y$.
\end{prop}
The proof is completely analogous to that of Theorem
\ref{chj}. We summarize the above discussion into the following:

\begin{thm}
Let $f$ be a $C^2$ functional on a Hilbert manifold $\mathfrak{M}$,
satisfying (i), (ii), (iii), (iv). Let $CM_k(f)$ be the free
Abelian group generated by its critical points of Morse index $k$.
Then the above construction produces a complex
\[
\partial_k : CM_k(f) \rightarrow CM_{k-1}(f), \quad \partial_k x =
\sum_{y\in \crit_{k-1}(f)} n(x ,y) y,
\]
uniquely determined up to isomorphism, whose homology groups are isomorphic to
the singular homology groups of $\mathfrak{M}$.
\end{thm}

Clearly, the Morse complex splits into subcomplexes, one for each
connected component of $\mathfrak{M}$, and the isomorphisms with
singular homology respects such a splitting. Moreover, the Morse
complex is filtered by the functional level: if $a\in \R$, the
boundary homomorphism maps $CM_k^a(f) := \Span (\crit_k (f) \cap
\{f<a\})$ into $CM_{k-1}^a(f)$. So $\{CM^a_*(f), \partial_*(f,g)\}$
is a subcomplex of $\{CM_*(f),\partial_*(f,g)\}$, and its
homology is seen to be isomorphic to the singular
homology of the sublevel $\{f<a\}$. Both the splitting into
subcomplexes and the $\R$ filtering are compatible with the
isomorphisms of Proposition \ref{indep}.

\subsection{The Morse complex of $\mathbb{\mathcal{E}}$}

By what we have seen, assumptions (L0), (L1),
and (L2) imply that both $\mathcal{E}_{\Lambda}$ and
$\mathcal{E}_{\Omega}$ satisfy the conditions (i), (ii), (iii), and
(iv) of section \ref{smc}. Therefore, if
$CM_k(\mathcal{E}_{\Lambda})$ (resp.\
$CM_k (\mathcal{E}_{\Omega})$) denotes
the free Abelian group generated by
the solutions in $\mathcal{P}_{\Lambda}(L)$ (resp.\
$\mathcal{P}_{\Omega}(L,q_0,q_1)$) of Morse index $i_{\Lambda}(q)=k$
(resp.\ $i_{\Omega}(q)=k$)
we get an isomorphism class of complexes
\[
\partial_k : CM_k(\mathcal{E}_{\Lambda}) \rightarrow
CM_{k-1}(\mathcal{E}_{\Lambda}), \quad \mbox{(resp.\ } \;  \partial_k :
CM_k(\mathcal{E}_{\Omega}) \rightarrow CM_{k-1}(\mathcal{E}_{\Omega})\mbox{),}
\]
whose homology is isomorphic to the singular homology of
$\Lambda^1(M)$ (respectively of $\Omega^1(M,q_0,q_1)$). Since the inclusions
\begin{eqnarray*}
\Lambda^1(M) \hookrightarrow \Lambda(M), \\
\Omega^1(M,q_0,q_1) \hookrightarrow \set{q\in
  C^0([0,1],M)}{q(0)=q_0, \; q(1)=q_1}
\end{eqnarray*}
are homotopy equivalence
and since the latter space is homotopically equivalent to
the based loop space $\Omega(M)$, we deduce
that the homology of the above complexes are isomorphic to the
singular homology of the free loop space of $M$, and of the based loop
space of $M$:
\[
H_k(\{CM_*(\mathcal{E}_{\Lambda}),\partial_*\}) \cong H_k(\Lambda(M)),
\quad H_k(\{CM_*(\mathcal{E}_{\Omega}),\partial_*\}) \cong
H_k(\Omega(M)).
\]
The complex $\{CM_*(\mathcal{E}_{\Lambda}),\partial_*\}$ (resp.\
$\{CM_*(\mathcal{E}_{\Omega}),\partial_*\}$) splits into subcomplexes, one for
each conjugacy class
of $\pi_1(M)$ (resp.\ one for each element of $\pi_1(M)$). Finally, if
$CM_k^a(\mathcal{E}_{\Lambda})$ (resp.\ $CM_k^a(\mathcal{E}_{\Omega})$
denotes the
subgroup of $CM_k(\mathcal{E}_{\Lambda})$ (resp.\
$CM_k(\mathcal{E}_{\Omega})$) generated by solutions of action less
than $a$, we obtain a filtering of the above complex by subcomplexes
such that
\begin{eqnarray*}
H_k(\{CM_*^a(\mathcal{E}_{\Lambda}),\partial_*\}) \cong
H_k(\{\mathcal{E}_{\Lambda} < a\}),
\quad \mbox{(resp.\ } H_k(\{CM_*^a(\mathcal{E}_{\Omega}),\partial_*\}) \cong
H_k(\{\mathcal{E}_{\Omega}<a\}) \mbox{ ).}
\end{eqnarray*}

\section{The isomorphism between the Morse and the Floer complex}
\label{tre}

We are now ready to state and prove the main result of this paper:

\begin{thm}
\label{prec} Assume that the Hamiltonian $H\in C^{\infty}(\T\times
T^*M)$ satisfies (H0), (H1), and (H2).
Assume also that $H$ is the Legendre transform of the
Lagrangian $L\in C^{\infty}(\T \times TM)$ satisfying (L0), (L1),
and (L2). Let $J$ be a $t$-dependent, $t\in \T$,
$\omega$-compatible almost complex structure on $T^* M$, belonging
to $\mathcal{J}_{j_0,\mathrm{reg}}(H)$. Let $g$ be a Riemannian structure on
$\Lambda^1(M)$ (resp.\ on $\Omega^1(M,q_0,q_1)$) uniformly
equivalent to $\langle\langle \cdot,\cdot \rangle\rangle_1$, and
having the Morse-Smale property with respect to $\mathcal{E}$.
Then there is a chain complex isomorphism
\[
\Theta : \{CM_*(\mathcal{E}),\partial_*(\mathcal{E},g)\}
\longrightarrow \{CF_*(H),\partial_*(H,J) \}
\]
of the form
\[
\Theta q =\sum_{\substack{x\in \mathcal{P}(H)\\ \mu(x) = m(q)}} n^+(q,x) \,
x, \quad \forall q\in \mathcal{P}(L),
\]
such that $n^+(q,x)=0$ if $\mathcal{E}(q)\leq \mathcal{A}(x)$,
unless $q$ and $x$ correspond to the same solution - that is
$\mathfrak{L}(t,q(t),\dot{q}(t))=(t,x(t))$ for every $t\in [0,1]$
- in which case $n^+(q,x)=\pm 1$.

In particular, $\Theta$ induces an isomorphism between the
subcomplexes $\{CM_*^a(\mathcal{E})\}$ and $\{CF_*^a(H)\}$, for every $a\in
\R$. Finally, $\Theta$ is compatible with the splitting of the
Floer and the Morse complex into the subcomplexes corresponding to
different conjugacy classes of $\pi_1(M)$ (resp.\ of different
elements of $\pi_1(M)$).
\end{thm}

Fix some $r\in ]2,4]$.
Given $q\in \mathcal{P}_{\Lambda}(L)$ and $x\in
\mathcal{P}_{\Lambda}(H)$, let $\mathcal{M}_{\Lambda}^+(q,x)$
be the set of all maps
\[
u\in C^{\infty}(]0,+\infty[\times\T,T^* M)\cap W^{1,r}(]0,1[\times
    \T,T^*M)
\]
such that
\begin{eqnarray*}
\partial_s u - J(t,u) (\partial_t u - X_H(t,u))=0 \quad \mbox{on }
  ]0,+\infty[ \times \T, \\
\tau^* u(0,\cdot) \in W^u(q) = W^u(q;\mathcal{E}_{\Lambda},g_{\Lambda}), \\
\lim_{s\rightarrow +\infty} u(s,t) = x(t) \quad \mbox{uniformly in }
  t\in \T.
\end{eqnarray*}
Similarly, if $q\in \mathcal{P}_{\Omega}(L)$ and $x\in
\mathcal{P}_{\Omega}(H)$, let $\mathcal{M}_{\Omega}^+(q,x)$
be the set of all maps
\[
u\in C^{\infty}(]0,+\infty[\times[0,1],T^* M)\cap W^{1,r}(]0,1[\times
    ]0,1[,T^*M)
\]
such that
\begin{eqnarray*}
\partial_s u - J(t,u) (\partial_t u - X_H(t,u))=0 \quad \mbox{on }
  ]0,+\infty[ \times [0,1], \\
\tau^* u(s,0) =q_0, \; \tau^*u(s,1) =q_1 \; \forall s\geq 0, \quad
  \mbox{and} \quad
\tau^* u(0,\cdot) \in W^u(q) = W^u(q;\mathcal{E}_{\Omega},g_{\Omega}), \\
\lim_{s\rightarrow +\infty} u(s,t) = x(t) \quad \mbox{uniformly in }
  t\in [0,1].
\end{eqnarray*}
It makes sense to look for solutions $u$ which are of Sobolev class
$W^{1,r}$ near $\{0\}\times [0,1]$, with $r\in ]2,4]$, because in this case
$u(0,\cdot)$ belongs to the Sobolev space $W^{1-1/r,r}(]0,1[,T^*M)$,
and the boundary condition $\tau^* u(0,\cdot)\in W^u(q)$ is
well-posed because $W^u(q)$ consists of curves in $W^{1,2}(]0,1[)$, which
continuously embeds into $W^{1-1/r,r}(]0,1[)$ for $r\leq 4$ (see
\cite{ada75}, section 7.5.8).

The action functional satisfies the estimates
\[
\mathcal{E}(\phi_{-\nabla_g \mathcal{E}}^{t_0}(\bar{q})) =
\mathcal{E}(q) - \int_{-\infty}^{t_0} g(\nabla_g
\mathcal{E}(\phi_{-\nabla_g \mathcal{E}}^t(\bar{q})),\nabla_g
\mathcal{E}(\phi_{-\nabla_g \mathcal{E}}^t(\bar{q})))  \, dt <
\mathcal{E}(q)
\]
for every $\bar{q}\in W^u(q) \setminus \{q\}$, and
\[
\mathcal{A}(u(s_0,\cdot)) = \mathcal{A}(x) + \int_{s_0}^{+\infty}
\int_0^1 |\nabla_J
\mathcal{A}(u(s,\cdot))(t)|^2_{J_t} \, dt \, ds >
\mathcal{A}(x)
\]
for every solution $u$ of $\partial_s u -
J(t,u)(\partial_t u - X_H(t,u))=0$ converging to $x$ for
$s\rightarrow +\infty$, different from the stationary solution
$u(s,t)=x(t)$. Then Lemma \ref{cmprsn} implies that if $u\in
\mathcal{M}^+ (q,x)$, $\bar{q}(t) := \tau^* u(0,t)$, and $t \leq 0
\leq s$, then
\begin{equation}
\label{aces} \mathcal{E}(q) \geq \mathcal{E}(\phi_{-\nabla
\mathcal{E}}^t (\bar{q})) \geq \mathcal{E}(\bar{q}) \geq
\mathcal{A} (u(0,\cdot)) \geq \mathcal{A} (u(s,\cdot)) \geq
\mathcal{A}(x) ,
\end{equation}
and $\mathcal{E}(q)=\mathcal{A}(x)$ if and only if $q$ and $x$
correspond to the same solution by the Legendre transform - that
is $\mathfrak{L}_L(t,q(t),\dot{q}(t)) = (t,x(t))$ for every $t$ -
in which case $\mathcal{M}^+ (q,x)$ consists of a single element,
the stationary solution $x$.

\subsection{The Fredholm theory}
\label{fts}

Let us describe the functional setting which allows to see
$\mathcal{M}^+(q,x)$ as the set of zeros of a smooth section of a
Banach bundle. In the $\Lambda$-case define
$\mathcal{B}_{\Lambda}^+ = \mathcal{B}_{\Lambda}^+ (q,x)$ to be
the set of maps $u: [0,+\infty[ \times \T \rightarrow T^* M$ which
are of Sobolev class $W^{1,r}$ on every compact subset of
$[0,+\infty[ \times \T$ and such that:
\begin{enumerate}
\item $\tau^* u(0,\cdot) \in W^u(q)$; \item there is $s_0\geq 0$
for which
\[
u(s,t) = \exp_{x(t)} (\zeta(s,t)) \quad \forall (s,t)\in
]s_0,+\infty[ \times \T,
\]
where $\zeta$ is a $W^{1,r}$ section of $x^*(TT^*M) \rightarrow
]s_0, +\infty[ \times \T$.
\end{enumerate}
The set $\mathcal{B}_{\Lambda}^+$ has a natural structure of
smooth Banach manifold, and its tangent space at $u\in
\mathcal{B}_{\Lambda}^+$ is identified with the space of $W^{1,r}$
sections $w$ of $u^*(TT^*M)$ such that
\begin{equation}
\label{tang} D \tau^* (u(0,\cdot)) [w(0,\cdot)] \in T_{\tau^* u(0,\cdot)}
W^u(q).
\end{equation}
Similarly, $\mathcal{B}_{\Omega}^+ = \mathcal{B}_{\Omega}^+ (q,x)$
will be the Banach manifold of all maps $u : [0,+\infty[ \times
[0,1]$ $\rightarrow T^* M$ which are of Sobolev class $W^{1,r}$ on
every compact subset of $[0,+\infty[ \times [0,1]$ such that
$u(s,0)\in T_{q_0}^* M$, $u(s,1)\in T_{q_1}^* M$ for every $s\geq
0$, and (i) and the analogous of (ii) hold.

We denote by $\mathcal{W}_{\Lambda}^+ = \mathcal{W}_{\Lambda}^+
(q,x)$ (resp.\ $\mathcal{W}_{\Omega}^+ = \mathcal{W}_{\Omega}^+
(q,x)$) the Banach bundle over $\mathcal{B}_{\Lambda}^+$ (resp.\
$\mathcal{B}_{\Omega}^+$) whose fiber $\mathcal{W}_u^+$ at $u$ is
the space of $L^r$ sections of $u^*(TT^*M)$. Standard elliptic
regularity results and the exponential convergence guaranteed by
(H0) imply that $\mathcal{M}^+(q,x)$ is the set of zeros of the
smooth section
\[
\partial_{J,H}^+ : \mathcal{B}^+ \rightarrow \mathcal{W}^+,
\quad u \mapsto \partial_s u - J(t,u)(\partial_t u -
X_H(t,u)).
\]
The aim of this section is to prove the following:

\begin{thm}
\label{fredth}
If $q\in \mathcal{P}(L)$, $x\in \mathcal{P}(H)$,
and $u\in \mathcal{M}^+(q,x)$, then the fiberwise derivative of
$\partial_{J,H}^+ : \mathcal{B}^+ (q,x) \rightarrow
\mathcal{W}^+ (q,x)$ is a Fredholm operator of index
\[
\ind D_f \partial_{J,H}^+ (u) = m(q) - \mu(x).
\]
\end{thm}

By (\ref{tang}), the space $W_u$ of $W^{1,r}$ sections $w$ of
$u^*(TT^*M)$ such that $w(0,\cdot)\in T_{u(0,\cdot)}^v T^*M$ - or
equivalently $D \tau^* (u(0,\cdot)) [w(0,\cdot)] = 0$ - is a closed linear
subspace of $T_u \mathcal{B}^+ (q,x)$ of codimension $m(q)$.
Therefore it is enough to prove that the restriction of
$D_f \partial_{J,H}^+ (u)$ to $W_u$ is a Fredholm operator of
index $-\mu(x)$.

We recall that $\lambda_0$ denotes the vertical Lagrangian
subspace of $(\R^{2n},\omega_0)$, $\lambda_0 = (0)\times \R^n$.
The proof of the following lemma is analogous to that of Lemma
\ref{ltriv} (it is actually simpler):

\begin{lem}
\label{trivbis} Let $u\in \mathcal{M}_{\Lambda}^+ (q,x)$ (resp.\
$\mathcal{M}_{\Omega}^+ (q,x)$), and let $\Phi^+: \T \times
\R^{2n} \rightarrow x^*(TT^*M)$ (resp.\ $\Phi^+ : [0,1] \times
\R^{2n} \rightarrow x^*(TT^*M)$) be a unitary trivialization such
that $\Phi^+(\cdot)\lambda_0 = T_{x(\cdot)}^v T^*M$. Then there
exists a unitary trivialization
\[
\Phi : [0,+\infty] \times \T \times \R^{2n} \rightarrow
u^*(TT^*M) \quad \mbox{(resp. } \Phi : [0,+\infty] \times [0,1]
\times \R^{2n} \rightarrow u^*(TT^*M) \mbox{)}
\]
which is smooth on $]0,+\infty] \times \T$ (resp.\ on $]0,+\infty]
    \times [0,1]$) and of class $W^{1,r}$ on $]0,1[\times \T$ (resp.\
    on $]0,1[ \times ]0,1[$), such that $\Phi(s,t) \lambda_0 =
    T_{u(s,t)}^v T^*M$ for every $(s,t)\in [0,+\infty]\times \T$
    (resp.\ $[0,+\infty] \times [0,1]$), and $\Phi(+\infty,\cdot)
    = \Phi^+(\cdot)$.
\end{lem}

The trivialization $\Phi$ given by Lemma \ref{trivbis} defines a
conjugacy between the restriction to $W_u$ of
$D_f\partial_{J,H}^+ (u): T_u \mathcal{B}^+ \rightarrow
\mathcal{W}^+_u$ and a bounded operator
\begin{eqnarray*}
D^+_{S,\Lambda} : W^{1,r}_{\lambda_0} (]0,+\infty[ \times
\T,\R^{2n}) \rightarrow L^r (]0,+\infty[ \times \T,\R^{2n}), \\
\mbox{(resp.\ } D^+_{S,\Omega} : W^{1,r}_{\lambda_0} (]0,+\infty[
\times ]0,1[ , \R^{2n}) \rightarrow L^r ( ]0,+\infty[ \times
]0,1[, \R^{2n}) \mbox{ )}
\end{eqnarray*}
of the form
\[
D_S^+ v = \partial_s v - J_0 \partial_t v - S(s,t) v.
\]
Here $S: [0,+\infty[ \times \T \rightarrow \mathfrak{gl} (2n,\R)$
(resp.\ $S : [0,+\infty[ \times [0,1] \rightarrow \mathfrak{gl}
(2n,\R)$) has the form
\[
S = \Phi^{-1} (\nabla_s \Phi - J(t,u)\nabla_t \Phi - \nabla_{\Phi}
J(t,u)\partial_t u + \nabla_{\Phi} \nabla H(t,u)).
\]
Therefore, $S$ is smooth on $]0,+\infty[ \times \T$ (resp.\ on
$]0,+\infty[ \times [0,1]$), and it is of class $L^r$ on
$]0,1[\times \T$ (resp.\ on $]0,1[\times ]0,1[$). Moreover,
\[
\lim_{s\rightarrow +\infty} S(s,t) = S(+\infty,t) \quad
\mbox{uniformly in } t,
\]
where $S(+\infty,t)$ is symmetric for every $t$, and the solution of
\begin{equation}
\label{stra} \gamma^{\prime}(t) = J_0 S(+\infty,t) \gamma(t) ,
\quad \gamma(0) = I,
\end{equation}
is conjugated to the differential of the Hamiltonian flow along $x$:
\[
\gamma(t) = \Phi(+\infty,t)^{-1} D \phi_H^t(x(0))
\Phi(+\infty,0).
\]
Hence (H0)$_{\Lambda}$ is translated into the condition
\begin{equation}
\label{unno} \mbox{1 is not an eigenvalue of } \gamma(1),
\end{equation}
while (H0)$_{\Omega}$ is translated into the condition
\begin{equation}
\label{ddue} \gamma(1) \lambda_0 \cap \lambda_0 = (0).
\end{equation}
Therefore Theorem \ref{fredth} is a consequence of the following:

\begin{thm}
\label{fftt}
If $r\in ]2,+\infty[$ the following facts hold:

\noindent (a) Assume that the matrix valued map
$S:]0,+\infty[\times \T \rightarrow
\mathfrak{gl} (2n,\R)$ is of the form $S=S_1+S_2$, where $S_1\in
L^{\infty} (]0,+\infty[ \times \T,\mathfrak{gl}(2n,\R))$, and $S_2\in
L^r (]0,+\infty[ \times \T,\mathfrak{gl}(2n,\R))$ has bounded
support. Moreover, assume that the limit
\[
S(+\infty,t) := \lim_{s\rightarrow +\infty} S(s,t) =
\lim_{s\rightarrow +\infty} S_1(s,t)
\]
is uniform in $t\in\T$, and that the solution $\gamma$ of (\ref{stra})
satisfies (\ref{unno}). Then the bounded linear operator
\[
D_{S,\Lambda}^+ : W^{1,r}_{\lambda_0} (]0,+\infty[ \times
\T,\R^{2n}) \rightarrow L^r ( ]0, +\infty[ \times \T, \R^{2n}) ,
\quad u \mapsto \partial_s u - J_0 \partial_t u - S(s,t)u,
\]
is Fredholm of index
\begin{equation}
\label{AA}
\ind D_{S,\Lambda}^+ = - \mu_{CZ} (\gamma).
\end{equation}
(b) Assume that the matrix-valued map $S:]0,+\infty[\times ]0,1[ \rightarrow
\mathfrak{gl} (2n,\R)$ is of the form $S=S_1+S_2$, where 
\[
S_1\in L^{\infty} (]0,+\infty[ \times ]0,1[,\mathfrak{gl}(2n,\R)),
\quad S_2\in L^r (]0,+\infty[ \times ]0,1[,\mathfrak{gl}(2n,\R)),
\]
and $S_2$ has bounded support. Moreover, assume that the limit
\[
S(+\infty,t) := \lim_{s\rightarrow +\infty} S(s,t) =
\lim_{s\rightarrow +\infty} S_1(s,t)
\]
is uniform in $t\in [0,1]$, and that the solution $\gamma$ of (\ref{stra})
satisfies (\ref{ddue}). Then the bounded linear operator
\begin{eqnarray*}
D_{S,\Omega}^+ : W^{1,r}_{\lambda_0} (]0,+\infty[ \times
]0,1[,\R^{2n}) \rightarrow L^r ( ]0, +\infty[ \times ]0,1[,
\R^{2n}) , \\ u \mapsto \partial_s u - J_0 \partial_t u -
S(s,t)u,
\end{eqnarray*}
is Fredholm of index
\begin{equation}
\label{BB}
\ind D_{S,\Omega}^+ = \frac{n}{2} -
\mu(\gamma(\cdot)\lambda_0,\lambda_0).
\end{equation}
\end{thm}

\begin{proof} The multiplication operator
\begin{eqnarray*}
W^{1,r}_{\lambda_0} (]0,+\infty[ \times
\T,\R^{2n}) \rightarrow L^r ( ]0, +\infty[ \times \T, \R^{2n}) ,
\quad u \mapsto S_2 u, \\
\mbox{(resp.\ } W^{1,r}_{\lambda_0} (]0,+\infty[ \times
]0,1[,\R^{2n}) \rightarrow L^r ( ]0, +\infty[ \times ]0,1[,
\R^{2n}) , \quad u \mapsto S_2 u \mbox{ )}
\end{eqnarray*}
is compact. Indeed, if $S_2$ has support in $[0,s_0] \times [0,1]$,
the above operator factorizes through
\[
W^{1,r} (]0,s_0[ \times ]0,1[) \hookrightarrow L^{\infty}(]0,s_0[ \times
  ]0,1[) \stackrel{S_2 \cdot}{\longrightarrow} L^r (]0,s_0[ \times ]0,1[),
\]
where the inclusion is compact because $r>2$, and the second map
is continuous because $S_2\in L^r$. Since a compact perturbation
of a Fredholm operator is still Fredholm with the same index, the
presence of the term $S_2$ does not change the Fredholm property
and the index. Therefore, in the remaining part of the proof we
may assume that $S_2=0$, hence that $S=S_1$ is in $L^{\infty}$. In
this case statements (a) and (b) will actually hold for every
$r\in ]1,+\infty[$.

\medskip

The proof that $D_{S,\Lambda}^+$ and
$D_{S,\Omega}^+$ are Fredholm is now standard; let us sketch the argument,
referring to \cite{sal99} for more details. The main ingredients are
the following facts:
\begin{enumerate}
\item By condition (\ref{unno}) (resp.\ (\ref{ddue})), the
translation invariant operator
\begin{eqnarray*}
D_{S(+\infty,\cdot), \Lambda} : W^{1,r} (\R \times \T,\R^{2n})
\rightarrow L^r (\R \times \T,\R^{2n}) \\
\mbox{(resp.\ } D_{S(+\infty,\cdot),\Omega} : W^{1,r}_{\lambda_0}
(\R \times ]0,1[, \R^{2n}) \rightarrow L^r(\R \times ]0,1[,
\R^{2n}) \mbox{ )}
\end{eqnarray*}
mapping $u$ into $\partial_s u - J_0 \partial_t u - S(+\infty,t)
u$, is invertible.

\item The Calderon-Zygmund inequality stated in Proposition
\ref{caldzyg} implies that there exists $c_0$ such that
\[
\|u\|_{W^{1,r}(]0,+\infty[ \times ]0,1[)} \leq c_0 (
\|D_S^+\|_{L^r(]0,+\infty[ \times ]0,1[)} + (1+\|S\|_{L^{\infty}})
\|u\|_{L^r(]0,+\infty[ \times ]0,1[)} ),
\]
for every $u\in W^{1,r}_{\lambda_0} (]0,+\infty[ \times \T,\R^{2n})$ (resp.\
$u\in W^{1,r}_{\lambda_0} (]0,+\infty[ \times ]0,1[,\R^{2n})$).

\item Let us identify the cokernel of $D_S^+$ with the annihilator
of the range of $D_S^+$: if $1/r + 1/r^{\prime} =1$ we have
\begin{eqnarray*}
\coker D_{S,\Lambda}^+ \cong (\ran D^+_{S,\Lambda})^{\circ} =
\Bigl\{v\in L^{r^{\prime}} (]0,+\infty[ \times \T,\R^{2n}) \, \Big| \\
\int_{]0,+\infty[\times \T} v \cdot D_{S,\Lambda}^+ u \, ds\, dt =
0 \;\;\; \forall u\in W^{1,r}_{\lambda_0} (]0,+\infty[\times
\T,\R^{2n})\Bigr\}, \\
\Bigl(\mbox{resp.\ } \coker D_{S,\Omega}^+ \cong (\ran
D^+_{S,\Omega})^{\circ} = \Bigl\{v\in L^{r^{\prime}} (]0,+\infty[
\times ]0,1[,\R^{2n})\, \Big| \\ \int_{]0,+\infty[\times ]0,1[} v
\cdot D_{S,\Lambda}^+ u \, ds\, dt = 0 \;\;\; \forall u\in
W^{1,r}_{\lambda_0} (]0,+\infty[\times ]0,1[,\R^{2n})\Bigr\}
\Bigr).
\end{eqnarray*}
Then the regularity theory for the weak solutions of the
Cauchy-Riemann operator implies the following facts:
\begin{enumerate}
\item the cokernel of $D_{S,\Lambda}^+$ consists of the maps
$v\in W^{1,r^{\prime}}(]0,+\infty[\times \T, \R^{2n})$ solving
\begin{equation}
\label{adj}
\partial_s v + J_0 \partial_t v + S(s,t)^T v = 0,
\end{equation}
and such that $v(0,t) \in \lambda_0^{\perp} = \R^n \times (0)$ for
every $t\in \T$;

\item the cokernel of $D_{S,\Omega}^+$ consists of the maps
$v\in W^{1,r^{\prime}}(]0,+\infty[\times ]0,1[,\R^{2n})$ solving
(\ref{adj}), and such that
\[
v(0,t) \in \lambda_0^{\perp} = \R^n \times (0) \;\; \forall t\in
[0,1], \quad v(s,0),\; v(s,1)\in \lambda_0 = (0) \times \R^n \;\;
\forall s\geq 0.
\]
\end{enumerate}

\end{enumerate}

By (i) and (ii) there exist constants $s_0$ and $c_1$ such that
\[
\|u\|_{W^{1,r}(]0,+\infty[\times ]0,1[)} \leq c_1 \left( \|D_S^+
u\|_{L^r(]0,+\infty[\times ]0,1[)} + \|u\|_{L^r(]0,s_0[ \times
]0,1[)} \right).
\]
Then the compactness of the embedding
\[
W^{1,r}(]0,s_0[\times ]0,1[) \hookrightarrow L^r(]0,s_0[\times
]0,1[)
\]
implies that $D_S^+$ has finite dimensional kernel and closed
range. Finally, (iii) implies that the cokernel of $D_S^+$
can be identified with the kernel of an operator of the same kind
(with different boundary conditions), which is finite dimensional
by the previous argument.

\medskip

There remains to compute the Fredholm index of $D_{S,\Lambda}^+$
and $D_{S,\Omega}^+$. We will make this computation when $n=1$ and $S$ is
a suitable constant matrix, and then we will use a homotopy argument to
pass to a general $S$. Notice that if $S$ is constant, the elements of
the kernel and cokernel of $D_S^+$ are smooth up to the boundary, and
the asymptotic conditions (\ref{unno}) and (\ref{ddue}) imply that
they decay exponentially fast for $s\rightarrow +\infty$.

\medskip

{\sc Claim 1.} If $n=1$ and $S(s,t) = Q_{\alpha} = \left(
\begin{array}{cc} 0 & \alpha \\ \alpha & 0 \end{array} \right)$,
with $\alpha\in \R \setminus \{0\}$, then
\[
\ind D_{S,\Lambda}^+ = 0 = -\mu_{CZ}(\gamma).
\]

\medskip

Notice that in this case $\gamma(t)= \left(
\begin{array}{cc} e^{\alpha t} & 0 \\ 0 & e^{-\alpha t} \end{array}
\right)$, so (\ref{unno}) is satisfied, and by (\ref{czi})
$\mu_{CZ}(\gamma)=0$.

Let $u\in \ker D_{Q_{\alpha},\Lambda}^+$. Then $u$ is a smooth
solution in $W^{1,r}(]0,+\infty[ \times \T,\R^2)$ of
\[
\begin{cases} \partial_s u - J_0 \partial_t u -
Q_{\alpha} u = 0 & \mbox{on } [0,+\infty[ \times \T, \\
u_1(0,t) = 0 & \forall t\in \T. \end{cases}
\]
Then the function
\[
w(s,t) = \begin{cases} u(s,t) & \mbox{if } s\geq 0, \\
(-u_1(-s,t),u_2(-s,t)) & \mbox{if } s<0, \end{cases}
\]
belongs to $W^{1,r}(\R \times \T,\R^2)$, and solves the same
problem. Hence, $w$ belongs to the kernel of the translation invariant
operator
$D_{Q_{\alpha},\Lambda}$, which is invertible by (i), so $w=0$.
By (iii-a), a similar argument shows that also the cokernel of
$D_{Q_{\alpha},\Lambda}^+$ is $(0)$. Therefore
$D_{Q_{\alpha},\Lambda}^+$ is invertible, and in particular its
index is 0.

\medskip

{\sc Claim 2.} If $n=1$ and $S(s,t)=\theta I$, with $\theta\in \R
\setminus 2\pi\Z$, then
\begin{equation}
\label{11}
\ind D_{\theta I,\Lambda}^+ = -2 \left\lfloor \frac{\theta}{2\pi}
\right\rfloor - 1 = - \mu_{CZ} (\gamma).
\end{equation}

\medskip

Notice that in this case $\gamma(t)=e^{t\theta J_0}$, so condition
(\ref{unno}) is equivalent to $\theta\notin 2\pi \Z$, and by (\ref{czi}),
$\mu_{CZ} (\gamma) = 2 \lfloor \theta/(2\pi) \rfloor + 1$.

By separating the variables, it is easily seen that the solutions
of the equation
\begin{equation}
\label{laeq}
\partial_s u - J_0 \partial_t u - \theta u = 0 \quad \mbox{on }
[0,+\infty[ \times \T
\end{equation}
have the form
\[
u(s,t) = \sum_{h\in \Z} e^{(\theta-2\pi h)s} e^{2\pi htJ_0}
\zeta_h, \quad \zeta_h = (\xi_h,\eta_h)\in \R^2.
\]
In order for such a function to decay for $s\rightarrow +\infty$,
it is necessary that $\zeta_h=0$ whenever $\theta-2\pi h>0$, so in
the sum above $h$ ranges from $\lceil \theta/(2\pi) \rceil$ to
$+\infty$. In particular, the first component of $u(0,t)$ is
\begin{equation}
\label{fcu} u_1(0,t) = \sum_{h= \lceil \theta/(2\pi)
\rceil}^{+\infty} ( \xi_h \cos 2\pi h t + \eta_h \sin 2\pi h t ).
\end{equation}
Recalling that $\{ 1, \sin 2\pi t, \cos 2\pi t, \sin 4 \pi t, \cos
4 \pi t, \dots \}$ is a complete orthogonal family in $L^2(\T)$,
we find that:
\begin{itemize}
\item if $\theta>0$, $u_1(0,\cdot)$ vanishes identically on $\T$
if and only if $\xi_h = \eta_h =0$ for every $h$; \item if
$\theta<0$, (\ref{fcu}) can be rewritten as
\begin{eqnarray*}
u_1(0,t) = \xi_0 + \sum_{h=1}^{-\lceil \theta/(2\pi)\rceil }
((\xi_h+ \xi_{-h})\cos 2\pi h t + (\eta_h - \eta_{-h} ) \sin 2\pi
h t ) \\ + \sum_{h=-\lceil \theta/(2\pi) \rceil}^{+\infty} ( \xi_h
\cos 2\pi h t + \eta_h \sin 2\pi h t ),
\end{eqnarray*}
so $u_1(0,\cdot)$ vanishes identically on $\T$ if and only if
$\xi_h = \eta_h = 0$ for $h\geq - \lceil \theta/(2\pi) \rceil+1$,
$\xi_0=0$, $\xi_{-h}=-\xi_h$ and $\eta_{-h} = \eta_h$ for $h\in \{
1,\dots, - \lceil \theta/ (2\pi) \rceil \}$.
\end{itemize}

We conclude that the kernel of $D_{\theta I,\Lambda}^+$ is $(0)$
when $\theta>0$, and it consists of the functions
\[
e^{\theta s} \binom{0}{\eta_0} + \sum_{h=1}^{-\lceil \theta/(2\pi)
\rceil} \left( e^{(\theta-2\pi h)s} e^{2\pi htJ_0}
\binom{\xi_h}{\eta_h} + e^{(\theta +2\pi h)s} e^{-2\pi ht J_0}
\binom{-\xi_h}{\eta_h} \right),
\]
when $\theta<0$. Therefore
\[
 \dim \ker D^+_{\theta I,\Lambda} = \begin{cases} 0 & \mbox{if }
 \theta>0, \\ 1 - 2 \left\lceil \frac{\theta}{2\pi} \right\rceil &
 \mbox{if } \theta<0. \end{cases}
\]
By (iii-a), the annihilator of the range of $D_{\theta
 I,\Lambda}^+$ consists of the smooth $\R^{2n}$-valued maps on
 $[0,+\infty[ \times \T$ which solve
 \[
 \begin{cases}\partial_s v + J_0 \partial_t v + \theta v = 0 &
 \mbox{on } ]0,+\infty[ \times \T, \\ v_1(0,t) = 0 & \forall t\in
 \T, \\ v(s,\cdot)\rightarrow 0 & \mbox{for } s\rightarrow +\infty.
 \end{cases}
 \]
 Since $v$ solves the above system if and only if $w(s,t) = J_0
 v(s,-t)$ solves
\[
 \begin{cases}\partial_s w - J_0 \partial_t w + \theta w = 0 &
 \mbox{on } ]0,+\infty[ \times \T, \\ w_2(0,t) = 0 & \forall t\in
 \T, \\ w(s,\cdot) \rightarrow 0 & \mbox{for } s\rightarrow +\infty,
 \end{cases}
 \]
we find
\[
\dim \coker D_{\theta I, \Lambda}^+ = \dim \ker D_{-\theta
I,\Lambda}^+ = \begin{cases} 1 - 2 \left\lceil -
\frac{\theta}{2\pi} \right\rceil & \mbox{if } \theta>0, \\ 0 &
\mbox{if } \theta<0, \end{cases}
\]
and the index formula (\ref{11}) follows.

\medskip

{\sc Proof of }(a). Now let $S$ be arbitrary. If
$\mu_{CZ}(\gamma)$ is odd, we can find $\theta\in \R \setminus
2\pi \Z$ such that
\[
2 \left\lfloor \frac{\theta}{2\pi} \right\rfloor + 1 =
\mu_{CZ}(\gamma).
\]
Reordering the coordinates $(q_1,\dots,q_n,p_1,\dots,p_n)$ as
$(q_1,p_1,\dots,q_n,p_n)$, we consider the symmetric matrix
\[
S_0 = \left( \begin{array}{cc} \theta & 0 \\ 0 & \theta
\end{array} \right) \oplus \left( \begin{array}{cc} 0 & 1 \\ 1 & 0
\end{array} \right) \oplus \dots \oplus \left( \begin{array}{cc} 0 & 1 \\ 1 & 0
\end{array} \right).
\]
By Claim 1 and Claim 2 we have
\[
\ind D_{S_0,\Lambda}^+ = - \mu_{CZ} ([0,1]\ni t \mapsto e^{tJ_0
S_0} ) = -2 \left\lfloor \frac{\theta}{2\pi} \right\rfloor - 1 = -
\mu_{CZ} (\gamma).
\]
If $\mu_{CZ} (\gamma)$ is even and $n\geq 2$, we can find
$\theta_1$, $\theta_2\in \R \setminus 2\pi \Z$ such that
\[
2 \left\lfloor \frac{\theta_1}{2\pi} \right\rfloor + 1 = 1 , \quad
2 \left\lfloor \frac{\theta_2}{2\pi} \right\rfloor + 1 =
\mu_{CZ}(\gamma)-1,
\]
and setting
\[
S_0 = \left( \begin{array}{cc} \theta_1 & 0 \\ 0 & \theta_1
\end{array} \right) \oplus \left( \begin{array}{cc} \theta_2 & 0 \\ 0 &
\theta_2 \end{array} \right) \oplus \left( \begin{array}{cc} 0 & 1
\\ 1 & 0 \end{array} \right) \oplus \dots \oplus \left( \begin{array}{cc}
0 & 1 \\ 1 & 0 \end{array} \right),
\]
we obtain again
\[
\ind D_{S_0,\Lambda}^+ = - \mu_{CZ} ([0,1]\ni t \mapsto e^{tJ_0
S_0} ) = -2 \left\lfloor \frac{\theta_1}{2\pi} \right\rfloor - 1
-2 \left\lfloor \frac{\theta_2}{2\pi} \right\rfloor - 1= -
\mu_{CZ} (\gamma).
\]
Since the Conley-Zehnder index labels the connected components of
the set (\ref{nd}), it is easy to construct a continuous homotopy
\[
H_r : [0,+\infty[ \times \T \rightarrow \mathfrak{gl} (2n,\R),
\quad r\in [0,1],
\]
such that $H_0=S_0$, $H_1=S$, $H_r(+\infty,t)$ is symmetric for
every $t\in \T$, and the solution $\gamma_r$ of
\[
\gamma_r^{\prime}(t) = J_0 H_r (+\infty,t) \gamma_r(t), \quad
\gamma_r(0)=I,
\]
satisfies (\ref{unno}). Then $r\mapsto D_{H_r,\Lambda}^+$ is a
continuous path of Fredholm operators, hence
\[
\ind D_{S,\Lambda}^+ = \ind D_{H_1,\Lambda}^+ = \ind
D_{H_0,\Lambda}^+ = \ind D_{S_0,\Lambda}^+ = - \mu_{CZ} (\gamma).
\]
This proves the index formula (\ref{AA}) in the case $\mu_{CZ}(\gamma)$ odd
or $n\geq 2$. The analysis is not complete in the case $n=1$, but
this case follows from the case $n=2$ by considering $S\oplus S$.

\medskip

{\sc Claim 3.} If $n=1$ and $S(s,t)=\theta I$, with $\theta\in \R
\setminus \pi \Z$, then
\begin{equation}
\label{22}
\ind D^+_{\theta I,\Omega} = - \left\lfloor \frac{\theta}{\pi}
\right\rfloor = \frac{1}{2} - \mu(\gamma \lambda_0,\lambda_0).
\end{equation}

\medskip

Notice that in this case $\gamma(t)=e^{\theta t J_0}$,
(\ref{ddue}) is equivalent to $\theta\notin \pi \Z$, and by
(\ref{rmi})
$\mu(\gamma \lambda_0,\lambda_0) = \lfloor \theta/\pi \rfloor + 1/2$.
The solutions of
\[
\begin{cases} \partial_s u - J_0 \partial_t u - \theta u = 0 &
\mbox{on } ]0,+\infty[ \times ]0,1[, \\ u_1(s,0)=u_1(s,1)=0 &
\forall s\geq 0, \end{cases}
\]
have the form:
\[
u(s,t) = \sum_{h\in \Z} e^{(\theta - h \pi)s} e^{h\pi t J_0}
\binom{0}{\eta_h}, \quad \eta_h\in \R.
\]
In order for such a function to decay for $s\rightarrow +\infty$
it is necessary that $\eta_h=0$ when $\theta-h \pi>0$, so in the
sum above $h$ ranges from $\lceil \theta/\pi \rceil$ to $+\infty$.
In particular,
\begin{equation}
\label{espr}
 u_1(0,t) = \sum_{h= \lceil \theta/\pi
\rceil}^{+\infty} \eta_h \sin h \pi t.
\end{equation}
Recalling that $\{\sin \pi t, \sin 2\pi t, \dots \}$ is a complete
orthogonal family in $L^2(]0,1[)$, we find that:
\begin{itemize}
\item if $\theta>0$, $u_1(0,\cdot)$ vanishes identically on
$[0,1]$ if and only if $\eta_h=0$ for every $h$; \item if
$\theta<0$, (\ref{espr}) can be rewritten as
\[
u_1(0,t) = \sum_{h=1}^{- \lceil \theta/\pi \rceil} (\eta_h -
\eta_{-h}) \sin h\pi t + \sum_{h= - \lceil \theta/\pi \rceil
+1}^{+\infty} \eta_h \sin h\pi t,
\]
so $u_1(0,\cdot)$ vanishes identically on $[0,1]$ if and only if
$\eta_h=0$ for every $h\geq - \lceil \theta/\pi \rceil +1$ and
$\eta_{-h}=\eta_h$ for every $h\in \{ 1, \dots, - \lceil
\theta/\pi \rceil \}$.
\end{itemize}

We conclude that the kernel of $D^+_{\theta I,\Omega}$ is $(0)$ is
$\theta>0$, and it consists of the functions
\[
u(s,t) = e^{\theta s} \binom{0}{\eta_h} + \sum_{h=1}^{- \lceil
\theta/\pi \rceil} \left( e^{(\theta-h \pi)s} e^{h\pi t J_0} +
e^{(\theta + h \pi)} e^{-h\pi t J_0} \right) \binom{0}{\eta_h}
\]
when $\theta<0$. Therefore
\begin{equation}
\label{aaa} \dim \ker D_{\theta I, \Omega}^+ = \begin{cases} 0 &
\mbox{if } \theta >0, \\ 1 - \left\lceil \frac{\theta}{\pi}
\right\rceil & \mbox{if } \theta<0. \end{cases}
\end{equation}
By (iii-b), the annihilator of the range of $D_{\theta I,\Omega}^+$
consists of the smooth $\R^{2n}$-valued maps on $[0,+\infty[ \times
[0,1]$ which solve
\begin{equation}
\label{croc}
\begin{cases} \partial_s v + J_0 \partial_t v + \theta v = 0 &
\mbox{on } ]0,+\infty[\times ]0,1[, \\ v_1(s,0)=v_1(s,1)=0 &
\forall s\geq 0, \\ v_2(0,t)=0 & \forall t\in [0,1], \\ v(s,\cdot)
\rightarrow 0 & \mbox{for } s\rightarrow +\infty.
\end{cases}
\end{equation}
The solutions of the first two equations of (\ref{croc}) have the
form
\[
v(s,t) = \sum_{h\in \Z} e^{(h\pi - \theta)s} e^{h \pi t J_0}
\binom{0}{\eta_h}, \quad \eta_h \in \R.
\]
In order for such a function to decay for $s\rightarrow +\infty$,
it is necessary that $\eta_h=0$ whenever $h\pi-\theta>0$, so in
the above sum $h$ ranges from $-\infty$ to $\lfloor \theta/\pi
\rfloor$. In particular, the second component of $v(0,t)$ is
\begin{equation}
\label{forma}
 v_2(0,t) = \sum_{h\leq \lfloor \theta/\pi \rfloor}
\eta_h \cos h \pi t.
\end{equation}
Recalling that $\{1,\cos \pi t, \cos 2\pi t, \dots\}$ is a
complete orthogonal family in $L^2(]0,1[)$, we find that:
\begin{itemize}
\item if $\theta<0$, $v_1(0,\cdot)$ vanishes identically on
$[0,1]$ if and only if $\eta_h=0$ for every $h$; \item if
$\theta>0$, (\ref{forma}) can be rewritten as
\[
v_2(0,t) = \eta_0 + \sum_{h=1}^{\lfloor  h/\pi \rfloor} (\eta_h +
\eta_{-h}) \cos h\pi t + \sum_{h< - \lfloor \theta/ \pi \rfloor}
\eta_h \cos h\pi t ,
\]
so $v_2(0,\cdot)$ vanishes identically on $[0,1]$ if and only if
$\eta_0=0$, $\eta_h=0$ for every $h< -\lfloor \theta/\pi \rfloor$,
and $\eta_{-h} = -\eta_h$ for every $h\in \{1,\dots,\lfloor
\theta/ \pi\rfloor\}$.
\end{itemize}

We conclude that the space of solutions of (\ref{croc}) is $(0)$
if $\theta<0$, it consists of the functions
\[
v(s,t) = \sum_{h=1}^{\lfloor \theta/ \pi\rfloor} \left( e^{h\pi
-\theta)s} e^{h\pi t J_0} - e^{-(h\pi +\theta)s} e^{-h\pi t J_0}
\right) \binom{0}{\eta_h},
\]
when $\theta<0$. Therefore
\begin{equation}
\label{bbb} \dim \coker D_{\theta I,\Omega}^+ = \begin{cases} 0 &
\mbox{if } \theta<0, \\ \left\lfloor \frac{\theta}{\pi}
\right\rfloor & \mbox{if } \theta>0. \end{cases}
\end{equation}
The index formula (\ref{22}) follows from (\ref{aaa}) and
(\ref{bbb}).

\medskip

{\sc Proof of }(b). Since $\gamma(0)=I$ and $\gamma(1) \lambda_0 \cap
\lambda_0 = (0)$, $\mu(\gamma \lambda_0,\lambda_0) - n/2$ is an
integer (see \cite{rs93} Corollary 4.12), so we can find numbers
$\theta_1,\dots\theta_n \in \R \setminus \pi \Z$ such that
\[
\sum_{h=1}^n \left\lfloor \frac{\theta_j}{\pi} \right\rfloor +
\frac{n}{2} = \mu(\gamma \lambda_0,\lambda_0).
\]
If $S_0$ is the symmetric matrix
\[
S_0 = \left( \begin{array}{cc} \theta_1 & 0 \\ 0 & \theta_1
\end{array} \right) \oplus \dots \oplus \left( \begin{array}{cc}
  \theta_2 & 0 \\ 0 & \theta_2 \end{array} \right),
\]
and $\gamma_0:[0,1]\rightarrow Sp(2n)$, $\gamma_0(t)=e^{tJ_0S_0}$,
by Claim 3 we have
\[
\ind D^+_{S_0,\Omega} = \frac{n}{2} - \mu(\gamma_0
\lambda_0,\lambda_0) =
- \sum_{j=1}^h \left\lfloor
\frac{\theta_j}{\pi} \right\rfloor = \frac{n}{2} - \mu(\gamma
\lambda_0,\lambda_0).
\]
By Corollary 4.11 of \cite{rs93}, two paths
$\gamma_0,\gamma_1:[0,1]\rightarrow Sp(2n)$ with $\gamma_j(0)=I$ and
$\gamma_j(1) \lambda_0 \cap \lambda_0 = (0)$ are
homotopic within this class if and only if
\[
\mu(\gamma_0 \lambda_0,\lambda_0) = \mu(\gamma_1 \lambda_0,\lambda_0).
\]
Therefore a homotopy argument analogous to the one used to prove (a) allows to
conclude the proof of (\ref{BB}).
\end{proof}

\subsection{Compatible orientations}
\label{orient}

The aim of this section is to orient the manifolds
$\mathcal{M}^+(q,x)$ in a way which is compatible with the
orientations of $\mathcal{M}(x,y)$ and of $W^u(q)$, for every
$x,y\in \mathcal{P}(H)$ and $q\in \mathcal{P}(L)$. The
construction will be analogous to the one described in section
\ref{scoh}.

Fix some $r\in ]2,4]$, and denote by $\Sigma^+_{\Lambda}$ the set
of operators
\[
D_{S,\Lambda}^+ : W^{1,r}_{\lambda_0} (]0,+\infty[ \times
\T,\R^{2n}) \rightarrow L^r (]0,+\infty[ \times \T,\R^{2n})
\]
of the form
\begin{equation}
\label{co11} v \mapsto \partial_s v - J_0 \partial_t v - S(s,t) v,
\end{equation}
where
\[
S\in C^0(]0,+\infty] \times \T, \mathfrak{gl}(2n)) \cap L^r(]0,1[
\times \T, \mathfrak{gl}(2n))
\]
is such that the loop $S(+\infty,\cdot) : \T \rightarrow Sym(2n)$
is a non-degenerate path, in the sense of section \ref{scoh}. By
Theorem \ref{fftt} (a), $D_{S,\Lambda}^+$ is a Fredholm operator
of index $-\mu_{CZ}(\gamma_{S(+\infty,\cdot)})$.

Similarly, $\Sigma^+_{\Omega}$ will denote the set of operators \[
D_{S,\Omega}^+ : W^{1,r}_{\lambda_0} (]0,+\infty[ \times
]0,1[,\R^{2n}) \rightarrow L^r (]0,+\infty[ \times ]0,1[,\R^{2n})
\]
of the form (\ref{co11}), where
\[
S\in C^0(]0,+\infty] \times [0,1], \mathfrak{gl}(2n)) \cap
L^r(]0,1[ \times ]0,1[, \mathfrak{gl}(2n))
\]
is such that the path $S(+\infty,\cdot) : [0,1] \rightarrow
Sym(2n)$ is non-degenerate. By Theorem \ref{fftt} (b),
$D_{S,\Omega}^+$ is a Fredholm operator of index $n/2 -
\mu(\gamma_{S(+\infty,\cdot)} (\cdot) \lambda_0,\lambda_0)$.

Therefore $\Sigma^+ \subset \mathrm{Fred} (W^{1,r}_{\lambda_0},
L^r)$ inherits the norm topology and it is the base space of the
restriction of the determinant bundle. As before (and essentially
for the same reasons), this line bundle is non-trivial on some
connected component.

If $S^+\in C^0(\T,Sym(2n))$ (resp.\ $C^0([0,1],Sym(2n))$) is a
non-degenerate path, we can consider the subset of $\Sigma^+$,
\[
\Sigma^+(S^+) := \set{D_S^+ \in \Sigma^+}{S(+\infty,\cdot) =
S^+(\cdot)}.
\]
The space $\Sigma^+(S^+)$ is contractible (it is actually
star-shaped), so the restriction of the determinant bundle to it -
denote it by $\Det(\Sigma^+(S^+))$ - is trivial.

Two orientations $o(S_1)$ of $\Det(\Sigma^+(S_1))$ and
$o(S_1,S_2)$ of $\Det(\Sigma(S_1,S_2))$ induce in a canonical way
an orientation
\[
o(S_1) \#\, o(S_1,S_2)
\]
of $\Det(\Sigma^+(S_2))$ (exactly as in section 3 of \cite{fh93}).
By construction,
\begin{equation}
\label{co12} o(S^+) \#\, o(S^+,S^+) = o(S^+).
\end{equation}
Associativity now reads as
\begin{equation}
\label{co13} (o(S_1) \#\, o(S_1,S_2)) \#\, o(S_2,S_3) = o(S_1) \#\,
(o(S_1,S_2)) \#\, o(S_2,S_3)).
\end{equation}

Let us fix a coherent orientation for $\Sigma$. A {\em compatible
orientation for $\Sigma^+$} consists of a set of orientations
$o(S^+)$ of $\Det(\Sigma^+(S^+))$, for every non-degenerate path
$S^+$, such that
\begin{equation}
\label{co14} o(S_1) \#\, o(S_1,S_2) = o(S_2),
\end{equation}
for every pair $(S_1,S_2)$ of non-degenerate paths. A compatible
orientation for $\Sigma^+$ can be constructed simply by choosing
an arbitrary non-degenerate path $S_0$, by fixing an arbitrary
orientation $o(S_0)$ of $\Det(\Sigma^+(S_0))$, and by setting
\[
o(S^+) := o(S_0) \#\, o(S_0,S^+),
\]
for every non-degenerate path $S^+$. The identity (\ref{co12})
implies that this is well-defined. The compatibility condition
(\ref{co14}) follows from the associativity property (\ref{co13})
and from the coherence of the orientation for $\Sigma$, i.e.\
(\ref{co4}). Actually, the above argument shows that there are
exactly two orientations for $\Sigma^+$ which are compatible with
a given coherent orientation for $\Sigma$.

For every $x\in \mathcal{P}(H)$ let $\Phi_x$ be the unitary
trivializations of $x^*(TT^*M)$ mapping $\lambda_0$ into the
vertical bundle, as chosen in section \ref{scoh}. Let $S_x$ be the
corresponding non-degenerate path. Let us fix a coherent
orientation for $\Sigma$, and a compatible orientation for
$\Sigma^+$. Let us fix orientations of the unstable manifolds
$W^u(q)$, for every $q\in \mathcal{P}(L)$. These data will now
determine an orientation of
\[
\Det (D_f \partial_{J,H}^+ (u)),
\]
the determinant of the fiberwise derivative of the section
\[
\partial_{J,H}^+ : \mathcal{B}^+ (q,x) \rightarrow
\mathcal{W}^+(q,x)
\] 
at every $u\in \mathcal{M}^+(q,x)$, for every
$q\in \mathcal{P}(L)$ and $x\in \mathcal{P}(H)$.

Namely, let $q\in \mathcal{P}(L)$, $x\in \mathcal{P}(H)$, and
$u\in \mathcal{M}^+(q,x)$. The closed finite codimensional
subspace
\[
W_u = \set{w\in T_u \mathcal{B}^+(q,x)}{D
\tau^*(u(0,\cdot))[w(0,\cdot)]=0}
\]
is the kernel of the continuous linear surjective operator
\[
T_u \mathcal{B}^+(q,x) \rightarrow T_{\tau^* \circ u(0,\cdot)}
W^u(q), \quad w \mapsto D \tau^*(u(0,\cdot))[w(0,\cdot)],
\]
so the orientation of $W^u(q)$ induces an orientation of the
quotient $T_u \mathcal{B}^+(q,x)/ W_u$, that is an
orientation of the line
\[
\Lambda^{\max} (T_u \mathcal{B}^+(q,x)/ W_u).
\]
By Lemma \ref{trivbis}, we can find a unitary trivialization
$\Phi_u$ of $u^*(TT^*M)$ which agrees with $\Phi_x$ for
$s=+\infty$, and maps $\lambda_0$ into the vertical subbundle. As
we have seen in section \ref{fts}, $\Phi_u$ conjugates the
restriction $D_f \partial_{J,H}^+(u)|_{W_u}$ to an operator
$D_S^+$ belonging to $\Sigma^+(S_x)$. Therefore,
\[
\Det (D_f \partial_{J,H}^+ (u)|_{W_u}) \cong \Det (D_S^+)
\]
inherits an orientation from $o(S_x)$. The analogous of Lemma 13
in \cite{fh93} implies that such an orientation does not depend on
the choice of the trivialization $\Phi_u$. From the canonical
isomorphism (see (\ref{co0}))
\[
\Det (D_f \partial_{J,H}^+ (u)) \cong \Det (D_f \partial_{J,H}^+
(u)|_{W_u}) \otimes \Lambda^{\max} (T_u \mathcal{B}^+(q,x) / W_u),
\]
we get the required orientation of $\Det (D_f \partial_{J,H}^+
(u))$.

Therefore, when the section $\partial_{J,H}^+: \mathcal{B}^+(q,x)
\rightarrow \mathcal{W}^+(q,x)$ is transverse to the zero section,
we obtain an orientation of $\mathcal{M}^+(q,x)$. In particular,
when $m(q)=\mu(x)$, the zero-dimensional manifold
$\mathcal{M}^+(q,x)$ is oriented, meaning that each point $u\in
\mathcal{M}^+(q,x)$ is given a number $\epsilon(u)\in \{-1,+1\}$.

\subsection{Compactness and convergence to broken trajectories}

The following result is now an easy consequence of the $L^{\infty}$
estimates of Theorem \ref{c0est}.

\begin{thm}
\label{newcomp}
Assume that $\|J-\widehat{J}\|_{\infty} < j_1$, where
$j_1$ is given by Theorem \ref{c0est}.
For every $q\in \mathcal{P}(L)$ and $x\in \mathcal{P}(H)$, the
space $\mathcal{M}^+_{\Lambda}(q,x)$ (resp.\
$\mathcal{M}^+_{\Omega}(q,x)$) is pre-compact in
$C^{\infty}_{\mathrm{loc}} ([0,+\infty[ \times \T,T^* M)$
(resp.\ $C^{\infty}_{\mathrm{loc}}
([0,+\infty[ \times [0,1],T^* M)$).
\end{thm}

\begin{proof} By (\ref{aces}),
\[
\mathcal{A}(x)  \leq \mathcal{A}(u(s,\cdot)) \leq \mathcal{E}(q)
\quad \forall u \in \mathcal{M}^+(q,x).
\]
Since $\tau^* u(0,\cdot)$ is an element in the unstable manifold of
$q$, which is pre-compact in the $W^{1,2}$ topology because of
Proposition \ref{mc}, the $W^{1,2}$ norm of $\tau^* u(0,\cdot)$
is uniformly bounded. The fact that $r\leq 4$ implies
that $W^{1,2}(]0,1[)$ continuously embeds into $W^{1-1/r,r}(]0,1[)$,
so also the $W^{1-1/r,r}$ norm of $\tau^* u(0,\cdot)$ is uniformly bounded.
By (H1), (H2), and the bound on $\|J-\widehat{J}\|_{\infty}$,
statements (iii) and (iv) of Theorem
\ref{c0est} imply that $\mathcal{M}^+(q,x)$
is bounded in $L^{\infty}$.

The fact that $\omega=d\theta$ implies that for every $t\in
[0,1]$, there are no non-constant
$J_t$-holomorphic spheres in $T^* M$, and no non
constant $J_t$-holomorphic discs having boundary on
some fiber of $T^* M$. Therefore a standard bubbling off argument
implies that $\mathcal{M}^+(q,x)$ is bounded in $C^1$, and then an
elliptic bootstrap produces bounds for the derivatives of every
order (see for instance \cite{flo88d} or \cite{sal90} for more
details). \end{proof}

The above result has the following consequence, which can be
proved by standard methods from Floer theory (see e.g.
\cite{sch93,sch95,sal99}).

\begin{prop}
\label{cubt}
Assume that $\|J-\widehat{J}\|_{\infty} < j_1$.
Let $(u_h)_{h\in \N}$ be a sequence in
$\mathcal{M}^+(q,x)$ and set $\bar{q}_h := u_h (0,\cdot)$.
Then there exist $q^0=q,q^1,\dots,q^a
\in \mathcal{P}(L)$, $x^0,x^1,\dots,x^b=x\in
\mathcal{P}(H)$, with $a,b\in \N$ and
\[
\mathcal{E}(q^0) > \mathcal{E}(q^1) > \dots > \mathcal{E}(q^a)
\geq \mathcal{A}(x^0) > \mathcal{A}(x^1) > \dots >
\mathcal{A}(x^b),
\]
curves
\[
\bar{q}^1 \in W^u(q^0) \cap W^s(q^1), \dots,\; \bar{q}^a \in
W^u(q^{a-1}) \cap W^s(q^a),
\]
and maps
\[
u^0 \in \mathcal{M}^+(q^a,x^0), \; u^1 \in \mathcal{M}(x^0,x^1),
\dots, \; u^b\in \mathcal{M}(x^{b-1},x^b),
\]
such that a subsequence $(\bar{q}_{h_k}, u_{h_k})_{k\in \N}$ converges to
$(\bar{q}^1,\dots,\bar{q}^a; u^0,\dots,u^b)$ in the following
sense: there are sequences $(t_k^j)_{k\in \N}\subset ]-\infty,0]$,
$j\in \{1,\dots,a\}$, and $(s_k^j)_{k\in \N}\subset
[0,+\infty[$, $j\in \{1,\dots,b\}$, such that
\[
\phi_{-\nabla_g \mathcal{E}}^{t^1_k}(\bar{q}_{h_k}) \rightarrow
\bar{q}^1, \dots, \; \phi_{-\nabla_g
\mathcal{E}}^{t^a_k}(\bar{q}_{h_k}) \rightarrow \bar{q}^a,
\quad \mbox{in }
\Lambda^1(M) \mbox{ (resp.\ in } \Omega^1(M,q_0,q_1) \mbox{)},
\]
and
\[
u_{h_k} \rightarrow u^0, \; u_{h_k}(\cdot + s^1_k,\cdot)
\rightarrow u^1, \dots, \; u_{h_k}(\cdot + s^b_k,\cdot)
\rightarrow u^b \quad \mbox{in } C^{\infty}_{\mathrm{loc}}.
\]
\end{prop}

\subsection{Transversality and gluing}

Transversality holds automatically at the stationary
solutions. Indeed, we have the following:

\begin{prop}
\label{stat}
Assume that $q\in \mathcal{P}(L)$ and $x\in \mathcal{P}(H)$
correspond to the same solution by the Legendre transform, meaning
that $\mathfrak{L}_L(t,q(t),\dot{q}(t))= (t,x(t))$ for every $t$.
Then
\[
D_f\partial_{J,H}^+ (x) : T_x \mathcal{B}^+(q,x) \rightarrow
\mathcal{W}^+  (q,x),
\]
the fiberwise derivative of $\partial_{J,H}^+$ at the
stationary solution $x$, is invertible.
\end{prop}

\begin{proof} We already know from Theorem \ref{fredth} that
$D_f\partial_{J,H}^+(x)$ is Fredholm of index 0, so it is
enough to show that its kernel is $(0)$.
If $\zeta_1$ and $\zeta_2$ are sections of $x^*(TT^*M) \rightarrow
[0,1]$, set
\[
\langle\langle \zeta_1,\zeta_2 \rangle\rangle_{L^2_J} = \int_0^1
\langle \zeta_1(t), \zeta_2(t) \rangle_{J_t}\, dt.
\]
Since $x$ is a critical point of $\mathcal{A}$, we have
\begin{equation}
\label{luno} d^2 \mathcal{A}(x) [\zeta_1,\zeta_2] = \langle\langle
\nabla_J^2 \mathcal{A}(x) \zeta_1,\zeta_2 \rangle\rangle_{L^2_J},
\end{equation}
where the operator
\[
\nabla^2_J \mathcal{A}(x) \zeta = -J(t,x) (\nabla_t \zeta - \nabla
X_H(t,x) \zeta)
\]
is $\langle\langle \cdot,\cdot \rangle\rangle_{L^2_J}$-symmetric.
Here $\nabla$ denotes the $t$-dependent Levi-Civita covariant
derivation corresponding to the $t$-dependent metric $\langle
\cdot,\cdot \rangle_{J_t}$. Moreover,
\[
D_f\partial_{J,H}^+ (x) v = \nabla_s v + \nabla^2_J
\mathcal{A}(x) v.
\]
Assume that $v\in \ker D\partial_{J,H}^+(x)$: $v$ is a
$W^{1,r}$ section of $x^*(TT^*M) \rightarrow [0,+\infty[ \times
\T$ (resp.\ $x^*(TT^*M) \rightarrow [0,+\infty[ \times [0,1]$ with
$v(s,0)\in T_{x(0)}^v T^*M$, $v(s,1)\in T_{x(1)}^v T^*M$
for every $s\geq 0$) such that
\begin{equation}
\label{lum}
\xi := D\tau^*(x(\cdot)) [v(0,\cdot)] \in T_q W^u(q),
\end{equation}
and
\begin{equation}
\label{ldue} \nabla_s v + \nabla_J^2 \mathcal{A}(x) v =0.
\end{equation}
By (\ref{lum}),
\begin{equation}
\label{ltre} d^2 \mathcal{E}(q) [\xi,\xi] \leq 0.
\end{equation}
Consider the function $\varphi: [0,+\infty[ \rightarrow \R$,
$\varphi(s) = \|v(s,\cdot)\|^2_{L^2_J}$. By (\ref{ldue}),
\begin{eqnarray}
\label{lquat} \varphi^{\prime}(s) = -2 \langle\langle v(s,\cdot),
\nabla^2_J \mathcal{A}(x)v(s,\cdot) \rangle\rangle_{L^2_J},
\\ \label{lc} \varphi^{\prime\prime}(s) = 4 \| \nabla_J^2
\mathcal{A}(x) v(s,\cdot) \|^2_{L^2_J} \geq 0.
\end{eqnarray}
Therefore $\varphi$ is convex. Since $\varphi(s)$ converges to 0
for $s\rightarrow +\infty$, either $\varphi$ is identically zero,
or $\varphi^{\prime}(0)$ is strictly negative. If by contradiction
the second case holds, (\ref{luno}) and (\ref{lquat}) imply that
\[
d^2 \mathcal{A}(x) [v(0,\cdot),v(0,\cdot)] >0.
\]
On the other hand, by (\ref{ltre}) and Lemma \ref{difdis} we obtain
\[
d^2 \mathcal{A}(x) [v(0,\cdot),v(0,\cdot)] \leq d^2 \mathcal{E}
(q) [\xi,\xi] \leq 0,
\]
a contradiction which proves that $v=0$.
\end{proof}

We shall denote by $\mathcal{J}_{\mathrm{reg}}(H,g)$ the set
of all almost complex structures $J\in
\mathcal{J}$ such that for every $q\in \mathcal{P}(L)$
and every $x\in \mathcal{P}(H)$, the section
\[
\partial_{J,H}^+ : \mathcal{B}^+(q,x) \rightarrow
\mathcal{W}^+(q,x)
\]
is transverse to the zero section.

\begin{thm}
\label{ttrans} The set
$\mathcal{J}_{\mathrm{reg}}(H,g)$ is residual in
$\mathcal{J}(H)$.
\end{thm}

Indeed, by Proposition \ref{stat} transversality is automatic when $q$
and $x$ correspond to the same solution. If this is not the case
and $u\in \mathcal{M}^+(q,x)$, $\tau^* u(0,\cdot)$ cannot be equal
to the projection $\tau^* \circ x$ (because $W^u(q) \cap
\crit(\mathcal{E}) = \{q\}$), so $u$ is not a stationary solution.
Therefore Theorem \ref{ttrans} can be proved by a standard
argument using the Sard-Smale theorem and the Carleman similarity
principle (see \cite{fhs96}).

The following gluing result can also be proved by standard
methods. The proof is completely analogous to the Floer gluing
theorem such as proven in \cite{sal99} or the
one proven in \cite{sch95}.

\begin{prop} \label{glue}
Assume that $J\in \mathcal{J}_{\mathrm{reg}}(H,g)$.

(a) Let $q^0,q^1\in \mathcal{P}(L)$, $x\in \mathcal{P}(H)$, with
$\mu(x)=m(q^1)=m(q^0)-1$, let $\bar{q}^1 \in W^u(q^0)\cap
W^s(q^1)$, and $u^0\in \mathcal{M}^+(q^1,x)$. Then there exists a
smooth curve $[0,1[ \rightarrow \mathcal{M}^+ (q^0,x)$, $r\mapsto
u(r)$, unique up to reparameterization and up to choice of its
value at $r=0$, which converges to $(\bar{q}^1;u^0)$ in the sense
of Proposition \ref{cubt} for $r\rightarrow 1$. Such a curve is
orientation preserving - with respect to the orientation of
$\mathcal{M}^+(q^0,x)$ defined in section \ref{orient} - if and
only if $\epsilon([\bar{q}^1]) \epsilon(u^0) = 1$.

(b) Let $q\in \mathcal{P}(L)$, $x^0,x^1\in \mathcal{P}(H)$, with
$m(q)=\mu(x^0) = \mu(x^1)+1$, let $u^0\in \mathcal{M}^+(q,x^0$),
and $u^1\in \mathcal{M}(x^0,x^1)$. Then there exists a smooth
curve $[0,1[ \rightarrow \mathcal{M}^+ (q,x^1)$, $r\mapsto u(r)$,
unique up to reparameterization and up to choice of its value at
$r=0$, which converges to $(\emptyset;u^0,u^1)$ in the sense of
Proposition \ref{cubt} for $r\rightarrow 1$. Such a curve is
orientation preserving - with respect to the orientation of
$\mathcal{M}^+(q,x^1)$ defined in section \ref{orient} - if and
only if $\epsilon(u^0) \epsilon([u^1]) = 1$.
\end{prop}

\subsection{The isomorphism}

Let us prove Theorem \ref{prec}. Since the isomorphism class of the
Floer complex $\{CF_*(H),\partial_*(H,J)\}$ does not depend on $J\in
\mathcal{J}_{\mathrm{reg}} (H)$
(Theorem \ref{chj}), by Theorem \ref{ttrans} we may assume that $J$
also belongs to $\mathcal{J}_{\mathrm{reg}} (H,g)$ and that
$\|J-\widehat{J}\|_{\infty} < j_1$, where $j_1$ is given by Theorem
\ref{c0est}.

Let $q\in \mathcal{P}(L)$ and $x\in \mathcal{P}(H)$ with
$m(q)=\mu(x)$. Then the zero-dimensional manifold
$\mathcal{M}^+(q,x)$ is compact: otherwise we could deduce a
violation of transversality from Proposition \ref{cubt}.
Therefore, $\mathcal{M}^+(q,x)$ is a finite set, and we can
indicate by $n^+(q,x)$ the integer
\[
n^+(q,x) = \sum_{u \in \mathcal{M}^+(q,x)} \epsilon(u),
\]
the numbers $\epsilon(u)$ having being defined in section
\ref{orient}. The sequence of homomorphisms
\[
\Theta_k = \Theta_k(H,J,g) : CM_k(\mathcal{E}) \rightarrow
CF_k(H), \quad k\in \N,
\]
can be defined in terms of the generators as
\[
\Theta_k q = \sum_{\substack{x\in \mathcal{P}(H)\\ \mu(x)=k}}
n^+(q,x) x , \quad \mbox{for } q\in \mathcal{P}(L), \; m(q)=k.
\]
A standard argument using Propositions \ref{cubt} and \ref{glue}
implies that $\Theta$ is a chain homomorphism, meaning that
\[
\Theta_{k-1} \partial_k (\mathcal{E},g) = \partial_k (H,J)
\Theta_k \quad \forall k\geq 1.
\]
Assume that $\mathcal{E}(q)\leq \mathcal{A}(x)$. Then (\ref{aces})
implies that $\mathcal{M}^+(q,x)$ is empty - hence $n^+(q,x)=0$ -
unless $q$ and $x$ correspond to the same solution by the Legendre
transform, in which case $\mathcal{M}^+(q,x) = \{(q,x)\}$ - hence
$n^+(q,x)=\pm 1$.
Let us order the generators of
$CM_k(\mathcal{E})$ and $CF_k(H)$ by increasing action, choosing any
order for subsets of solutions with identical action (but keeping the
same order for the solutions of the Lagrangian system and the
corresponding solutions of the Hamiltonian system). 
Then the homomorphism $\Theta_k$ is represented by a (possibly infinite)
square matrix which is lower triangular and has the entries $\pm
1$ on the diagonal. Such a homomorphism is necessarily invertible,
hence $\Theta$ is a chain complex isomorphism.

Finally, $\mathcal{M}(q,x)$ is necessarily empty - hence
$n^+(q,x)=0$ - if $q$ and $\tau^* \circ x$ are not homotopic
within the space of free loops (resp.\ within the space of curves
joining $q_0$ to $q_1$). Therefore $\Theta$ is compatible with the
splitting of the Morse and the Floer complexes corresponding to
the partition of $\pi_1(M)$ into its conjugacy classes (resp.\
into its elements). This completes the proof of Theorem
\ref{prec}.

\providecommand{\bysame}{\leavevmode\hbox to3em{\hrulefill}\thinspace}
\providecommand{\MR}{\relax\ifhmode\unskip\space\fi MR }
\providecommand{\MRhref}[2]{%
  \href{http://www.ams.org/mathscinet-getitem?mr=#1}{#2}
}
\providecommand{\href}[2]{#2}


\begin{thebibliography}{10}

\bibitem{ama04m}
A.~Abbondandolo and P.~Majer, \emph{Lectures on the {M}orse complex for
  infinite dimensional manifolds}, Morse theoretic methods in non-linear
  analysis and symplectic topology (Montreal), Kluwert, 2004, to appear.

\bibitem{ada75}
R.~A. Adams, \emph{Sobolev spaces}, Pure and Applied Mathematics, vol.~65,
  Academic Press, New York, 1975.

\bibitem{ben86}
V.~Benci, \emph{Periodic solutions of {L}agrangian systems on a compact
  manifold}, J. Diff. Eq. \textbf{63} (1986), 135--161.

\bibitem{cie94}
K.~Cieliebak, \emph{Pseudo-holomorphic curves and periodic orbits on cotangent
  bundles}, J. Math. Pures Appl. (9) \textbf{73} (1994), 251--278.

\bibitem{coh04}
R.~Cohen, \emph{Morse theory, graphs, and string topology}, {\tt
  arXiv:math.GT/0411272} (2004).

\bibitem{cr03}
O.~Cornea and A.~Ranicki, \emph{Rigidity and glueing for {M}orse and {N}ovikov
  complexes}, J. Eur. Math. Soc. \textbf{5} (2003), 343--394.

\bibitem{dol80}
A.~Dold, \emph{Lectures on algebraic topology}, Springer, Berlin, 1980.

\bibitem{dui76}
J.~J. Duistermaat, \emph{On the {M}orse index in variational calculus},
  Advances in Math. \textbf{21} (1976), 173--195.

\bibitem{flo88d}
A.~Floer, \emph{The unregularized gradient flow of the symplectic action},
  Comm. Pure Appl. Math. \textbf{41} (1988), 775--813.

\bibitem{fh93}
A.~Floer and H.~Hofer, \emph{Coherent orientations for periodic orbit problems
  in symplectic geometry}, Math. Z. \textbf{212} (1993), 13--38.

\bibitem{fhs96}
A.~Floer, H.~Hofer, and D.~Salamon, \emph{Transversality in elliptic {M}orse
  theory for the symplectic action}, Duke Math. J. \textbf{80} (1996),
  251--292.

\bibitem{hof85b}
H.~Hofer, \emph{Lagrangian embeddings and critical point theory}, Ann. Inst. H.
  Poincar\'e Anal. Non Lin\'eaire \textbf{5} (1985), 407--462.

\bibitem{hz94}
H.~Hofer and E.~Zehnder, \emph{Symplectic invariants and {H}amiltonian
  dynamics}, Birkh\"auser Advanced Texts: Basler Lehrb\"ucher, Birkh\"auser,
  Basel, 1994.

\bibitem{kli82}
W.~Klingenberg, \emph{Riemannian geometry}, Walter de Gruyter \& Co., Berlin,
  1982.

\bibitem{man91}
R.~Ma{\~{n}}\'e, \emph{Global variational methods in conservative dynamics},
  IMPA, Rio de Janeiro, 1991.

\bibitem{maz85}
V.~G. Maz'ja, \emph{Sobolev spaces}, Springer-Verlag, New York, 1985.

\bibitem{ms04}
D.~McDuff and D.~Salamon, \emph{{$J$}-holomorphic curves and symplectic
  topology}, Colloquium Publications, vol.~52, American Mathematical Society,
  Providence, R.I., 2004.

\bibitem{mo97}
D.~Milinkovic and Y.~G. Oh, \emph{Floer homology as a stable {M}orse homology},
  J. Korean Math. Soc. \textbf{34} (1997), 1065--1087.

\bibitem{mo98}
D.~Milinkovic and Y.~G. Oh, \emph{Generating function versus the
  action functional}, Geometry,
  Topology and Dynamics, CRM Proc. Lecture Notes, vol.~15, Amer. Math. Soc.,
  Providence, RI, 1998, pp.~107--125.

\bibitem{oh97}
Y.~G. Oh, \emph{Symplectic topology as the topology of action functional. {I} -
  {R}elative {F}loer theory on the cotangent bundle}, J. Diff. Geom.
  \textbf{46} (1997), 499--577.

\bibitem{pal63}
R.~S. Palais, \emph{Morse theory on {H}ilbert manifolds}, Topology \textbf{2}
  (1963), 299--340.

\bibitem{qui85}
D.~Quillen, \emph{Determinants of {C}auchy-{R}iemann operators over a {R}iemann
  surface}, Functional Anal. Appl. \textbf{19} (1985), 31--34.

\bibitem{ram04}
A.~Ramirez, Ph.D. thesis, Stanford University, 2004, in preparation.

\bibitem{rs93}
J.~Robbin and D.~Salamon, \emph{Maslov index theory for paths}, Topology
  \textbf{32} (1993), 827--844.

\bibitem{rs95}
J.~Robbin and D.~Salamon, \emph{The spectral flow and the {M}aslov 
index}, Bull. London Math. Soc. \textbf{27} (1995), 1--33.

\bibitem{sal90}
D.~Salamon, \emph{Morse theory, the {C}onley index and {F}loer homology}, Bull.
  London Math. Soc. \textbf{22} (1990), 113--140.

\bibitem{sal99}
D.~Salamon, \emph{Lectures on {F}loer homology}, Symplectic geometry
and topology (Y.~Eliashberg and L.~Traynor, eds.), IAS/Park City 
Mathematics Series, Amer. Math. Soc., 1999, pp.~143--225.

\bibitem{sw03}
D.~Salamon and J.~Weber, \emph{Floer homology and the heat flow}, e-print {\tt
  arXiv:math.SG/0304383}, 2003.

\bibitem{sz92}
D.~Salamon and E.~Zehnder, \emph{Morse theory for periodic solutions of
  {H}amiltonian systems and the {M}aslov index}, Comm. Pure Appl. Math.
  \textbf{45} (1992), 1303--1360.

\bibitem{sch93}
M.~Schwarz, \emph{Morse homology}, Birkh\"auser, Basel, 1993.

\bibitem{sch95}
M.~Schwarz, \emph{Cohomology operations from {$S^1$}-cobordisms in {F}loer
  homology}, Ph.D. thesis, Swiss Federal Inst.~of Techn. Zurich, Zurich,
  Diss.~ETH No.~11182, 1995.

\bibitem{tra94}
L.~Traynor, \emph{Generating function homology}, Geom. and Funct. Anal.
  \textbf{4} (1994), 718--748.

\bibitem{vit95}
C.~Viterbo, \emph{Generating functions in symplectic topology and
  applications}, Proceedings ICM 94, Z\"urich (Basel), vol.~1, Birkh\"auser,
  1995, pp.~537--547.

\bibitem{vit96}
C.~Viterbo, \emph{Functors and computations in {F}loer homology with 
applications, {P}art {II}}, preprint, 1996.

\bibitem{web02}
J.~Weber, \emph{Perturbed closed geodesics are periodic orbits: index and
  transversality}, Math. Z. \textbf{241} (2002), 45--82.

\end{thebibliography}
\end{document}